\documentclass[12pt,reqno]{amsart}
\usepackage{times}

\newtheorem{thm}{Theorem}
\newtheorem{lemma}[thm]{Lemma}
\newtheorem{cor}[thm]{Corollary}
\newtheorem{prop}[thm]{Proposition}

\newcommand{\bn}{{{\mathbb B}_n}}
\newcommand{\sn}{{{\mathbb S}_n}}
\newcommand{\C}{{\mathbb C}}
\newcommand{\D}{{\mathbb D}}
\newcommand{\cn}{{\mathbb C}^n}
\newcommand{\inb}{\int_\bn}
\newcommand{\ins}{\int_\sn}
\newcommand{\apa}{A^p_\alpha}
\newcommand{\dva}{dv_\alpha}
\newcommand{\lpa}{L^p(\bn,\dva)}
\newcommand{\bloch}{{\mathcal B}}
\newcommand{\re}{{\rm Re}\,}

\begin{document}

\title[Bergman Spaces]
{Theory of Bergman Spaces\\
  in the Unit Ball of $\cn$}

\author{Ruhan Zhao and Kehe Zhu}
\thanks{The second author is partially supported by a grant from the NSF}
\address{Department of Mathematics\\
                SUNY\\
                Brockport, NY 14420, USA}
\email{rzhao@brockport.edu}
\address{Department of Mathematics\\
              SUNY\\
              Albany, NY 12222, USA}
\email{kzhu@math.albany.edu}

\begin{abstract}
There has been a great deal of work done in recent years on weighted Bergman 
spaces $\apa$ on the unit ball $\bn$ of $\cn$, where $0<p<\infty$ and $\alpha>-1$. 
We extend this study in a very natural way to the case where $\alpha$ is {\em any} 
real number and $0<p\le\infty$. This unified treatment covers all classical Bergman 
spaces, Besov spaces, Lipschitz spaces, the Bloch space, the Hardy space $H^2$, 
and the so-called Arveson space. Some of our results about integral representations, 
complex interpolation, coefficient multipliers, and Carleson measures are new 
even for the ordinary (unweighted) Bergman spaces of the unit disk. 
\end{abstract}

\maketitle

\section{Introduction}

Throughout the paper we fix a positive integer $n$ and let 
$$\cn=\C\times\cdots\times\C$$
denote the $n$ dimensional complex Euclidean space. For
$z=(z_1,\cdots,z_n)$ and $w=(w_1,\cdots,w_n)$ in $\cn$ we write
$$\langle z,w\rangle=z_1\overline w_1+\cdots+z_n\overline w_n$$
and
$$|z|=\sqrt{|z_1|^2+\cdots+|z_n|^2}.$$
The open unit ball in $\cn$ is the set
$$\bn=\{z\in\cn:|z|<1\}.$$
We use $H(\bn)$ to denote the space of all holomorphic functions in $\bn$.

For any $-\infty<\alpha<\infty$ we consider the positive measure
$$\dva(z)=(1-|z|^2)^\alpha\,dv(z),$$
where $dv$ is volume measure on $\bn$. It is easy to see that $\dva$ is
finite if and only if $\alpha>-1$. When $\alpha>-1$, we normalize $dv_\alpha$
so that it is a probability measure.

Bergman spaces with standard weights are defined as follows:
$$\apa=H(\bn)\cap L^p(\bn,\dva),$$
where $p>0$ and $\alpha>-1$. Here the assumption that $\alpha>-1$ is essential, 
because the space $L^p(\bn,\dva)$ does not contain any holomorphic function other 
than $0$ when $\alpha\le-1$. When $\alpha=0$, we use $A^p$ to denote the 
ordinary unweighted Bergman spaces. Bergman spaces with standard weights on
the unit ball have been studied by numerous authors in recent years. See 
Aleksandrov \cite{Alex}, Beatrous-Burbea \cite{BB2}, Coifman-Rochberg \cite{cr}, 
Rochberg \cite{Roch}, Rudin \cite{rudin}, Stoll \cite{Stoll}, and Zhu \cite{zhu6} for 
results and references.

In this paper we are going to extend the definition of $\apa$ to the case in which 
$\alpha$ is any real number and develop a theory for the extended family of spaces. 
More specifically, we study the following topics about the generalized spaces $\apa$: 
various characterizations, integral representations, atomic decomposition, complex 
interpolation, optimal pointwise estimates, duality, reproducing kernels when $p=2$, 
Carleson type measures, and various special cases. 

Some results in the paper are straightforward consequences or generalizations of known 
results in the case $\alpha>-1$ , thanks to the isomorphism between $\apa$ and $A^p$ via 
fractional integral and differential operators, and several results here have appeared 
before in the literature in different forms or different contexts. We included 
them here with full proofs for the sake of a complete and coherent theory. But we will
try our best to reference the appropriate early sources whenever possible. On the
other hand, most of our results require new techniques and reveal new properties.
Although our paper is expository in nature, several of our results are new even 
in the case of ordinary Bergman spaces of the unit disk.

Our starting point is the observation that, for $p>0$ and $\alpha>-1$, a holomorphic
function $f$ in $\bn$ belongs to $\apa$ if and only if the function
$(1-|z|^2)Rf(z)$ belongs to $L^p(\bn,\dva)$, where
$$Rf(z)=\sum_{k=1}^nz_k\frac{\partial f}{\partial z_k}(z)$$
is the radial derivative of $f$. This result is well known to experts in the field and
is sometimes referred to as a theorem of Hardy and Littlewood (especially in the
one-dimensional case). See Beatrous \cite{Beatrous}, Pavlovic \cite{Pav}, or 
Theorem~2.16 of \cite{zhu6}. More generally, we can repeatedly apply this result and 
show that, for any positive integer $k$, a holomorphic function $f$ is in $\apa$ if and 
only if the function $(1-|z|^2)^kR^kf(z)$ belongs to $L^p(\bn,\dva)$.

Now for $p>0$ and $-\infty<\alpha<\infty$ we fix a nonnegative integer $k$
with $pk+\alpha>-1$ and define $\apa$ as the space of holomorphic functions $f$ in 
$\bn$ such that the function $(1-|z|^2)^kR^kf(z)$ belongs to $L^p(\bn,\dva)$.
As was mentioned in the previous paragraph, this definition of $\apa$ is consistent
with the traditional definition when $\alpha>-1$. Also, it is easy to show (see
Section 3) that the definition of $\apa$ is independent of the integer $k$.

Our family of generalized Bergman spaces $\apa$, with $p>0$ and $\alpha$ real,
covers any space (except $H^\infty$) of holomorphic functions that is defined in terms of 
membership in $L^p(\bn,dv)$, $0<p\le\infty$, for any combination of partial derivatives
and powers of $1-|z|^2$. This family of spaces has appeared before
in the literature under different names. For example, for any positive $p$
and real $s$ there is the classical diagonal Besov space $B^s_p$ consisting of
holomorphic functions $f$ in $\bn$ such that $(1-|z|^2)^{k-s}R^kf(z)$ belongs to
$L^p(dv_{-1})$, where $k$ is any positive integer greater than $s$. It is clear that
$B^s_p=A^p_\alpha$ with $\alpha=-(ps+1)$; and $A^p_\alpha=B^s_p$ with
$s=-(\alpha+1)/p$. Thus our spaces $\apa$ are exactly the diagonal
Besov spaces. See Ahern-Cohen \cite{AB}, Arazy-Fisher-Janson-Peetre \cite{AFJP}, 
Arcozzi-Rochberg-Sawyer \cite{ars2}, Frazier-Jawerth \cite{FJ},
Hahn-Youssfi \cite{HY1}\cite{HY2}, Kaptanoglu \cite{Kap1}\cite{Kap2}, Nowark \cite{N}, 
Peloso \cite{Peloso}, and Zhu \cite{zhu6} for some recent results on such Besov spaces 
and more references. In particular, our spaces $\apa$ are the same as
the spaces $B^p_q$ (with $q=\alpha$) in Kaptanoglu \cite{Kap1}, although an
unnecessary condition $-qp+q>-1$ was imposed in \cite{Kap1}. 

On the other hand, if $s$ is a positive integer, $p$ is positive, and $\alpha$ is real, 
then there is the Sobolev space $W^p_{s,\alpha}$ consisting of holomorphic functions 
$f$ in $\bn$ such that the partial derivatives of $f$ of order up to $N$ all belong to 
$L^p(\bn,dv_\alpha)$. It is easy to see that our generalized Bergman spaces are exactly 
the holomorphic Sobolev spaces. See Ahern-Cohen \cite{AB}, Aleksandrov \cite{Alex},
Beatrous-Burbea \cite{BB2} for results and more references. 

Therefore, for those who are more familiar or more comfortable with Besov or
Sobolev spaces, our paper can be considered a unified theory for such spaces
as well. However, we believe that most people nowadays are familiar and comfortable 
with the term ``Bergman spaces'', and our theory here is almost identical to the
theory of ordinary Bergman spaces (as presented in Zhu \cite{zhu6} for example), so 
it is also reasonable for us to call $\apa$ weighted Bergman spaces.

\section{Preliminaries}

In this section we present preliminary material on Bergman kernel functions and
fractional differential and integral operators. This material will be heavily used
in later sections.

Throughout the paper we use
$$m=(m_1,\cdots,m_n)$$
to denote an $n$-tuple of nonnegative integers. It is customary to write
$$|m|=m_1+\cdots+m_n$$
and
$$m!=m_1!\cdots m_n!.$$
If $z=(z_1,\cdots,z_n)$ is a point in $\cn$, we write
$$z^m=z_1^{m_1}\cdots z_n^{m_n}.$$
The following multi-nomial formula will be used (implicitly) several times
later on:
\begin{equation}
\langle z,w\rangle^k=\sum_{|m|=k}\frac{k!}{m!}z^m\overline w^m.
\label{eq1}
\end{equation}

If $f$ is a holomorphic function in $\bn$, it has a unique Taylor series,
$$f(z)=\sum_ma_mz^m.$$
If we define
$$f_k(z)=\sum_{|m|=k}a_mz^m,\qquad k=0,1,2,\cdots,$$
then each $f_k$ is a homogeneous polynomial of degree $k$,
and we can rearrange the Taylor series of $f$ as follows:
$$f(z)=\sum_{k=0}^\infty f_k(z).$$
This is called the homogeneous expansion of $f$.

Using homogeneous expansion of $f$ we can write the radial derivative $Rf$ as 
$$Rf(z)=\sum_{k=1}^{\infty}kf_k(z).$$
More general, for any real number $t$, we can define
the following fractional radial derivative for a holomorphic
function $f$ in $\bn$:
$$R^tf(z)=\sum_{k=1}^{\infty}k^tf_k(z).$$

When we work with partial derivatives, we will use the following notation:
$$\partial^mf=\frac{\partial^{|m|}f}{\partial z_1^{m_1}\cdots\partial z_n^{m_n}},$$
where $m$ is any $n$-tuple of nonnegative integers.

An important tool in the study of holomorphic function spaces is the notion of
fractional differential and integral operators. There are numerous types of
fractional differential and integral operators, we introduce one that is intimately
related to and interacts well with the Bergman kernel functions. More specifically, 
for any complex parameters $s$ and $t$ with the property that neither $n+s$ nor 
$n+s+t$ is a negative integer, we define two operators $R^{s,t}$ and $R_{s,t}$ on 
$H(\bn)$ as follows. If
$$f(z)=\sum_{k=0}^\infty f_k(z)$$
is the homogeneous expansion of a holomorphic function in $\bn$, we define
$$R^{s,t}f(z)=\sum_{k=0}^\infty\frac{\Gamma(n+1+s)\Gamma(n+1+k
+s+t)}{\Gamma(n+1+s+t)\Gamma(n+1+k+s)}f_k(z).$$
If $H(\bn)$ is equipped with the topology of ``uniform convergence on compact sets",
it is easy to see that each $R^{s,t}$ is a continuous invertible operator on
$H(\bn)$. We use $R_{s,t}$ to denote the inverse of $R^{s,t}$ on
$H(\bn)$. Thus
$$R_{s,t}f(z)=\sum_{k=0}^\infty\frac{\Gamma(n+1+s+t)\Gamma(n+1+k+s)}
{\Gamma(n+1+s)\Gamma(n+1+k+s+t)}f_k(z).$$

When $s$ is real and $t>0$, it follows from Stirling's formula that
$$\frac{\Gamma(n+1+s)\Gamma(n+1+k+s+t)}{\Gamma(n+1+s+t)\Gamma(n+1+k+s)}
\sim k^t$$
as $k\to\infty$. In this case, $R^{s,t}$ is indeed a fractional radial differential operator
of order $t$ and $R_{s,t}$ is a fractional radial integral operator of order $t$.

The operators $R^{s,t}$ and $R_{s,t}$ seem to have first appeared in Peloso \cite{Peloso},
and independently in Zhu \cite{zhu3}\cite{zhu4}\cite{zhu5}, as a way to define and study 
holomorphic function spaces on the unit ball, and more generally, on bounded symmetric
domains. This type of fractional differential and integral operators also
became an important tool in the books by Arcozzi-Rochberg-Sawyer \cite{ars2} and 
Zhu \cite{zhu6}. 

Kaptanoglu \cite{Kap1}\cite{Kap2}\cite{Kap3} used these operators in a slightly more 
general way. More specifically, the technical conditions that $n+s$ and $n+s+t$ should not 
be negative integers can be removed if one is willing to make a separate definition
for $R^{s,t}$ (and $R_{s,t}$) in this case. However, since these operators are meant
to transform the kernel function $(1-\langle z,w\rangle)^{-(n+1+s)}$ to
$(1-\langle z,w\rangle)^{-(n+1+s+t)}$, it is clear that the technical conditions mentioned
above are natural. Otherwise, these functions would become polynomials
and the corresponding reproducing Hilbert spaces would become finite dimensional.
Besides, in all our applications, it always involves in choosing a sufficiently large parameter 
$s$, and with the technical conditions imposed on $s$ and $t$, there is never a lack
of $s$ for such choices. Also, the use of complex parameters does not present any 
extra difficulty and will be more convenient for us on several occasions.

\begin{lemma}
Suppose neither $n+s$ nor $n+s+t$ is a negative integer. Then $R_{s,t}=R^{s+t,-t}$.
\label{1}
\end{lemma}

\begin{proof}
This follows directly from the definition of these operators.
\end{proof}

\begin{lemma}
Suppose $s$, $t$, and $\lambda$ are complex parameters such that none of
$n+\lambda$, $n+\lambda+t$, and $n+\lambda+s+t$ is a negative
integer. Then $R^{\lambda,t}R^{\lambda+t,s}=R^{\lambda,s+t}$.
\label{2}
\end{lemma}

\begin{proof}
This also follows from the definition of these operators.
\end{proof}

As was mentioned earlier, the main advantage of the operators $R^{s,t}$ and
$R_{s,t}$ is that they interact well with Bergman kernel functions. This is made
precise by the following result.

\begin{prop}
Suppose neither $n+s$ nor $n+s+t$ is a negative integer. Then
$$R^{s,t}\frac1{(1-\langle z,w\rangle)^{n+1+s}}=\frac1{(1-\langle z,w\rangle)^{n+1+s+t}},$$
and
$$R_{s,t}\frac1{(1-\langle z,w\rangle)^{n+1+s+t}}=\frac1{(1-\langle z,w\rangle)^{n+1+s}}.$$
Furthermore, the above relations uniquely determine the operators $R^{s,t}$
and $R_{s,t}$ on $H(\bn)$.
\label{3}
\end{prop}

\begin{proof}
See Proposition 1.14 of \cite{zhu6}. The proof there is for the case when $s$ and
$t$ are real. But obviously the same proof works for complex parameters as well.
\end{proof}

Most of the time we use the above proposition as follows. If a holomorphic function
$f$ in $\bn$ has an integral representation
$$f(z)=\inb\frac{d\mu(w)}{(1-\langle z,w\rangle)^{n+1+s}},$$
then
$$R^{s,t}f(z)=\inb\frac{d\mu(w)}{(1-\langle z,w\rangle)^{n+1+s+t}}.$$
In particular, if $\alpha>-1$ and $n+\alpha+t$ is not a negative integer, then
$$R^{\alpha,t}f(z)=\lim_{r\to1^-}\inb\frac{f(rw)\,dv_\alpha(w)}{(1-\langle z,w\rangle)^{n
+1+\alpha+t}}$$
for every function $f\in H(\bn)$. See Corollary 2.3 of \cite{zhu6}.

\begin{prop}
Suppose $N$ is a positive integer and $s$ is a complex number such that
$n+s$ is not a negative integer. Then the operator $R^{s,N}$ is a linear partial
differential operator on $H(\bn)$ of order $N$ with polynomial coefficients, that is,
$$R^{s,N}f(z)=\sum_{|m|\le N}p_m(z)\partial^mf(z),$$
where each $p_m$ is a polynomial.
\label{4}
\end{prop}

\begin{proof}
The proof of Proposition 1.15 of \cite{zhu6} works for complex parameters as well.
\end{proof}

\begin{prop}
Suppose $s$ and $t$ are complex parameters such that neither $n+s$ nor
$n+s+t$ is a negative integer. If $\alpha=s+N$ for some positive integer $N$,
then
$$R^{s,t}\frac1{(1-\langle z,w\rangle)^{n+1+\alpha}}=\frac{h(\langle z,w\rangle)}
{(1-\langle z,w\rangle)^{n+1+\alpha+t}},$$
where $h$ is a certain one-variable polynomial of degree $N$. Similarly,
there exists a one-variable polynomial $q$ of degree $N$ such that
$$R_{s,t}\frac1{(1-\langle z,w\rangle)^{n+1+\alpha+t}}=\frac{q(\langle z,w\rangle)}
{(1-\langle z,w\rangle)^{n+1+\alpha}}.$$
\label{5}
\end{prop}

\begin{proof}
See the proof of Lemma 2.18 of \cite{zhu6} for the result concerning $R^{s,t}$.
Combining this with Lemma~\ref{1}, the result for $R_{s,t}$ follows as well.

Alternatively, we can use Proposition~\ref{3} to write
$$R^{s,t}\frac1{(1-\langle z,w\rangle)^{n+1+\alpha}}=R^{s,t}R^{s,N}\frac1{(1-
\langle z,w\rangle)^{n+1+s}}.$$
Since $R^{s,N}$ and $R^{s,t}$ commute, another application of Proposition~\ref{3}
gives
$$R^{s,t}\frac1{(1-\langle z,w\rangle)^{n+1+\alpha}}=R^{s,N}\frac1{(1-\langle z,
w\rangle)^{n+1+s+t}}.$$
The desired result then follows from Proposition~\ref{4}.
\end{proof}

We also include an easy but important fact concerning the radial derivative.

\begin{lemma}
For any positive integer $k$ the operator $R^k$ is a $k$th order partial differential
operator on $H(\bn)$ with polynomial coefficients.
\label{6}
\end{lemma}

\begin{proof}
Obvious.
\end{proof}

We are going to need two integral estimates involving Bergman
kernel functions.

\begin{prop}
Suppose $s$ and $t$ are real numbers with $s>-1$. Then the integral
$$I(z)=\inb\frac{(1-|w|^2)^s\,dv(w)}{|1-\langle z,w\rangle|^{n+1+s+t}}$$
has the following asymptotic behavior as $|z|\to1^-$.
\begin{enumerate}
\item[(a)] If $t<0$, then $I(z)$ is continuous on $\overline{\bn}$. In particular,
$I(z)$ is bounded for $z\in\bn$.
\item[(b)] If $t>0$, then $I(z)$ is comparable to $(1-|z|^2)^{-t}$.
\item[(c)] If $t=0$, then $I(z)$ is comparable to $-\log(1-|z|^2)$.
\end{enumerate}
\label{7}
\end{prop}

\begin{proof}
See Proposition 1.4.10 of Rudin \cite{rudin}.
\end{proof}

\begin{prop}
Suppose $a$ and $b$ are complex parameters. If $S$ and $T$ are integral
operators defined by
$$Sf(z)=(1-|z|^2)^a\inb\frac{(1-|w|^2)^bf(w)\,dv(w)}{(1-\langle z,w\rangle)^{n+1+a+b}},$$
and
$$Tf(z)=(1-|z|^2)^a\inb\frac{(1-|w|^2)^bf(w)\,dv(w)}{|1-\langle z,w\rangle|^{n+1+a+b}},$$
then for any $1\le p<\infty$ and $\alpha$ real, the following conditions are
equivalent.
\begin{enumerate}
\item[(a)] The operator $S$ is bounded on $L^p(\bn,dv_\alpha)$.
\item[(b)] The operator $T$ is bounded on $L^p(\bn,dv_\alpha)$.
\item[(c)] The parameters satisfy $-p\,\re a<\alpha+1<p(\re b+1)$.
\end{enumerate}
\label{8}
\end{prop}

\begin{proof}
See \cite{KZ} or Theorem 2.10 of \cite{zhu6}. Once again, those proofs are given
for real parameters, but the proof for the complex case is essentially the same. The 
only extra attention to pay is this: when $\lambda=u+iv$ is a complex constant, we have
$$(1-\langle z,w\rangle)^\lambda=|1-\langle z,w\rangle|^\lambda\exp(iu\theta-v\theta),$$
where $\theta$ is the argument of $1-\langle z,w\rangle$, say $\theta\in[0,2\pi)$. It
follows that
$$|(1-\langle z,w\rangle)^\lambda|=|1-\langle z,w\rangle|^u\exp(-v\theta).$$
Since $v$ is a constant and $\theta\in[0,2\pi)$, we see that
$$|(1-\langle z,w\rangle)^\lambda|\sim|1-\langle z,w\rangle|^u
=|1-\langle z,w\rangle|^{\re\lambda}.$$
\end{proof}

Note that certain special cases of the above proposition can be found in
Forelli-Rudin \cite{fr} and Rudin \cite{rudin}.

\begin{prop}
Suppose $\re\alpha>-1$. Then there exists a constant $c_\alpha$ such that
$$f(z)=c_\alpha\inb\frac{f(w)(1-|w|^2)^\alpha\,dv(w)}{(1-\langle z,w\rangle)^{n+1+
\alpha}},\qquad z\in\bn,$$
where $f$ is any holomorphic function in $\bn$ such that
$$\inb|f(z)(1-|z|^2)^\alpha|\,dv(z)<\infty.$$
\label{9}
\end{prop}

\begin{proof}
See Theorem 7.1.4 of Rudin \cite{rudin} or Theorem 2.2 of Zhu \cite{zhu6}.
\end{proof}

\section{Isomorphism of Bergman Spaces}

Our first main result shows that for fixed $p$, the spaces $\apa$ are all isomorphic.
A word of caution is necessary here: while the isomorphism among $\apa$ reduces
the topological structure of $\apa$ to that of the ordinary Bergman space $A^p$, it does
not help too much when the properties of individual functions are concerned. This
is clear in the Hilbert space case: the Hardy space $H^2$, the Bergman space $A^2$,
and the Dirichlet space $B_2$ are all isomorphic as Hilbert spaces, but their respective
function theories behave much differently from one to another.

There is a good amount of overlap between the material in this and the next section with 
the results in Beatrous-Burbea \cite{BB1}\cite{BB2}, Kaptanoglu \cite{Kap2}, and 
Peloso \cite{Peloso}. We present independent proofs here in order to 
achieve a complete and coherent theory. As was mentioned in the introduction, the 
spaces $B^p_q$ in Kaptanoglu \cite{Kap2} and our spaces $\apa$ are actually the same 
(with the identification of $\alpha$ and $q$), while the family of spaces $A^p_q$ in 
Beatrous-Burbea \cite{BB2} covers ordinary Bergman spaces
(our $A^p_\alpha$ with $\alpha>-1$) and Hardy spaces $H^p$.

\begin{thm}
Suppose $p>0$ and $\alpha$ is real. If $s$ is a complex parameter such that neither 
$n+s$ nor $n+s+(\alpha/p)$ is a negative integer, then a holomorphic function 
$f$ in $\bn$ is in $\apa$ if and only if $R_{s,\alpha/p}f$ is in $A^p$. Equivalently, 
$R_{s,\alpha/p}$ is an invertible operator from $\apa$ onto $A^p$.
\label{10}
\end{thm}

\begin{proof}
Recall that a holomorphic function $f$ in $\bn$ is in $\apa$ if and only if there 
exists a nonnegative integer $k$ with $pk+\alpha>-1$ such that the function
$(1-|z|^2)^kR^kf(z)$ is in $\lpa$. Obviously, this is equivalent to the condition
that
\begin{equation}
R^kf\in L^p(\bn,dv_{pk+\alpha}).
\label{eq2}
\end{equation}

By Theorem 2.16 of \cite{zhu6}, the condition that $R_{s,\alpha/p}f\in A^p$
is equivalent to
$$(1-|z|^2)^kR^kR_{s,\alpha/p}f(z)\in L^p(\bn,dv).$$
Since $R^k$ commutes with $R_{s,\alpha/p}$, the above condition is equivalent to
$$R_{s,\alpha/p}R^kf\in L^p(\bn,dv_{pk}).$$
If $\alpha>0$, then by Theorem 2.19 of \cite{zhu6}, the above condition is equivalent to
$$(1-|z|^2)^{\alpha/p}R^{s,\alpha/p}R_{s,\alpha/p}R^kf\in L^p(\bn,dv_{pk}).$$
Since $R^{s,\alpha/p}$ is the inverse of $R_{s,\alpha/p}$, the above condition is
equivalent to
$$(1-|z|^2)^{\alpha/p}R^kf\in L^p(\bn,dv_{pk}),$$
which is the same as (\ref{eq2}). This proves the theorem for $\alpha>0$.

If $\alpha=0$, the operator $R_{s,\alpha/p}$ becomes the identity operator,
and the desired result is trivial.

If $\alpha<0$, then by Lemma~\ref{1}, we have $R_{s,\alpha/p}f\in A^p$ if and
only if $R^{s+\alpha/p,-\alpha/p}f\in A^p$, which, according to Theorem 2.16 of
\cite{zhu6}, is equivalent to
$$(1-|z|^2)^kR^kR^{s+\alpha/p,-\alpha/p}f\in L^p(\bn,dv).$$
Since $R^k$ commutes with $R^{s+\alpha/p,-\alpha/p}$, the above condition is
equivalent to
$$R^{s+\alpha/p,-\alpha/p}R^kf\in L^p(\bn,dv_{pk}),$$
or
$$(1-|z|^2)^{-\alpha/p}R^{s+\alpha/p,-\alpha/p}R^kf\in L^p(\bn,dv_{pk+\alpha}).$$
Since $\alpha<0$, it follows from Theorem 2.19 of \cite{zhu6} that the above condition 
is equivalent to (\ref{eq2}). This proves the desired result for $\alpha<0$ and
completes the proof of the theorem.
\end{proof}

As a consequence, we obtain the following result which shows that the definition 
of $\apa$ is actually independent of the integer $k$ used. This is of course a
phenomenon that has been well known to experts in the field.

\begin{cor}
Suppose $p>0$ and $\alpha$ is real. Then the following conditions are
equivalent for holomorphic functions $f$ in $\bn$.
\begin{enumerate}
\item[(a)] $f\in\apa$, that is, for some positive integer $k$ with $kp+\alpha>-1$ the 
function $(1-|z|^2)^kR^kf(z)$ is in $\lpa$.
\item[(b)] For every positive integer $k$ with $kp+\alpha>-1$ the function
$(1-|z|^2)^kR^kf(z)$ is in $\lpa$.
\end{enumerate}
\label{11}
\end{cor}

\begin{proof}
This follows from the proof of Theorem~\ref{10}. This also follows from the 
equivalence of (a) and (d) in Theorem 2.16 of \cite{zhu6}.
\end{proof}

Since the polynomials are dense in $A^p$, and since the operators $R^{s,t}$
and $R_{s,t}$ map the set of polynomials onto the set of polynomials, we
conclude from Theorem~\ref{10} that the polynomials are dense in each
space $\apa$.

The following result is a generalization of Theorem~\ref{10}.

\begin{thm}
Suppose $\alpha$ is real, $\beta$ is real, and $p>0$. Let $t=(\alpha-\beta)/p$
and let $s$ be a complex parameter such that neither $n+s$ nor $n+s+t$ is a negative
integer. Then the operator $R_{s,t}$ maps $A^p_\alpha$ boundedly onto
$A^p_\beta$.
\label{12}
\end{thm}

\begin{proof}
We can approximate $s$ by a sequence $\{s_k\}$ of complex numbers such that
each of the operators $R_{s_k,t}$, $R^{s_k+t,\beta/p}$, and $R_{s_k,\alpha/p}$
is well defined. According to Lemmas~\ref{1} and \ref{2}, we have
$$R_{s_k,t}=R^{s_k+t,\beta/p}R_{s_k,\alpha/p}.$$
Since $R^{s_k+t,\beta/p}$ is the inverse of $R_{s_k+t,\beta/p}$, it follows
from Theorem~\ref{10} that each $R_{s_k,t}$ maps $A^p_\alpha$ boundedly onto
$A^p_\beta$. Since $R_{s,t}$ is well defined, an easy limit argument then shows
that $R_{s,t}$ maps $A^p_\alpha$ boundedly onto $A^p_\beta$.
\end{proof}

For any positive $p$ and real $\alpha$ we let $N$ be the smallest nonnegative
integer such that $pN+\alpha>-1$ and define
\begin{equation}
\|f\|_{p,\alpha}=|f(0)|+\left[\inb(1-|z|^2)^{pN}|R^Nf(z)|^p\,\dva(z)\right]^{1/p}
\label{eq3}
\end{equation}
for $f\in\apa$. Then $\apa$ becomes a Banach space when $p\ge1$. For
$0<p<1$ the space $\apa$ is a topological vector space with a complete metric
\begin{equation}
d(f,g)=\|f-g\|_{p,\alpha}^p.
\label{eq4}
\end{equation}
The metric $d$ is invariant in the sense that 
$$d(f,g)=d(f-g,0).$$ 
In particular, $\apa$ is an $F$-space. One of the properties of an $F$-space
that we will use later is that the closed graph theorem is valid for it.

\section{Several Characterizations of $\apa$}

In this section we obtain various characterizations of $\apa$ in terms of
fractional differential operators and in terms of higher order derivatives.

\begin{thm}
Suppose $p>0$ and $\alpha$ is real. Then the following conditions are
equivalent for holomorphic functions $f$ in $\bn$.
\begin{enumerate}
\item[(a)] $f\in\apa$.
\item[(b)] For some nonnegative integer $k$ with $kp+\alpha>-1$ the functions
$(1-|z|^2)^{|m|}\partial^mf(z)$, where $|m|=k$, all belong to $\lpa$.
\item[(c)] For every nonnegative integer $k$ with $kp+\alpha>-1$ the functions
$(1-|z|^2)^{|m|}\partial^mf(z)$, where $|m|=k$, all belong to $\lpa$.
\end{enumerate}
\label{13}
\end{thm}

\begin{proof}
Fix a nonnegative integer $k$ with $pk+\alpha>-1$ and assume that
$$(1-|z|^2)^k\partial^mf(z)\in\lpa$$
for all $|m|=k$, then
$$(1-|z|^2)^k\partial^mf(z)\in\lpa$$
for all $|m|\le k$; see Theorem 2.17 of \cite{zhu6}. Since $R^k$ is a linear partial 
differential operator on $H(\bn)$ with polynomial coefficients (see Lemma~\ref{6}), 
we have
$$(1-|z|^2)^kR^kf(z)\in\lpa,$$
or $f\in\apa$. This proves that condition (b) implies (a). That condition (c)
implies (b) is obvious.

Next assume that $f\in\apa$. Then by Theorem~\ref{10}, the function
$g=R_{\beta,\alpha/p}f$ is in $A^p$, where $\beta$ is a sufficiently large (to be
specified later) positive number. By Proposition~\ref{9}, we have
$$R_{\beta,\alpha/p}f(z)=\inb\frac{g(w)\,dv_{\beta}(w)}{(1-\langle z,w\rangle)^{n
+1+\beta}}.$$
Apply $R^{\beta,\alpha/p}$ to both sides and use Proposition~\ref{3}. We obtain
\begin{equation}
f(z)=\inb\frac{g(w)\,dv_{\beta}(w)}{(1-\langle z,w\rangle)^{n+1+\beta+\alpha/p}}.
\label{eq5}
\end{equation}

If $p\ge1$ and $k$ is any nonnegative integer such that $pk+\alpha>-1$, then
we choose $\beta$ large enough so that
\begin{equation}
-pk<\alpha+1<p\left(\beta+\frac\alpha p\right).
\label{eq6}
\end{equation}
Rewrite the reproducing formula (\ref{eq5}) as
$$f(z)=\inb\frac{(1-|w|^2)^{\beta+\alpha/p}h(w)\,dv(w)}{(1-\langle z,w
\rangle)^{n+1+\beta+\alpha/p}},$$
where
$$h(z)=(1-|z|^2)^{-\alpha/p}g(z).$$
Differentiating under the integral sign, we obtain a positive
constant $C$ (depending on the parameters but not on $f$ and $z$) such that
$$(1-|z|^2)^k|\partial^mf(z)|\le C(1-|z|^2)^k\inb\frac{(1-|w|^2)^{\beta+
\alpha/p}|h(w)|\,dv(w)}{|1-\langle z,w\rangle|^{n+k+1+\beta+\alpha/p}},$$
where $|m|=k$. Since $h\in\lpa$, it follows from (\ref{eq6}) and Proposition~\ref{8}
that the functions $(1-|z|^2)^k\partial^mf(z)$, where $|m|=k$, all
belong to $\lpa$.

The case $0<p<1$ calls for a different proof. In this case, we differentiate under
the integral sign in (\ref{eq5}) and obtain a constant $C>0$ (depending on the
parameters but not on $f$ and $z$) such that
$$(1-|z|^2)^k|\partial^mf(z)|\le C(1-|z|^2)^k\inb\frac{|g(w)|(1-|w|^2)^\beta
\,dv(w)}{|1-\langle z,w\rangle|^{n+k+1+\beta+\alpha/p}},$$
where $|m|=k$. We write
$$\beta=\frac{n+1+\beta'}p-(n+1)$$
and assume that $\beta$ is large enough so that $\beta'>0$. Then we can
apply Lemma 2.15 of \cite{zhu6} to show that the integral
$$\inb\frac{|g(w)|(1-|w|^2)^\beta\,dv(w)}{|1-\langle z,w\rangle|^{n+k+1
+\beta+\alpha/p}}$$
is less than or equal to a positive constant times
$$\left[\inb\left|\frac{g(w)}{(1-\langle z,w\rangle)^{n+k+1
+\beta+\alpha/p}}\right|^p(1-|w|^2)^{\beta'}\,dv(w)\right]^{1/p}.$$
It follows that there exists a positive constant $C'$ such that
$$(1-|z|^2)^{kp}|\partial^mf(z)|^p\le C'(1-|z|^2)^{kp}\inb\frac{|g(w)|^p
(1-|w|^2)^{\beta'}\,dv(w)}{|1-\langle z,w\rangle|^{p(n+k+1+\beta)+\alpha}},$$
where $|m|=k$. Integrate both sides against the measure $dv_\alpha$ and 
apply Fubini's theorem. We see that the integral
$$\inb(1-|z|^2)^{kp}|\partial^mf(z)|^p\,dv_\alpha(z)$$
is less than or equal to $C'$ times
$$\inb|g(w)|^p(1-|w|^2)^{\beta'}\,dv(w)\inb\frac{(1-|z|^2)^{kp+\alpha}\,dv(z)}{|1
-\langle z,w\rangle|^{p(n+k+1+\beta)+\alpha}}.$$
Estimating the inner integral above according to Proposition~\ref{7},
we find another constant $C''>0$ such that
$$\inb(1-|z|^2)^{kp}|\partial^mf(z)|^p\dva(z)\le C''\inb|g(w)|^p\,dv(w)$$
for all $|m|=k$. This proves that (a) implies (c), and completes the proof of
the theorem.
\end{proof}

Note that several special cases of the above theorem are well known. See
Beatrous-Burbea \cite{BB2} or Pavlovic \cite{Pav} for example. In fact,
any nontangential partial differential operator of order $k$ with $C^\infty$
coefficients may be used in place of $R^k$; see Peloso \cite{Peloso}. The 
proof above uses several techniques developed in Zhu \cite{zhu6}.

\begin{thm}
Suppose $p>0$, $\alpha$ is real, and $f$ is holomorphic in $\bn$. Then the
following conditions are equivalent.
\begin{enumerate}
\item[(a)] $f\in\apa$.
\item[(b)] There exists some real $t$ with $pt+\alpha>-1$ such that the function
$(1-|z|^2)^tR^{s,t}f(z)$ is in $\lpa$, where $s$ is any real parameter
such that neither $n+s$ nor $n+s+t$ is a negative integer.
\item[(c)] For every real $t$ with $pt+\alpha>-1$ the function
$(1-|z|^2)^tR^{s,t}f(z)$ is in $\lpa$, where $s$ is any real parameter
such that neither $n+s$ nor $n+s+t$ is a negative integer.
\end{enumerate}
\label{14}
\end{thm}

\begin{proof}
It is obvious that condition (c) implies (b).

To show that condition (b) implies (a), we fix a sufficiently large positive number
$\beta$ and apply Proposition~\ref{9} to write
$$R^{s,t}f(z)=c_{t+\beta}\inb\frac{R^{s,t}f(w)(1-|w|^2)^{t+\beta}\,dv(w)}
{(1-\langle z,w\rangle)^{n+1+t+\beta}},$$
where $c_{t+\beta}$ is a positive constant such that $c_{t+\beta}\,dv_\beta$
is a probability measure on $\bn$. Apply $R^k$ to both sides, where $k$ is a
nonnegative integer such that $kp+\alpha>-1$. Then there exists a polynomial
$h$ of degree $k$ such that
$$R^kR^{s,t}f(z)=\inb\frac{h(\langle z,w\rangle)R^{s,t}f(w)\,dv_{t+\beta}(w)}
{(1-\langle z,w\rangle)^{n+1+k+t+\beta}}.$$
If $\beta$ is chosen so that $\beta-s$ is a sufficiently large positive integer,
we first write
$$h(\langle z,w\rangle)=\sum_{j=0}^kc_j(1-\langle z,w\rangle)^j,$$
then apply the operator $R_{s,t}$ to every term according to the second part of 
Proposition~\ref{5}, and then combine the various terms. The result is that
$$R_{s,t}R^kR^{s,t}f(z)=\inb\frac{g(z,w)R^{s,t}f(w)\,dv_{t+\beta}(w)}
{(1-\langle z,w\rangle)^{n+1+k+\beta}},$$
where $g$ is a polynomial. Since the operators $R_{s,t}$, $R^k$, and $R^{s,t}$
commute with each other, and since $R_{s,t}$ is the inverse of $R^{s,t}$, we 
obtain a constant $C>0$ such that
$$(1-|z|^2)^k|R^kf(z)|\le C(1-|z|^2)^k\inb\frac{(1-|w|^2)^t|R^{s,t}f(w)|
\,dv_\beta(w)}{|1-\langle z,w\rangle|^{n+1+k+\beta}}.$$
We then follow the same arguments as in the proof of Theorem~\ref{13} to show
that the condition
$$(1-|z|^2)^t|R^{s,t}f(z)|\in\lpa$$
implies
$$(1-|z|^2)^kR^kf(z)\in\lpa.$$
This proves that condition (b) implies (a).

To show that condition (a) implies (c), we fix a function $f\in\apa$ and
choose a sufficiently large positive number $\beta$ such that the
function $g=R_{\beta,\alpha/p}f$ is in $A^p$. We then follow the same
arguments as in the proof of Theorem~\ref{13} to finish the proof. The only
adjustment to make here is this: instead of differentiating under the
integral sign, we apply the operator $R^{s,t}$ inside the integral sign
and take advantage of Proposition~\ref{5} (assuming that $\beta$
is chosen so that $\beta-s$ is a positive integer). We leave the details
to the interested reader.
\end{proof}

The above theorem has appeared in several papers before, at least in various
special forms. See Kaptanoglu \cite{Kap1}\cite{Kap2}, Peloso \cite{Peloso},
and Zhu \cite{zhu6}. The book \cite{BB2} of Beatrous and Burbea also contains a 
version of the result for $\alpha>-1$ which is based on a different family of
fractional radial differential operators.

\section{Holomorphic Lipschitz spaces}

The classical Lipschitz space $\Lambda_\alpha$, $0<\alpha<1$, consists
of holomorphic functions $f$ in $\bn$ such that
$$|f(z)-f(w)|\le C|z-w|^\alpha,\qquad z,w\in\bn,$$
where $C$ is a positive constant depending on $f$. It is well known that
a holomorphic function $f$ is in $\Lambda_\alpha$ if and only if there exists
a positive constant $C$ such that
$$(1-|z|^2)^{1-\alpha}|Rf(z)|\le C,\qquad z\in\bn.$$
See Rudin \cite{rudin} and Zhu \cite{zhu6}.

In this section we extend the theory of Lipschitz spaces $\Lambda_\alpha$
to the full range $-\infty<\alpha<\infty$. More specifically, for any real number 
$\alpha$ we let $\Lambda_\alpha$ denote the space of
holomorphic functions $f$ in $\bn$ such that for some nonnegative integer 
$k>\alpha$ the function $(1-|z|^2)^{k-\alpha}R^kf(z)$ is bounded in $\bn$.
We first prove that the definition of $\Lambda_\alpha$ is independent of the 
integer $k$ used.

\begin{lemma}
Suppose $f$ is holomorphic in $\bn$. Then the following conditions are
equivalent.
\begin{enumerate}
\item[(a)] There exists some nonnegative integer $k>\alpha$ such that the function
$(1-|z|^2)^{k-\alpha}R^kf(z)$ is bounded in $\bn$.
\item[(b)] For every nonnegative integer $k>\alpha$ the function
$$(1-|z|^2)^{k-\alpha}R^kf(z)$$ 
is bounded in $\bn$.
\end{enumerate}
\label{15}
\end{lemma}

\begin{proof}
Suppose $k$ is a nonnegative integer with $k>\alpha$. Let $N=k+1$. 

If the function $(1-|z|^2)^{N-\alpha}R^Nf(z)$ is bounded in $\bn$, then an 
elementary integral estimate based on the identity
$$R^kf(z)-R^kf(0)=\int_0^1\frac{R^Nf(tz)}t\,dt$$
shows that the function $(1-|z|^2)^{k-\alpha}R^kf(z)$ is bounded in $\bn$.

Conversely, if the function $(1-|z|^2)^{k-\alpha}R^kf(z)$ is bounded, then
there exists a constant $c>0$ such that
$$R^kf(z)=c\inb\frac{(1-|w|^2)^{k-\alpha}R^kf(w)\,dv(w)}{(1-\langle z,w\rangle)^{n+1
+k-\alpha}};$$
see Proposition~\ref{9}. Taking the radial derivative on both sides, we get
$$R^Nf(z)=C\inb\frac{\langle z,w\rangle(1-|w|^2)^{k-\alpha}R^kf(w)\,dv(w)}{(1-\langle z,
w\rangle)^{n+1+N-\alpha}},$$
where $C=c(n+1+k-\alpha)$. This combined with Proposition~\ref{7}
shows that the function $(1-|z|^2)^{N-\alpha}R^Nf(z)$ is bounded in $\bn$.

Therefore, the function $(1-|z|^2)^{k-\alpha}R^kf(z)$ is bounded if and only if the
function $(1-|z|^2)^{k+1-\alpha}R^{k+1}f(z)$ is bounded, where $k$ is any
nonnegative integer satisfying $k>\alpha$. This clearly proves the desired
result.
\end{proof}

The above lemma is most likely known to experts in the field, although we
could not find a precise reference. In the case $\alpha>0$, the above result
as well as everything else in this section can be found in Zhu \cite{zhu6}.

In what follows we let $k$ be the smallest nonnegative integer greater
than $\alpha$ and define a norm on $\Lambda_\alpha$ by
$$\|f\|_\alpha=|f(0)|+\sup_{z\in\bn}(1-|z|^2)^{k-\alpha}|R^kf(z)|.$$
It is then easy to check that $\Lambda_\alpha$ becomes a nonseparable
Banach space when equipped with this norm.

We write $\bloch=\Lambda_0$. This is called the Bloch space. It is clear that
$f\in\bloch$ if and only if
$$\sup_{z\in\bn}(1-|z|^2)|Rf(z)|<\infty.$$
See \cite{zhu6} for more information about $\bloch$. Our next result shows that 
all the spaces $\Lambda_\alpha$ are isomorphic to the Bloch space.

\begin{thm}
Suppose $s$ is complex and $\alpha$ is real such that neither $n+s$ nor
$n+s+\alpha$ is a negative integer. Then the operator $R^{s,\alpha}$ maps
$\Lambda_\alpha$ onto $\bloch$.
\label{16}
\end{thm}

\begin{proof}
Suppose $f$ is holomorphic in $\bn$. Then $R^{s,\alpha}f$ is in the Bloch
space if and only if the function $(1-|z|^2)^kR^kR^{s,\alpha}f(z)$ is bounded in 
$\bn$, where $k$ is any positive integer. See Lemma~\ref{15} above.

If $f\in\Lambda_\alpha$, then the function 
$$g(z)=(1-|z|^2)^{k-\alpha}R^kf(z)$$ 
is bounded in $\bn$, where $k$ is any positive integer greater than $\alpha$. 
Let $N$ be a sufficiently large positive integer such that the number $\beta$ 
defined by
$$k-\alpha+\beta=s+N$$
has real part greater than $-1$. Then we use Proposition~\ref{9} to write
$$R^kf(z)=c\inb\frac{g(w)\,dv_\beta(w)}{(1-\langle z,w\rangle)^{n+1+k-\alpha+\beta}}.$$
Applying Proposition~\ref{5}, we obtain a polynomial $h$ such that
$$R^{s,\alpha}R^kf(z)=\inb\frac{g(w)h(\langle z,w\rangle)\,dv_\beta(w)}{(1-
\langle z,w\rangle)^{n+1+k+\beta}}.$$
By Proposition~\ref{7}, the function 
$$(1-|z|^2)^kR^{s,\alpha}R^kf(z)=(1-|z|^2)^kR^kR^{s,\alpha}f(z)$$
is bounded in $\bn$, so $R^{s,\alpha}f$ is in the Bloch space.

On the other hand, if $R^{s,\alpha}f$ is in the Bloch space, then by Lemma~\ref{1},
the function $R_{s+\alpha,-\alpha}f$ is in the Bloch space. We fix a suffiently large
positive integer $N$ such that $\beta=N+s+\alpha$ has real part greater than $-1$. 
By part (d) of Theorem 3.4 in \cite{zhu6} (the result there was stated and proved for
real $\beta$, it is clear that the complex case holds as well), there exists a function 
$g\in L^\infty(\bn)$ such that
$$R_{s+\alpha,-\alpha}f(z)=\inb\frac{g(w)\,dv_\beta(w)}{(1-\langle z,w\rangle)^{n+1
+\beta}}.$$
We apply the operator $R^{s+\alpha,-\alpha}$ to both sides and use 
Proposition~\ref{5} to obtain
$$f(z)=\inb\frac{p(\langle z,w\rangle)g(w)\,dv_\beta(w)}{(1-\langle z,w\rangle)^{n+
1+\beta-\alpha}},$$
where $p$ is a polynomial. An easy computation then shows that
$$R^kf(z)=\inb\frac{q(\langle z,w\rangle)g(w)\,dv_\beta(w)}{(1-\langle z,w\rangle)^{n
+1+\beta+k-\alpha}},$$
where $k$ is any positive integer greater than $\alpha$ and $q$ is another 
polynomial. By Proposition~\ref{7}, the function 
$(1-|z|^2)^{k-\alpha}R^kf(z)$ is bounded in $\bn$, namely, $f\in\Lambda_\alpha$.
This completes the proof of the theorem.
\end{proof}

More generally, if $s$ is any complex number such that neither $n+s$ nor
$n+s+\alpha-\beta$ is a negative integer, then the operator $R^{s,\alpha-\beta}$
is bounded invertible operator from $\Lambda_\alpha$ onto $\Lambda_\beta$.
See the proof of Theorem~\ref{12}.

\begin{thm}
Suppose $f$ is holomorphic in $\bn$ and $\alpha$ is real. If $\re\beta>-1$ and
$n+\beta-\alpha$ is not a negative integer, then $f\in\Lambda_\alpha$ if
and only if there exists a function $g\in L^\infty(\bn)$ such that
\begin{equation}
f(z)=\inb\frac{g(w)\,dv_\beta(w)}{(1-\langle z,w\rangle)^{n+1+\beta-\alpha}}
\label{eq7}
\end{equation}
for $z\in\bn$.
\label{17}
\end{thm}

\begin{proof}
If $f$ admits the integral representation (\ref{eq7}), then for any nonnegative
integer $k>\alpha$ we have
$$R^kf(z)=\inb\frac{p(\langle z,w\rangle)g(w)\,dv_\beta(w)}{(1-\langle z,w
\rangle)^{n+1+\beta+k-\alpha}},$$
where $p(z)$ is a certain polynomial of degree $k$. An application of
Proposition~\ref{7} shows that the function 
$(1-|z|^2)^{k-\alpha}R^kf(z)$ is bounded in $\bn$.

On the other hand, if $f\in\Lambda_\alpha$, then by Theorem~\ref{16}, the
function $R^{\beta-\alpha,\alpha}f$ is in the Bloch space. According to
the classical integral representation of functions in the Bloch space (see
Choe \cite{Choe} or part (d) of Theorem 3.4 in Zhu's book \cite{zhu6}), there 
exists a function $g\in L^\infty(\bn)$ such that
$$R^{\beta-\alpha,\alpha}f(z)=\inb\frac{g(w)\,dv_\beta(w)}{(1-\langle z,
w\rangle)^{n+1+(\beta-\alpha)+\alpha}}.$$
Applying the operator $R_{\beta-\alpha,\alpha}$ to both sides and using
Proposition~\ref{3}, we conclude that
$$f(z)=\inb\frac{g(w)\,dv_\beta(w)}{(1-\langle z,w\rangle)^{n+1+\beta-\alpha}}.$$
This completes the proof of the theorem.
\end{proof}

Since the proof of Theorem 3.4 in Zhu \cite{zhu6} is constructive, it follows that 
there exists a bounded linear operator
$$L:\Lambda_\alpha\to L^\infty(\bn)$$
such that the integral representation in (\ref{eq7}) can be given by choosing $g=L(f)$.

\begin{thm}
Suppose $\alpha$ is real and $k$ is a nonnegative integer greater than $\alpha$.
Then a holomorphic function $f$ in $\bn$ belongs to $\Lambda_\alpha$ if and
only if the functions
$$(1-|z|^2)^{k-\alpha}\partial^mf(z),\qquad |m|=k,$$
are all bounded in $\bn$.
\label{18}
\end{thm}

\begin{proof}
If $f\in\Lambda_\alpha$, we apply Theorem~\ref{17} to represent $f$ in the form
$$f(z)=\inb\frac{g(w)\,dv_\beta(w)}{(1-\langle z,w\rangle)^{n+1+\beta-\alpha}},$$
where $g\in L^\infty(\bn)$, $\beta>-1$, and $n+\beta-\alpha$ is not a negative
integer. Differentiate under the integral sign and apply Proposition~\ref{7}. We 
see that the functions $(1-|z|^2)^{k-\alpha}\partial^mf(z)$, where $|m|=k$, are all 
bounded in $\bn$.

Conversely, if the function $(1-|z|^2)^{k-\alpha}\partial^mf(z)$ is bounded in $\bn$
for every $|m|=k$, then it is easy to see that the function 
$(1-|z|^2)^{k-\alpha}\partial^mf(z)$ is bounded in $\bn$ for every $|m|\le k$. Since
$R^k$ is a $k$th order linear partial differential operator on $H(\bn)$ with 
polynomial coefficients (see Lemma~\ref{6}), we see that the function
$(1-|z|^2)^{k-\alpha}R^kf(z)$ is bounded in $\bn$, namely, $f\in\Lambda_\alpha$.
\end{proof}

Various special cases (such as the Bloch space and the case $\alpha\in(0,1)$) of 
the above theorem and the next have been well known. See 
Aleksandrov \cite{Alex}, Choe \cite{Choe}, Nowark \cite{N}, 
Ouyang-Yang-Zhao \cite{OYZ}, Pavlovic \cite{Pav}, Peloso \cite{Peloso}, and 
Zhu \cite{zhu6}, among many others.

\begin{thm}
Suppose $\alpha$ and $t$ are real with $t>\alpha$. If $s$ is a complex parameter
such that neither $n+s$ nor $n+s+t$ is a negative integer, then a holomorphic
function $f$ in $\bn$ belongs to $\Lambda_\alpha$ if and only if the function
$(1-|z|^2)^{t-\alpha}R^{s,t}f(z)$ is bounded in $\bn$.
\label{19}
\end{thm}

\begin{proof}
First assume that the function
$$g(z)=(1-|z|^2)^{t-\alpha}R^{s,t}f(z)$$
is bounded in $\bn$.  By Proposition~\ref{9}, there exists a positive
constant $c$ such that
$$R^{s,t}f(z)=c\inb\frac{g(w)\,dv_\beta(w)}{(1-\langle z,w
\rangle)^{n+1+t+\beta-\alpha}},$$
where $\beta$ is a sufficiently large positive number with $\beta-\alpha=s+N$
for some positive integer $N$. If $k$ is a nonnegative integer greater than
$\alpha$, it is easy to see that there exists a polynomial $p$ of degree $k$
such that
\begin{equation}
R^kR^{s,t}f(z)=\inb\frac{p(\langle z,w\rangle)g(w)\,dv_\beta(w)}{(1-\langle z,
w\rangle)^{n+1+k+t+\beta-\alpha}}.
\label{eq8}
\end{equation}
We decompose
$$p(\langle z,w\rangle)=\sum_{j=0}^kc_j(1-\langle z,w\rangle)^j,$$
apply the operator $R_{s,t}$ to both sides of (\ref{eq8}), use 
Proposition~\ref{5}, and combine the terms. The result is that
$$R_{s,t}R^kR^{s,t}f(z)=\inb\frac{h(z,w)g(w)\,dv_\beta(w)}{(1-\langle z,w
\rangle)^{n+1+k+\beta-\alpha}},$$
where $h$ is a certain polynomial. Since all radial differential operators 
commute, we have
$$R_{s,t}R^kR^{s,t}=R^k.$$
This together with Proposition~\ref{7} shows that
$$|R^kf(z)|\le\frac C{(1-|z|^2)^{k-\alpha}}$$
for some constant $C>0$, that is, $f\in\Lambda_\alpha$.

Next assume that $f\in\Lambda_\alpha$. Let $N$ be a sufficiently large
positive integer and write $\beta-\alpha=s+N$. By Theorem~\ref{17}, there
exists a function $g\in L^\infty(\bn)$ such that
$$f(z)=\inb\frac{g(w)\,dv_\beta(w)}{(1-\langle z,w\rangle)^{n+1+\beta-\alpha}}.$$
According to Proposition~\ref{5}, there exists a polynomial $h$ such that
$$R^{s,t}f(z)=\inb\frac{h(z,w)g(w)\,dv_\beta(w)}{(1-\langle z,w\rangle)^{n+1+\beta
-\alpha+t}}.$$
An application of Proposition~\ref{7} then shows that
$$|R^{s,t}f(z)|\le\frac C{(1-|z|^2)^{t-\alpha}}$$
for some constant $C>0$, that is, the function $(1-|z|^2)^{t-\alpha}R^{s,t}f(z)$
is bounded in $\bn$.
\end{proof}

All results in this section so far are in terms of a certain function being bounded
in $\bn$. We mention that these results remain true when the big oh conditions
are replaced by the corresponding little oh conditions. More specifically, for
each real number $\alpha$, we let $\Lambda_{\alpha,0}$ denote the space of
holomorphic functions $f$ in $\bn$ such that there exists a nonnegative integer
$k>\alpha$ such that the function $(1-|z|^2)^{k-\alpha}R^kf(z)$ is in $\C_0(\bn)$.
Here $\C_0(\bn)$ denotes the space of continuous functions $f$ in $\bn$ with
the property that
$$\lim_{|z|\to1^-}f(z)=0.$$
It can be shown that the definition of $\Lambda_{\alpha,0}$ is independent of 
the integer $k$ used. The special case $\Lambda_{0,0}$ is denoted by $\bloch_0$
and is called the little Bloch space of $\bn$. Clearly, $f\in\bloch_0$ if and only if
$$\lim_{|z|\to1^-}(1-|z|^2)Rf(z)=0.$$

An alternative description of $\Lambda_{\alpha,0}$ is that it is the closure of the
set of polynomials in $\Lambda_\alpha$, or the closure in $\Lambda_\alpha$ of 
the set of functions holomorphic on the closed unit ball.

It is then clear how to state and prove the little oh analogues of all results of
this section. It is also well known that when dealing with the little oh
type results of this section, the space $\C_0(\bn)$ can be replaced by
$\C(\overline\bn)$, the space of functions that are continuous on the closed
unit ball.  We leave out the routine details.

\section{Pointwise Estimates}

We often need to know how fast a function in $\apa$ grows near the boundary.
Using results from the previous section, we obtain optimal pointwise estimates 
for functions in $\apa$.

\begin{thm}
Suppose $p>0$ and $n+1+\alpha>0$. Then there exists a constant $C>0$
(depending on $p$ and $\alpha$) such that
$$|f(z)|\le\frac{C\|f\|_{p,\alpha}}{(1-|z|^2)^{(n+\alpha+1)/p}}$$
for all $f\in\apa$ and $z\in\bn$.
\label{20}
\end{thm}

\begin{proof}
Suppose $f\in\apa$. Then $R^Nf\in A^p_{pN+\alpha}$, where $pN+\alpha>-1$. 
By Theorem 2.1 of \cite{zhu6},
$$(1-|z|^2)^{(n+1+pN+\alpha)/p}|R^Nf(z)|\le C\|f\|_{p,\alpha}$$
for some positive constant $C$ (depending only on $\alpha$). Since
$$\frac{n+1+pN+\alpha}p=N+\frac{n+1+\alpha}p,$$
it follows from Lemma~\ref{15} that there exists a
constant $C'>0$ (depending on $p$ and $\alpha$) such that
$$(1-|z|^2)^{(n+1+\alpha)/p}|f(z)|\le C'\|f\|_{p,\alpha}$$
for all $z\in\bn$.
\end{proof}

In the case $\alpha>-1$ the above theorem can be found in numerous papers
in the literature, including Beatrous-Burbea \cite{BB2} and Vukoti\'c \cite{V}.

It is not hard to see that the estimate given in Theorem~\ref{20} above is
optimal, namely, the exponent $(n+\alpha+1)/p$ cannot be improved.
However,  using polynomial approximations, we can show that
$$\lim_{|z|\to1^-}(1-|z|^2)^{(n+\alpha+1)/p}f(z)=0$$
whenever $f\in\apa$ with $n+1+\alpha>0$. Also, if $\alpha>-1$, then the
constant $C$ can be taken to be $1$; see Theorem 2.1 in \cite{zhu6}.

\begin{thm}
Suppose $p>0$ and $n+1+\alpha<0$. Then every function in $\apa$ is continuous
on the closed unit ball and so is bounded in $\bn$.
\label{21}
\end{thm}

\begin{proof}
Given $f\in\apa$, Theorem~\ref{10} tells us that we can find a function $g\in A^p$ 
such that $f=R^{s,\alpha/p}g$, where $s$ is any real parameter such that neither
$n+s$ nor $n+s+(\alpha/p)$ is a negative integer. By Theorem~\ref{20} and
the remark following it, the function $(1-|z|^2)^{(n+1)/p}g(z)$ is in $\C_0(\bn)$, which, 
according to the little oh version of Lemma~\ref{15}, is the same as 
$g\in\Lambda_{-(n+1)/p,0}$. Let $\beta$ be
a sufficiently large positive number such that
$$\beta+\frac{n+1}p=s+N$$
for some positive integer $N$. We first apply the little oh version of 
Theorem~\ref{17} to find a function $h\in\C_0(\bn)$ such that
$$g(z)=\inb\frac{h(w)\,dv_\beta(w)}{(1-\langle z,w\rangle)^{n+1+\beta+(n+1)/p}}.$$
We then apply the operator $R^{s,\alpha/p}$ to both sides and make use of
Proposition~\ref{5}. The result is
$$f(z)=\inb\frac{p(z,w)h(w)\,dv_\beta(w)}{(1-\langle z,w\rangle)^{n+1+\beta
+(n+1+\alpha)/p}},$$
where $p$ is a polynomial. By part (a) of Proposition~\ref{7}, the integral above 
converges uniformly for $z\in\bn$ and so the function $f(z)$ is continuous on the 
closed unit ball.
\end{proof}

When $n+1+\alpha<0$, functions in $\apa$ are actually much better than just
being continuous on the closed unit ball. For example, it follows from
Theorems~\ref{12}, \ref{19}, and \ref{20} that every function in $\apa$, $n+1+\alpha<0$,
actually belongs to a Lipschitz space $\Lambda_\beta$ for some $\beta>0$.
See Corollary 5.5 of Beatrous-Burbea \cite{BB2} for a slightly different version
of this observation.

Theorems \ref{20} and \ref{21} also follow from Lemmas 5.4 and 5.6 of 
Beatrous-Burbea \cite{BB2}. However, as our next result shows, the estimates in 
\cite{BB2} for the remaining case $\alpha=-(n+1)$ do not seem to be optimal.

\begin{thm}
Suppose $n+1+\alpha=0$ and $f\in\apa$.
\begin{enumerate}
\item[(a)] If $0<p\le 1$, then $f(z)$ is continuous on the closed unit ball. In
particular, $f$ is bounded in $\bn$.
\item[(b)] If $1<p<\infty$ and $1/p+1/q=1$, then there exists a positive constant
$C$ (depending on $p$) such that
$$|f(z)|\le C\left[\log\frac2{1-|z|^2}\right]^{1/q}$$
for all $z\in\bn$.
\end{enumerate}
\label{22}
\end{thm}

\begin{proof}
Note that $A^p_{-(n+1)}=B_p$, the diagonal Besov spaces on $\bn$; see
Chapter 7 of \cite{zhu6}. If $0<p\le 1$, the Besov space $B_p$ is contained in 
$B_1$ (this is well known, and follows easily from Theorem~\ref{32} and the fact
that $l^p\subset l^1$ for $0<p\le 1$).  Since $B_1$ is contained in the
ball algebra (see Theorem 6.8 of \cite{zhu6} for example), we conclude that
$B_p$ is contained in the ball algebra whenever $0<p\le1$.

If $p>1$, we use Theorem 6.7 of \cite{zhu6} to find a function $g\in L^p(\bn,d\tau)$ such that
$$f(z)=\inb\frac{g(w)\,dv(w)}{(1-\langle z,w\rangle)^{n+1}},$$
where 
$$d\tau(z)=\frac{dv(z)}{(1-|z|^2)^{n+1}}$$ 
is the M\"obius invariant measure on $\bn$. Rewrite the above integral
representation as
$$f(z)=\inb\left(\frac{1-|w|^2}{1-\langle z,w\rangle}\right)^{n+1}g(w)\,d\tau(w)$$
and apply H\"older's inequality. We obtain
$$|f(z)|\le\left[\inb|g(w)|^p\,d\tau(w)\right]^{\frac1p}\left[\inb\frac{(1-|w|^2)^{(n+1)(q-1)}}
{|1-\langle z,w\rangle|^{(n+1)q}}\,dv(w)\right]^{\frac1q}.$$
An application of Proposition~\ref{7} to the last integral above yields
the desired estimated for $f(z)$.
\end{proof}

\section{Duality}

A linear functional
$$F:\apa\to\C$$
is said to be bounded if there exists a positive constant $C$ such that
\begin{equation}
|F(f)|\le C\|f\|_{p,\alpha}
\label{eq9}
\end{equation}
for all $f\in\apa$. The dual space of $\apa$, denoted by $\left(\apa\right)^*$, 
is the vector space of all bounded linear functionals on $\apa$. For any
bounded linear functional $F$ on $\apa$ we use $\|F\|$ to denote the smallest
constant $C$ satisfying (\ref{eq9}). It is then easy to check that $\left(\apa\right)^*$ 
becomes a Banach space with this norm, regardless of $p\ge1$ or $p<1$.

By results of the previous section, the point evaluation at any $z\in\bn$ is a
bounded linear functional on $\apa$. Therefore, $\left(\apa\right)^*$ is a nontrivial
Banach space for all $p>0$ and all real $\alpha$.

Results of this section for the case $p>1$ are motivated by the well-known duality 
relation $(A^p)^*=A^q$ for ordinary Bergman spaces under the ordinary volume
integral pairing. Results of this section in the case $0<p\le1$ are motivated by and
are generalizations of various special cases obtained in the papers 
Duren-Romberg-Shields \cite{DRS}, Rochberg \cite{Roch}, Shapiro \cite{shapiro}, 
and Zhu \cite{zhu3}. We also mention that this section overlaps with part of
Section 7 of Kaptanoglu \cite{Kap2}.

\begin{thm}
Suppose $1<p<\infty$, $\alpha$ is real, and $\beta$ is real. If 
$$\frac1p+\frac1q=1,$$
and if $s_1$ and $s_2$ are complex parameters such that both $R_{s_1,\alpha/p}$ 
and $R_{s_2,\beta/q}$ are well-defined operators, then $(\apa)^*=A^q_\beta$ 
(with equivalent norms) under the integral pairing
$$\langle f,g\rangle=\inb R_{s_1,\alpha/p}f\,\overline{R_{s_2,\beta/q}g}\,dv,$$
where $f\in\apa$ and $g\in A^q_\beta$.
\label{23}
\end{thm}

\begin{proof}
This follows from the identities
$$R^{s_1,\alpha/p}A^p=\apa,\qquad R^{s_2,\beta/q}A^q=A^q_\beta,$$
and the well-known duality $(A^p)^*=A^q$ under the integral pairing
$$\langle f,g\rangle=\inb f\,\bar g\,dv,$$
where $f\in A^p_\alpha$ and $g\in A^q_\beta$.
\end{proof}

If $\alpha>-1$ and $\beta>-1$, then the integral pairing
$$\inb R_{s_1,\alpha/p}f\,\overline{R_{s_2,\beta/q}g}\,dv$$
can be replaced by the integral pairing
$$\inb f\,\bar g\,dv_\gamma,\qquad f\in\apa,g\in A^q_\beta,$$
where
\begin{equation}
\gamma=\frac\alpha p+\frac\beta q.
\label{eq10}
\end{equation}
See Theorem 2.12 of \cite{zhu6}. For arbitrary $\alpha$ and $\beta$, we can
also use the integral pairing
$$\lim_{r\to1^-}\inb R_{s,\gamma}f(rz)\,\overline{g(rz)}\,dv(z),
\qquad f\in A^p_\alpha,g\in A^q_\beta,$$
where $\gamma$ is defined by (\ref{eq10}) and $s$ is any complex parameter 
such that the operator $R_{s,\gamma}$ is well defined.

More generally, if $\gamma$ is given by (\ref{eq10}) and if $k$ is a sufficiently 
large positive integer, then the duality $\left(\apa\right)^*=A^q_\beta$ can be
realized with the following integral pairing
$$\langle f,g\rangle_\gamma=f(0)\overline{g(0)}+\inb(1-|z|^2)^kR^kf(z)\,
\overline{(1-|z|^2)^kR^kg(z)}\,dv_\gamma(z),$$
where $f\in\apa$ and $g\in A^q_\beta$. Many other different, but equivalent,
integral pairings are possible.

\begin{thm}
Suppose $0<p\le1$, $\alpha$ is real, and $\beta$ is real. If $s_1$ and $s_2$
are complex parameters such that the operators $R_{s_1,\alpha/p}$ and
$R^{s_2,\beta}$ are well-defined, then $\left(\apa\right)^*=\Lambda_\beta$ under 
the integral pairing
$$\langle f,g\rangle=\lim_{r\to1^-}\inb R_{s_1,\alpha/p}f(rz)\,\overline{R^{s_2,\beta}
g(rz)}\,dv_\gamma(z),$$
where $f\in\apa$, $g\in\Lambda_\beta$, and $\gamma=(n+1)(1/p-1)$.
\label{24}
\end{thm}

\begin{proof}
This follows from the identities
$$R^{s_1,\alpha/p}A^p=\apa,\quad R^{s_2,\beta}\Lambda_\beta=\bloch,$$
and the well-known duality $(A^p)^*=\bloch$ under the integral pairing
$$\langle f,g\rangle=\lim_{r\to1^-}\inb f(rz)\overline{g(rz)}\,dv_\gamma(z).$$
See Theorem 3.17 of \cite{zhu6}.
\end{proof}

Once again, it is easy to come up with other different (but equivalent) duality pairings.
We state two special cases.

\begin{cor}
For any real $\alpha$ we have $(A^1_\alpha)^*=\Lambda_\alpha$ (with equivalent
norms) under the integral pairing
$$\langle f,g\rangle=\lim_{r\to1^-}\inb f(rz)\overline{g(rz)}\,dv(z),$$
where $f\in A^1_\alpha$ and $g\in\Lambda_\alpha$.
\label{25}
\end{cor}

\begin{proof}
Simply choose $s_1=s_2$, $\alpha=\beta$, and $p=1$ in the theorem.
\end{proof}

\begin{cor}
Suppose $\alpha$ is real and $s$ is any complex parameter such that $R^{s,\alpha}$
is well defined. Then $(A^1_\alpha)^*=\bloch$ (with equivalent norms) under
the integral pairing
$$\langle f,g\rangle=\lim_{r\to1^-}\inb R_{s,\alpha}f(rz)\overline{g(rz)}\,dv(z),$$
where $f\in A^1_\alpha$ and $g\in\bloch$.
\label{26}
\end{cor}

\begin{proof}
Simply choose $\beta=0$ in the theorem.
\end{proof}

\begin{thm}
Suppose $\alpha$ and $\beta$ are real. If $s_1$ and $s_2$ are complex parameters
such that the operators $R^{s_1,\alpha}$ and $R^{s_2,\beta}$ are both well defined,
then $(\Lambda_{\beta,0})^*=A^1_\alpha$ (with equivalent norms) under the
integral pairing
$$\langle f,g\rangle=\lim_{r\to1^-}\inb R_{s_1,\alpha}f(rz)\,\overline{R^{s_2,\beta}g(rz)}
\,dv(z),$$
where $f\in A^1_\alpha$ and $g\in\Lambda_{\beta,0}$.
\label{27}
\end{thm}

\begin{proof}
See the proof of Theorem~\ref{24}.
\end{proof}

We also mention two special cases.

\begin{cor}
For any real $\alpha$ we have $(\Lambda_{\alpha,0})^*=A^1_\alpha$ (with
equivalent norms) under the integral pairing
$$\langle f,g\rangle=\lim_{r\to1^-}\inb f(rz)\overline{g(rz)}\,dv(z),$$
where $f\in\Lambda_{\alpha,0}$ and $g\in A^1_\alpha$.
\label{28}
\end{cor}

\begin{proof}
See the proof of Corollary \ref{25}.
\end{proof}

\begin{cor}
Suppose $\alpha$ is real and $s$ is a complex parameter such that the operator 
$R^{s,\alpha}$ is well defined. Then $(\bloch_0)^*=A^1_\alpha$ (with equivalent 
norms) under the integral pairing
$$\langle f,g\rangle=\lim_{r\to1^-}\inb R_{s,\alpha}f(rz)\overline{g(rz)}\,dv(z),$$
where $f\in\bloch_0$ and $g\in A^1_\alpha$.
\label{29}
\end{cor}

\begin{proof}
Just set $\beta=0$ in the theorem.
\end{proof}

\section{Integral Representations}

In this section we focus on the case $1\le p<\infty$ and show that each
space $\apa$ is a quotient space of $L^p(\bn,dv_\beta)$. We do this using
Bergman type projections.

Integral representations of functions in Bergman spaces of $\bn$ started in
Forelli-Rudin \cite{fr} and have seen several generalizations; see Choe \cite{Choe}, 
Kaptanoglu \cite{Kap2}, and Zhu \cite{zhu6}. The next result appears to be new 
even in the case of unweighted Bergman spaces of the unit disk.

\begin{thm}
Suppose $p\ge1$ and $\alpha$ is real. If $\gamma$ and $\lambda$ are
complex parameters satisfying the following two conditions,
\begin{enumerate}
\item[(a)] $p(\re\gamma+1)>\re\lambda+1$,
\item[(b)] $n+\gamma+(\alpha-\lambda)/p$ is not a negative integer,
\end{enumerate}
then a holomorphic function $f$ in $\bn$ belongs to $\apa$ if and only if
\begin{equation}
f(z)=\inb\frac{g(w)\,dv_\gamma(w)}{(1-\langle z,w\rangle)^{n+1+\gamma
+(\alpha-\lambda)/p}}
\label{eq11}
\end{equation}
for some $g\in L^p(\bn,dv_\lambda)$.
\label{30}
\end{thm}

\begin{proof}
Suppose that the parameters satisfy conditions (a) and (b). Let
$$\beta=\gamma-\frac{\lambda}p.$$
Then $\lambda=p(\gamma-\beta)$. Note that condition (a) is equivalent to 
$p(\re\beta+1)>1$. In particular, $\re\beta>-1$ and $n+\beta$ is not a negative integer.
Also, condition (b) is equivalent to the condition that $n+(\alpha/p)+\beta$ is not
a negative integer. So the operators $R^{\beta,\alpha/p}$ and $R_{\beta,\alpha/p}$
are well defined.

If $f\in\apa$, then by Theorem~\ref{10}, the function $R_{\beta,\alpha/p}f$ is
in $A^p$. It follows from Theorem 2.11 of \cite{zhu6} (note that the result there
was stated and proved for real parameters, the case of complex parameters
is proved in exactly the same way) and the condition
$p(\re\beta+1)>1$ that there exists a function $h\in L^p(\bn,dv)$ such that
$$R_{\beta,\alpha/p}f(z)=\inb\frac{h(w)\,dv_\beta(w)}{(1-\langle z,w\rangle)^{n+1+
\beta}}.$$
Apply the operator $R^{\beta,\alpha/p}$ to both sides and use Proposition~\ref{3}. 
Then
$$f(z)=\inb\frac{h(w)\,dv_\beta(w)}{(1-\langle z,w\rangle)^{n+1+\beta+(\alpha/p)}}.$$
Let
$$g(w)=(1-|w|^2)^{\beta-\gamma}h(w).$$
Then $g\in L^p(\bn,dv_\lambda)$ and
$$f(z)=\inb\frac{g(w)\,dv_\gamma(w)}{(1-\langle z,w\rangle)^{n+1+\gamma
+(\alpha-\lambda)/p}}.$$

The above arguments can be reversed. So any function represented by (\ref{eq11})
is necessarily a function in $\apa$. This completes the proof of the theorem.
\end{proof}

Once again, the proof of Theorem 2.11 of \cite{zhu6} is constructive. So there exists
a bounded linear operator
$$L:\apa\to L^p(\bn,dv_\lambda)$$
such that the integral representation in (\ref{eq11}) can be achieved with the
choice $g=L(f)$.

If condition (b) above is not satisfied, then
$$\frac1{(1-\langle z,w\rangle)^{n+1+\gamma+(\alpha-\lambda)/p}}=
(1-\langle z,w\rangle)^k$$
for some nonnegative integer $k$, and any function represented by (\ref{eq11}) is
a polynomial of degree less than or equal to $k$. In this case, the integral
representation (\ref{eq11}) cannot possibly give rise to all functions in $\apa$.
This shows that condition (b) is essential for the theorem.

We can also show that condition (a) is essential. In fact, if every function
$g\in L^p(\bn,dv_\lambda)$ gives rise to a function $f$ in $\apa$ via the
integral representation (\ref{eq11}), then we can apply the operator
$R^{\gamma+(\alpha-\lambda)/p,k}$ to both sides of (\ref{eq11}) and use
Theorem~\ref{14} to infer that the operator
$$Tg(z)=(1-|z|^2)^k\inb\frac{g(w)\,dv_\gamma(w)}{(1-\langle z,w\rangle)^{n+1+k
+\gamma+(\alpha-\lambda)/p}}$$
maps $L^p(\bn,dv_\lambda)$ boundedly into $L^p(\bn,dv_\alpha)$, where $k$
is any nonnegative integer such that $pk+\alpha>-1$. Write
$$g(w)=(1-|w|^2)^{(\alpha-\lambda)/p}h(w).$$
Then $g\in L^p(\bn,dv_\lambda)$ if and only if $h\in L^p(\bn,dv_\alpha)$. 
It follows that the operator
$$Sh(z)=(1-|z|^2)^k\inb\frac{(1-|w|^2)^{\gamma+(\alpha-\lambda)/p}h(w)\,dv(w)}
{(1-\langle z,w\rangle)^{n+1+k+\gamma+(\alpha-\lambda)/p}}$$
maps $L^p(\bn,dv_\alpha)$ boundedly into $L^p(\bn,\dva)$. By Proposition~\ref{8}, 
the parameters must satisfy the conditions
$$-pk<\alpha+1<p\,\re\left(\gamma+\frac{\alpha-\lambda}p+1\right).$$
It is easy to see that these two conditions are the same as the following
two conditions:
$$pk+\alpha>-1,\qquad p(\re\gamma+1)>\re\lambda+1.$$
Therefore, the conditions in Theorem~\ref{30} above are best possible.

\begin{cor}
Suppose $p\ge1$ and $\alpha$ is real. If $\beta$ is any complex parameter 
such that 
\begin{enumerate}
\item[(a)] $p(\re\beta+1)>\alpha+1$,
\item[(b)] $n+\beta$ is not a negative integer,
\end{enumerate}
then a holomorphic function $f$ in $\bn$ belongs to $\apa$ if and only if
$$f(z)=\inb\frac{g(w)\,dv_\beta(w)}{(1-\langle z,w\rangle)^{n+1+\beta}}$$
for some $g\in L^p(\bn,\dva)$.
\label{31}
\end{cor}

\begin{proof}
Simply set $\gamma=\beta$ and $\lambda=\alpha$ in the theorem.
\end{proof}

\section{Atomic Decomposition}

Atomic decomposition for the Bergman spaces $\apa$ was first obtained in
Coifman-Rochberg \cite{cr} in the case $\alpha>-1$. This turns out to be a powerful 
theorem in the theory of Bergman spaces. We now generalize the result to all $\apa$. 
We will also obtain atomic decomposition for the generalized holomorphic
Lipschitz spaces $\Lambda_\alpha$.

\begin{thm}
Suppose $p>0$, $\alpha$ is real, and $b$ is real. If $b$ is neither $0$ nor a 
negative integer, and
\begin{equation}
b>n\max\left(1,\frac1p\right)+\frac{\alpha+1}p,
\label{eq12}
\end{equation}
then there exists a sequence $\{a_k\}$ in $\bn$ such that a holomorphic function 
$f$ in $\bn$ belongs to $\apa$ if and only if
\begin{equation}
f(z)=\sum_{k=1}^\infty c_k\frac{(1-|a_k|^2)^{b-(n+1+\alpha)/p}}{(1-\langle z,
a_k\rangle)^b}
\label{eq13}
\end{equation}
for some sequence $\{c_k\}\in l^p$.
\label{32}
\end{thm}

\begin{proof}
Note that the condition in (\ref{eq12}) implies that
$$b-\frac{\alpha}p>n\max\left(1,\frac1p\right)+\frac1p>n.$$
This, together with the assumption that $b$ is neither $0$ nor a negative integer,
shows that the operators $R^{s,\alpha/p}$ and $R_{s,\alpha/p}$ are well defined,
where $s$ is determined by
$$b=n+1+s+\frac{\alpha}p.$$

Also note that the condition in (\ref{eq12}) implies that
$$b'>n\max\left(1,\frac1p\right)+\frac1p,$$
where $b'=b-(\alpha/p)$. By Theorem 2.30 of \cite{zhu6}, there exists a sequence 
$\{a_k\}$ such that $f\in A^p$ if and only if
\begin{equation}
f(z)=\sum_{k=1}^\infty c_k\frac{(1-|a_k|^2)^{b'-(n+1)/p}}{(1-\langle z,a_k\rangle)^{b'}}
\label{eq14}
\end{equation}
for some sequence $\{c_k\}\in l^p$.

If $f$ is given by (\ref{eq13}), then
$$R_{s,\alpha/p}f(z)=\sum_{k=1}^\infty c_k\frac{(1-|a_k|^2)^{b-(n+1+\alpha)/p}}
{(1-\langle z,a_k\rangle)^{b-\alpha/p}},$$
or
$$R_{s,\alpha/p}f(z)=\sum_{k=1}^\infty c_k\frac{(1-|a_k|^2)^{b'-(n+1)/p}}
{(1-\langle z,a_k\rangle)^{b'}}.$$
According to the previous paragraph, we have $R_{s,\alpha/p}f\in A^p$. Combining
this with Theorem~\ref{10}, we conclude that $f\in\apa$.

The above arguments can be reversed, showing that every function $f\in\apa$
admits an atomic decomposition (\ref{eq13}). This completes the proof of the theorem.
\end{proof}

Recall that when $\alpha=-(n+1)$, the resulting spaces $\apa$ are nothing but the
diagonal Besov spaces $B_p$. Atomic decompositions for Besov spaces have
also been obtained in Frazier-Jawerth \cite{FJ} and Peloso \cite{Peloso}.

It can be shown that the assumptions on the parameters in the above theorem are
optimal. It can also be shown that for $f\in\apa$, we have
$$\|f\|^p_{p,\alpha}\sim\inf\sum_{k=1}^\infty|c_k|^p,$$
where the infimum is taken over all sequences $\{c_k\}$ satisfying the representation
(\ref{eq13}).

The atomic decomposition for functions in the Bloch space was first obtained in 
Rochberg \cite{Roch}. As a consequence of atomic decomposition for the Bloch 
space, we now obtain an atomic decomposition for functions in the generalized 
Lipschitz spaces.

\begin{thm}
Suppose $\alpha$ and $b$ are real parameters with the following two properties:
\begin{enumerate}
\item[(a)] $b+\alpha>n$.
\item[(b)] $b$ is neither $0$ nor a negative integer.
\end{enumerate}
Then there exists a sequence $\{a_k\}$ in $\bn$ such that a holomorphic function $f$ 
in $\bn$ belongs to $\Lambda_\alpha$ if and only if
\begin{equation}
f(z)=\sum_{k=1}^\infty c_k\frac{(1-|a_k|^2)^{b+\alpha}}{(1-\langle z,a_k\rangle)^b}
\label{eq15}
\end{equation}
for some sequence $\{c_k\}\in l^\infty$.
\label{33}
\end{thm}

\begin{proof}
Choose $s$ so that $b=n+1+s$. Then the operators $R^{s,\alpha}$ and
$R_{s,\alpha}$ are well defined. 

Let $b'=b+\alpha$. Then a function $f$ is represented by (\ref{eq15}) if and only if
$$R^{s,\alpha}f(z)=\sum_{k=1}^\infty c_k\frac{(1-|a_k|^2)^{b'}}{(1-\langle z,a_k
\rangle)^{b'}}$$
for some $\{c_k\}\in l^\infty$. Since $R^{s,\alpha}\Lambda_\alpha=\bloch$, the 
desired result then follows from the atomic decomposition for the Bloch space; 
see Theorem 3.23 of \cite{zhu6}.
\end{proof}

Once again, the assumptions on the parameters $b$ and $\alpha$ in the above
theorem are best possible.

A little oh version of Theorem~\ref{33} also holds, giving the atomic decomposition 
for the space $\Lambda_{\alpha,0}$. The only adjustment to be made is to replace
the sequence space $l^\infty$ by $c_0$ (consisting of sequences that tend to $0$).
We omit the details.

As a corollary of atomic decomposition, we prove the following embedding of
weighted Bergman spaces which is well known and very useful in the special
case $\alpha>-1$; see Aleksandrov \cite{Alex}, Beatrous-Burbea \cite{BB2},
 Rochberg \cite{Roch}, and Lemma 2.15 of Zhu \cite{zhu6}.

\begin{thm}
Suppose $0<p\le1$ and $\alpha$ is real. If 
$$\beta=\frac{n+1+\alpha}p-(n+1),$$
then $\apa$ is continuously contained in $A^1_\beta$.
\label{34}
\end{thm}

\begin{proof}
Suppose $0<p\le1$ and fix any positive integer $b$ such that $b>(n+1+\alpha)/p$. 
If $f\in\apa$, then there exists a sequence $\{c_k\}\in l^p\subset l^1$ such that
$$f(z)=\sum_{k=1}^\infty c_k\frac{(1-|a_k|^2)^{b-(n+1+\alpha)/p}}{(1-\langle z,
a_k\rangle)^b},$$
where $\{a_k\}$ is a certain sequence in $\bn$. For any $k\ge1$ write
$$f_k(z)=\frac1{(1-\langle z,a_k\rangle)^b}.$$
Then
$$\|f\|_{1,\beta}\le\sum_{k=1}^\infty|c_k|(1-|a_k|^2)^{b-(n+1+\alpha)/p}\|f_k\|_{1,\beta}.$$
An easy computation shows that
$$R^Nf_k(z)=\frac{h(\langle z,a_k\rangle)}{(1-\langle z,a_k\rangle)^{b+N}},$$
where $N$ is the smallest nonnegative integer with $N+\beta>-1$ and $h$ is a 
polynomial of degree $N$. It follows that
$$\inb(1-|z|^2)^N|R^Nf_k(z)|\,dv_\beta(z)\le C\inb\frac{(1-|z|^2)^{N+\beta}\,dv(z)}
{|1-\langle z,a_k\rangle|^{b+N}},$$
where $C$ is a positive constant (independent of $k$). Estimating the second integral 
above by Proposition~\ref{7}, we obtain
$$\inb(1-|z|^2)^N|R^Nf_k(z)|\,dv_\beta(z)\le\frac{C'}{(1-|a_k|^2)^{b-(n+1+\alpha)/p}},$$
where $C'$ is another positive constant independent of $k$. This shows that
$$\|f\|_{1,\beta}\le C'\sum_{k=1}^\infty|c_k|<\infty,$$
completing the proof of the theorem.
\end{proof}

The above theorem can also be proved without appealing to atomic decomposition.
In fact, if $k$ is a sufficiently large positive integer (such that $kp+\alpha>-1$ and
$k+\beta>-1$), then the condition $f\in\apa$, $0<p\le1$, implies that 
$R^kf\in A^p_{\alpha'}$, where $\alpha'=kp+\alpha$. By Lemma 2.15 of \cite{zhu6}, 
we have $R^kf\in A^1_{\beta'}$, where
$$\beta'=\frac{n+1+\alpha'}p-(n+1)=k+\frac{n+1+\alpha}p-(n+1)=k+\beta,$$
or equivalently, the function $(1-|z|^2)^kR^kf(z)$ belongs to $L^1(\bn,dv_\beta)$,
that is, $f\in A^1_\beta$.

\begin{thm}
Suppose $p>0$ and $\alpha$ is real. If $q$ and $r$ are positive numbers satisfying
$$\frac1p=\frac1q+\frac1r,$$
then every function $f\in\apa$ admits a decomposition
$$f(z)=\sum_{k=1}^\infty g_k(z)h_k(z),$$
where each $g_k$ is in $A^q_\alpha$ and each $h_k$ is in $A^r_\alpha$. Furthermore,
if $0<p\le1$, then
$$\sum_{k=1}^\infty\|g_k\|_{q,\alpha}\|h_k\|_{r,\alpha}\le C\|f\|_{p,\alpha},$$
where $C$ is a positive constant independent of $f$.
\label{35}
\end{thm}

\begin{proof}
Consider the function
$$f(z)=\frac1{(1-\langle z,a\rangle)^b},$$
where $a\in\bn$ and $b$ is the constant from Theorem~\ref{32}. We can write
$f=gh$, where
$$g(z)=\frac1{(1-\langle z,a\rangle)^{bp/q}},\quad h(z)=\frac1{(1-\langle z,a\rangle)^{bp/r}}.$$
If $k$ is a sufficiently large positive integer, then it follows from Proposition~\ref{7} that
\begin{eqnarray*}
\|f\|_{p,\alpha}&\sim&\left[\inb\frac{(1-|z|^2)^{pk}\,dv_\alpha(z)}{|1-\langle z,a\rangle|^{p(b
+k)}}\right]^{1/p}\\
&\sim&\left[\frac1{(1-|a|^2)^{pb-(n+1+\alpha)}}\right]^{1/p}\\
&=&\frac1{(1-|a|^2)^{b-(n+1+\alpha)/p}}.
\end{eqnarray*}
Similar computations show that
$$\|g\|_{q,\alpha}\sim\frac1{(1-|a|^2)^{(bp-n-1-\alpha)/q}}$$
and
$$\|h\|_{r,\alpha}\sim\frac1{(1-|a|^2)^{(bp-n-1-\alpha)/r}}.$$
It follows that
$$\|f\|_{p,\alpha}\sim\|g\|_{q,\alpha}\|h\|_{r,\alpha}.$$
The desired result then follows from Theorem~\ref{32} and the fact that
$$\sum_{k=1}^\infty|c_k|\le C\left(\sum_{k=1}^\infty|c_k|^p\right)^{1/p}$$
when $0<p\le1$. See the proof of Corollary 2.33 in \cite{zhu6} as well.
\end{proof}

When $\alpha>-1$, the above theorem can be found in Coifman-Rochberg \cite{cr}
and Rochberg \cite{Roch}.

\section{Complex Interpolation}

In this section we determine the complex interpolation space of two generalized
weighted Bergman spaces. We also determine the complex interpolation space
between a weighted Bergman space and a Lipschitz space.

Throughout this section we let
$$S=\{z=x+iy\in\C:0<x<1\},\quad \overline S=\{z=x+iy\in\C:0\le x\le1\}.$$
Thus $S$ is an open strip in the complex plane and $\overline S$ is its closure.
We denote the two boundary lines of $S$ by
$$L(S)=\{z=x+iy\in\C:x=0\},\quad R(S)=\{z=x+iy\in\C:x=1\}.$$

The complex method of interpolation is based on Hadamard's
three lines theorem, which states that if $f$ is a function that is continuous
on $\overline S$, bounded on $\overline S$, and analytic in $S$, then
$$\sup_{\re z=\theta}|f(z)|\le(\sup_{\re z=0}|f(z)|)^{1-\theta}(\sup_{\re z=1}|f(z)|)^\theta$$
for any $\theta\in(0,1)$.

Let $X$ and $Y$ be two Banach spaces of holomorphic functions in $\bn$. Then 
$X+Y$ becomes a Banach space with the norm
$$\|f\|_{X+Y}=\inf(\|g\|_X+\|h\|_Y),\qquad f\in X+Y,$$
where the infimum is taken over all decompositions $f=g+h$ with $g\in X$ and 
$h\in Y$. If $\theta\in(0,1)$, the complex interpolation space $[X,Y]_\theta$
consists of holomorphic functions $f$ in $\bn$ with the following properties:
\begin{enumerate}
\item[(1)] There exists a function $\zeta\mapsto f_\zeta$ from $\overline S$ into
the Banach space $X+Y$ that is analytic in $S$, continuous on $\overline S$, 
and bounded on $\overline S$.
\item[(2)] $f_\theta=f$.
\item[(3)] The function $\zeta\mapsto f_\zeta$ is bounded and continuous
from $L(S)$ into $X$.
\item[(4)] The function $\zeta\mapsto f_\zeta$ is bounded and continuous
from $R(S)$ into $Y$.
\end{enumerate}
The space $[X,Y]_\theta$ is a Banach space with the following norm:
$$\|f\|_\theta=\inf\max(\sup_{\re\zeta=0}\|f_\zeta\|_X,\sup_{\re\zeta=1}\|f_\zeta\|_Y),$$
where the infimum is taken over all $f_\zeta$ satisfying conditions (1) through (4)
above. See Bergh-L\"ofstr\"om \cite{bl} and Bennett-Sharpley \cite{bs} for more 
information about complex interpolation.

The complex method of interpolation spaces is functorial in the sense that if
$$T:X+Y\to X'+Y'$$
is a linear operator with the property that $T$ maps $X$ boundedly into $X'$ and
$T$ maps $Y$ boundedly into $Y'$, then $T$ also maps $[X,Y]_\theta$ boundedly
into $[X',Y']_\theta$ for each $\theta\in(0,1)$.

The most classical example of complex interpolation spaces concerns $L^p$ spaces 
(over any measure space). More specifically, if $1\le p_0<p_1\le\infty$ and
$$\frac1p=\frac{1-\theta}{p_0}+\frac\theta{p_1}$$
for some $0<\theta<1$, then
$$\left[L^{p_0},L^{p_1}\right]_\theta=L^p$$
with equal norms.  

More generally, if $w_0$ and $w_1$ are weight functions of a measure $\mu$,
and if $1\le p_0<p_1<\infty$, then for any $\theta\in(0,1)$ we have
$$\left[L^{p_0}(w_0),L^{p_1}(w_1)\right]_\theta=L^p(w)$$
with equal norms, provided that
$$\frac1p=\frac{1-\theta}{p_0}+\frac\theta{p_1}$$
and 
$$w^{\frac1p}=w_0^{\frac{1-\theta}{p_0}}w_1^{\frac\theta{p_1}}.$$
This is usually referred to as the Stein-Weiss interpolation theorem. 
See Stein-Weiss \cite{sw}.

\begin{thm}
Suppose $\alpha$ and $\beta$ are real. If $1\le p_0\le p_1<\infty$ and
$$\frac1p=\frac{1-\theta}{p_0}+\frac\theta{p_1}$$
for some $\theta\in(0,1)$, then
$$\left[A^{p_0}_\alpha,A^{p_1}_\beta\right]_\theta=A^p_\gamma$$
with equivalent norms, where $\gamma$ is determined by
$$\frac\gamma p=\frac\alpha{p_0}(1-\theta)+\frac\beta{p_1}\,\theta.$$
\label{36}
\end{thm}

\begin{proof}
It is clear that $1\le p<\infty$. We fix a large positive number $s$ such that
\begin{equation}
p(s+1)>\gamma+1,\quad p_0(s+1)>\alpha+1,\quad p_1(s+1)>\beta+1.
\label{eq16}
\end{equation}
Then by Corollary~\ref{31}, the integral operator
$$Tg(z)=\inb\frac{g(w)\,dv_s(w)}{(1-\langle z,w\rangle)^{n+1+s}}$$
maps $L^p(\bn,dv_\gamma)$ boundedly onto $A^p_\gamma$; it maps
$L^{p_0}(\bn,dv_\alpha)$ boundedly onto $A^{p_0}_\alpha$; and it maps
$L^{p_1}(\bn,dv_\beta)$ boundedly onto $A^{p_1}_\beta$.
It follows from the functorial property of complex interpolation and the
Stein-Weiss interpolation theorem that $T$ maps the space
$$\left[L^{p_0}(\bn,dv_\alpha),L^{p_1}(\bn,dv_\beta)\right]_\theta=
L^p(\bn,dv_\gamma)$$
boundedly into $\left[A^{p_0}_\alpha,A^{p_1}_\beta\right]_\theta$. Since
$TL^p(\bn,dv_\gamma)=A^p_\gamma$, we conclude that
$$A^p_\gamma\subset\left[A^{p_0}_\alpha,A^{p_1}_\beta\right]_\theta,$$
and the inclusion is continuous.

On the other hand, if $k$ is a sufficiently large positive integer, the operator
$L$ defined by
$$L(f)(z)=(1-|z|^2)^kR^kf(z),\qquad f\in H(\bn),$$
maps $A^{p_0}_\alpha$ boundedly into $L^{p_0}(\bn,dv_\alpha)$; and it maps
$A^{p_1}_\beta$ boundedly into $L^{p_1}(\bn,dv_\beta)$. By the functorial
property of complex interpolation and the Stein-Weiss interpolation theorem,
the operator $L$ also maps $\left[A^{p_0}_\alpha,A^{p_1}_\beta\right]_\theta$ boundedly
into $L^p(\bn,dv_\gamma)$. Equivalently, if $f\in\left[A^{p_0}_\alpha,A^{p_1}_\beta
\right]_\theta$, then the function $(1-|z|^2)^kR^kf(z)$ belongs to $L^p(\bn,dv_\gamma)$,
that is, $f\in A^p_\gamma$. We conclude that
$$\left[A^{p_0}_\alpha,A^{p_1}_\beta\right]_\theta\subset A^p_\gamma,$$
and the inclusion is continuous. This completes the proof of the theorem.
\end{proof}

\begin{cor}
Suppose $\alpha$ is real, $\beta$ is real, $1\le p<\infty$, and $0<\theta<1$. Then
$$\left[A^p_\alpha,A^p_\beta\right]_\theta=A^p_\gamma,$$
where $\gamma=\alpha(1-\theta)+\beta\,\theta$.
\label{37}
\end{cor}

\begin{thm}
Suppose $\alpha$ and $\beta$ are real. If $1\le p<\infty$ and $0<\theta<1$, then
$$\left[A^p_\alpha,\Lambda_\beta\right]_\theta=A^q_\gamma$$
with equivalent norms, where $q=p/(1-\theta)$ and $\gamma=\alpha-q\beta\theta$.
\label{38}
\end{thm}

\begin{proof}
First we consider the linear operator
$$Tf(z)=\inb\frac{f(w)\,dv_{s+\beta}(w)}{(1-\langle z,w\rangle)^{n+1+s}},$$
where $s$ is a fixed and sufficiently large positive number. By Theorem~\ref{17}
and Theorem~\ref{30}, the operator $T$ maps $L^\infty(\bn)$ boundedly onto
$\Lambda_\beta$; and it maps $L^p(\bn,dv_{\alpha+p\beta})$ boundedly onto
$A^p_\alpha$. Since
$$\frac1q=\frac{1-\theta}p+\frac\theta\infty,$$
it follows that $T$ maps the space
$$L^q(\bn,dv_{\alpha+p\beta})=\left[L^p(\bn,dv_{\alpha+p\beta}),L^\infty(\bn)\right]_\theta$$
boundedly into $\left[A^p_\alpha,\Lambda_\beta\right]_\theta$. But according to 
Theorem~\ref{30} again, we have $TL^q(\bn,dv_{\alpha+p\beta})=A^q_\gamma$. So
$$A^q_\gamma\subset\left[A^p_\alpha,\Lambda_\beta\right]_\theta,$$
and the inclusion is continuous.

Next we consider the linear operator
$$Lf(z)=(1-|z|^2)^{k-\beta}R^kf(z),\qquad f\in H(\bn),$$
where $k$ is a fixed and sufficiently large positive integer. The operator $L$ maps
$A^p_\alpha$ boundedly into $L^p(\bn,dv_{\alpha+p\beta})$; and it maps $\Lambda_\beta$
boundedly into $L^\infty(\bn)$. Therefore, $L$ also maps $\left[A^p_\alpha,
\Lambda_\beta\right]_\theta$ boundedly into $L^q(\bn,dv_{\alpha+p\beta})$, that is,
$f\in[A^p_\alpha,\Lambda_\beta]_\theta$ implies that the function
$(1-|z|^2)^{k-\beta}R^kf(z)$ is in $L^q(\bn,dv_{\alpha+p\beta})$, which is the same as
the function $(1-|z|^2)^kR^kf(z)$ being in $L^q(\bn,dv_\gamma)$, or $f\in A^q_\gamma$.
We conclude that
$$\left[A^p_\alpha,\Lambda_\beta\right]_\theta\subset A^q_\gamma,$$
and the inclusion is continuous. This completes the proof of the theorem.
\end{proof}

\begin{thm}
Suppose $\alpha$ is real, $\beta$ is real, and $0<\theta<1$. Then
$$\left[\Lambda_\alpha,\Lambda_\beta\right]_\theta=\Lambda_\gamma$$
with equivalent norms, where $\gamma=\alpha(1-\theta)+\beta\,\theta$.
\label{39}
\end{thm}

\begin{proof}
Fix a sufficiently large positive number $s$. If $f\in\Lambda_\gamma$, there 
exists a function $g\in L^\infty(\bn)$ such that
$$f(z)=\inb\frac{g(w)\,dv_s(w)}{(1-\langle z,w\rangle)^{n+1+s-\gamma}};$$
see Theorem~\ref{17}. For any $\zeta\in\overline S$ we define
$$f_\zeta(z)=\inb\frac{g(w)(1-|w|^2)^{\alpha(1-\zeta)+\beta\zeta-\gamma}
\,dv_s(w)}{(1-\langle z,w\rangle)^{n+1+s-\gamma}}.$$
Let $k$ be a sufficiently large positive integer. Then it follows easily from
Proposition~\ref{7} that the map $\zeta\mapsto f_\zeta$ is a bounded continuous
function from $\overline S$ into 
$$\Lambda_\alpha+\Lambda_\beta=\Lambda_{\min(\alpha,\beta)},$$
and its restriction to $S$ is analytic. Also, the map $\zeta\mapsto f_\zeta$
is a bounded continuous function from $L(S)$ into $\Lambda_\alpha$, and
from $R(S)$ into $\Lambda_\beta$. Since $f_\theta=f$, we conclude from the
definition of complex interpolation that 
$f\in [\Lambda_\alpha,\Lambda_\beta]_\theta$. This shows that
$$\Lambda_\gamma\subset\left[\Lambda_\alpha,\Lambda_\beta\right]_\theta$$
and the inclusion is continuous.

On the other hand, if $f\in\left[\Lambda_\alpha,\Lambda_\beta\right]_\theta$, then
there exists a family of functions $f_\zeta$, where $\zeta\in\overline S$, such that
\begin{enumerate}
\item[(a)] $\zeta\mapsto f_\zeta$ is a bounded continuous function from 
$\overline S$ into $\Lambda_{\min(\alpha,\beta)}$ whose restriction to $S$ is
analytic.
\item[(b)] $\zeta\mapsto f_\zeta$ is a bounded continuous function from $L(S)$
into $\Lambda_\alpha$.
\item[(c)] $\zeta\mapsto f_\zeta$ is a bounded continuous function from $R(S)$
into $\Lambda_\beta$.
\item[(d)] $f_\theta=f$.
\end{enumerate}
Let $k$ be a positive integer with $k>\max(\alpha,\beta)$ and consider the
functions
$$g_\zeta(z)=(1-|z|^2)^{k-\alpha(1-\zeta)-\beta\zeta}R^kf_\zeta(z),
\qquad z\in\bn,\zeta\in\overline S.$$

By conditions (b) and (c) of the previous paragraph, there exist finite positive 
constants $M_0$ and $M_1$ such that
\begin{equation}
\sup_{z\in\bn,\zeta\in L(S)}|g_\zeta(z)|=M_0,\quad
\sup_{z\in\bn,\zeta\in R(S)}|g_\zeta(z)|=M_1.
\label{eq17}
\end{equation}
For any fixed point $z\in\bn$, it follows from condition (a) of the previous
paragraph that the function $F(\zeta)=g_\zeta(z)$ is a bounded
continuous function on $\overline S$ whose restriction to $S$ is analytic.
Moreover, it follows from (\ref{eq17}) that $|F(\zeta)|\le M_0$ for $\zeta\in L(S)$ 
and $|F(\zeta)|\le M_1$ for $\zeta\in R(S)$. By Hadamard's three lines theorem, 
we must have
$$|F(\theta)|\le M_0^{1-\theta}M_1^\theta,$$
or 
$$(1-|z|^2)^{k-\gamma}|R^kf(z)|\le M_0^{1-\theta}M_1^\theta.$$
Since the constant on the right-hand side is independent of $z$, we have
shown that $f\in \Lambda_\gamma$. Therefore,
$$\left[\Lambda_\alpha,\Lambda_\beta\right]_\theta\subset\Lambda_\gamma,$$
and the inclusion is continuous. This completes the proof of the theorem.
\end{proof}

\section{Reproducing Kernels}

In this section we focus on the Hilbert space case $p=2$. We are going to
obtain a characterization of $A^2_\alpha$ in terms of Taylor coefficients, and
we are going to define a canonical inner product on $A^2_\alpha$ so that
the associated reproducing kernel can be calculated in closed form.

Reproducing kernels for $A^2_\alpha$ are also calculated in 
Beatrous-Buebea \cite{BB2} in terms of a certain family of hypergeometric functions.
Our approach here is different. We wish to write the reproducing kernel of $A^2_\alpha$
as something that is as close to $(1-\langle z,w\rangle)^{-(n+1+\alpha)}$ as possible.

\begin{thm}
Suppose $\alpha$ is real and
$$f(z)=\sum_ma_mz^m.$$
Then $f\in A^2_\alpha$ if and only if its Taylor coefficients satisfy the
condition
\begin{equation}
\sum_{|m|>0}\frac{m!\,e^{|m|}}{|m|^{n+|m|+\alpha+\frac12}}|a_m|^2<\infty.
\label{eq18}
\end{equation}
\label{40}
\end{thm}

\begin{proof}
Fix a positive integer $k$ such that $2k+\alpha>-1$. If $f(z)=\sum_ma_mz^m$
is the Taylor series of $f$ in $\bn$, then
$$R^kf(z)=\sum_{|m|>0}a_m|m|^kz^m.$$
It follows that the integral
$$I_{k,\alpha}(f)=\inb\left|(1-|z|^2)^kR^kf(z)\right|^2\,dv_\alpha(z)$$
is equal to
$$\sum_{|m|>0}|a_m|^2|m|^{2k}\inb|z^m|^2(1-|z|^2)^{2k+\alpha}\,dv(z).$$
By integration in polar coordinates (see 1.4.3 of Rudin \cite{rudin} or Lemma 1.11 
of Zhu \cite{zhu6}), there exists a constant $C>0$ (independent of $f$)
such that
$$I_{k,\alpha}(f)=C\sum_{|m|>0}\frac{|m|^{2k}m!}{\Gamma(n+|m|+2k+\alpha+1)}\,|a_m|^2.$$
Since $n$, $k$, and $\alpha$ are all constants, an application of Stirling's formula
shows that
$$\Gamma(n+|m|+2k+\alpha+1)\sim|m|^{n+|m|+2k+\alpha+\frac12}e^{-|m|}$$
as $|m|\to\infty$. We conclude that the integral $I_{k,\alpha}(f)$ is finite if and only
if the condition in (\ref{eq18}) holds, and the proof of the theorem is complete.
\end{proof}

An immediate consequence of the condition in (\ref{eq18}) is that the space $A^2_\alpha$
is independent of the integer $k$ used in the definition of $A^2_\alpha$. Of course,
we already knew this from Section 3.

\begin{thm}
Suppose $\alpha>-(n+1)$. Then $A^2_\alpha$ can be equipped with an inner
product such that the associated reproducing kernel is given by
\begin{equation}
K_\alpha(z,w)=\frac1{(1-\langle z,w\rangle)^{n+1+\alpha}}.
\label{eq19}
\end{equation}
\label{41}
\end{thm}

\begin{proof}
It follows from Stirling's formula again that the condition in (\ref{eq18}) is equivalent to
\begin{equation}
\sum_m\frac{m!\,\Gamma(n+1+\alpha)}{\Gamma(n+|m|+\alpha+1)}\,|a_m|^2<\infty.
\label{eq20}
\end{equation}
Now define an inner product on $A^2_\alpha$ as follows:
$$\langle f,g\rangle_\alpha=\sum_m\frac{m!\,\Gamma(n+1+\alpha)}{\Gamma(n+|m|+
\alpha+1)}\,a_m\overline b_m,$$
where
$$f(z)=\sum_ma_mz^m,\qquad g(z)=\sum_mb_mz^m.$$
Then $A^2_\alpha$ becomes a separable Hilbert space with the following functions
forming an orthonormal basis:
$$e_m(z)=\sqrt{\frac{\Gamma(n+|m|+\alpha+1)}{m!\,\Gamma(n+\alpha+1)}}\,z^m,$$
where $m$ runs over all $n$-tuples of nonnegative integers. It follows from the
multinomial formula (\ref{eq1}) that the reproducing kernel of $A^2_\alpha$ is given by
\begin{eqnarray*}
K_\alpha(z,w)&=&\sum_me_m(z)\overline{e_m(w)}\\
&=&\sum_m\frac{\Gamma(n+|m|+\alpha+1)}{m!\,\Gamma(n+\alpha+1)}
\,z^m\overline w^m\\
&=&\sum_{k=0}^\infty\frac{\Gamma(n+k+\alpha+1)}{k!\,\Gamma(n+1+\alpha)}
\sum_{|m|=k}\frac{k!}{m!}z^m\overline w^m\\
&=&\sum_{k=0}^\infty\frac{\Gamma(n+k+1+\alpha)}{k!\,\Gamma(n+\alpha+1)}
\langle z,w\rangle^k\\
&=&\frac1{(1-\langle z,w\rangle)^{n+1+\alpha}}.
\end{eqnarray*}
This proves the desired result.
\end{proof}

When $\alpha>-1$, the reproducing kernel for $A^2_\alpha$ is of course well
known. See Rudin \cite{rudin} or Zhu \cite{zhu6}. When $\alpha\le -1$, the point
here is that you need to use an appropriate inner product on $A^2_\alpha$ so
that its reproducing kernel is computable in closed form.

\begin{thm}
Suppose $\alpha=-(n+1)$. Then $A^2_\alpha$ can be equipped with an inner product
such that the associated reproducing kernel is
\begin{equation}
K_{-(n+1)}(z,w)=1+\log\frac1{1-\langle z,w\rangle}.
\label{eq21}
\end{equation}
\label{42}
\end{thm}

\begin{proof}
If $\alpha=-(n+1)$, then Theorem~\ref{40} tells us that a holomorphic function
$f(z)=\sum_ma_mz^m$ in $\bn$ belongs to $A^2_{-(n+1)}$ if and only if
$$\sum_{|m|>0}\frac{m!e^{|m|}}{|m|^{|m|-\frac12}}\,|a_m|^2<\infty,$$
which, according to Stirling's formula, is equivalent to
$$\sum_{|m|>0}|m|\frac{m!}{|m|!}\,|a_m|^2<\infty.$$
If we define an inner product on $A^2_{-(n+1)}$ by
\begin{equation}
\langle f,g\rangle_{-(n+1)}=f(0)\overline{g(0)}+\sum_{|m|>0}|m|\frac{m!}{|m|!}\,
a_m\overline b_m,
\label{eq22}
\end{equation}
where
$$f(z)=\sum_ma_mz^m,\qquad g(z)=\sum_mb_mz^m,$$
then $A^2_{-(n+1)}$ becomes a separable Hilbert space with the following
functions forming an orthonormal basis:
$$1,\quad e_m(z)=\sqrt{\frac{|m|!}{m!|m|}}\,z^m,$$
where $m$ runs over all $n$-tuples of nonnegative integers with $|m|>0$.
It follows from (\ref{eq1}) that the reproducing kernel of $A^2_{-(n+1)}$ is given by
\begin{eqnarray*}
K_{-(n+1)}(z,w)&=&1+\sum_{|m|>0}\frac{|m|!}{m!|m|}\,z^m\overline w^m\\
&=&1+\sum_{k=1}^\infty\frac1k\sum_{|m|=k}\frac{k!}{m!}z^m\overline w^m\\
&=&1+\sum_{k=1}^\infty\frac{\langle z,w\rangle^k}k\\
&=&1+\log\frac1{1-\langle z,w\rangle},
\end{eqnarray*}
completing the proof of the theorem.
\end{proof}

The space $A^2_{-(n+1)}$ can be thought of as the high dimensional analog of
the classical Dirichlet space in the unit disk. It is the unique space of holomorphic
functions in the unit ball that can be equipped with a semi-inner product that is
invariant under the action of the automorphism group. See Zhu \cite{zhu2}.
The formula in Theorem~\ref{42} above also appeared in Peloso \cite{Peloso}
and Zhu \cite{zhu2}.

\begin{thm}
Suppose $-N<n+1+\alpha<-N+1$ for some positive integer $N$. Then for any
polynomial
$$Q(z,w)=\sum_{|m|\le N}\omega_mz^m\overline w^m$$
with the property that
$$\omega_m>(-1)^{N+1}\frac{\Gamma(n+|m|+\alpha+1)}{m!\,\Gamma(n+\alpha+1)},$$
we can equip $A^2_\alpha$ with an inner product such that the associated reproducing
kernel is given by
\begin{equation}
K_\alpha(z,w)=Q(z,w)+\frac{(-1)^N}{(1-\langle z,w\rangle)^{n+1+\alpha}}.
\label{eq23}
\end{equation}
\label{43}
\end{thm}

\begin{proof}
By Theorem~\ref{40} and Stirling's formula again, a function $f(z)=\sum_ma_mz^m$
is in $A^2_\alpha$ if and only if
$$\sum_{|m|>0}\left|\frac{m!\,\Gamma(n+1+\alpha)}{\Gamma(n+|m|+\alpha+1)}
\right|\,|a_m|^2<\infty.$$

If $-N<n+1+\alpha<-N+1$, it follows from the identity
$$\frac{\Gamma(n+\alpha+1)}{\Gamma(n+|m|+\alpha+1)}=\frac1{(n+1+\alpha)
(n+2+\alpha)\cdots(n+|m|+\alpha)}$$
that for any $|m|>N$ we have
$$\left|\frac{\Gamma(n+\alpha+1)}{\Gamma(n+|m|+\alpha+1)}\right|=(-1)^N
\frac{\Gamma(n+\alpha+1)}{\Gamma(n+|m|+\alpha+1)}.$$
Therefore, for any positive coefficients $c_m$, where $|m|\le N$, we can define an
inner product on $A^2_\alpha$ as follows:
\begin{equation}
\langle f,g\rangle_\alpha=\sum_{|m|\le N}c_ma_m\overline b_m+(-1)^N
\sum_{|m|>N}\frac{m!\,\Gamma(n+\alpha+1)}{\Gamma(n+|m|+\alpha+1)}
\,a_m\overline b_m,
\label{eq24}
\end{equation}
where
$$f(z)=\sum_ma_mz^m,\qquad g(z)=\sum_mb_mz^m.$$
Then $A^2_\alpha$ becomes a separable Hilbert space with the following
functions forming an orthonormal basis:
$$e_m(z)=\frac1{\sqrt{c_m}}\,z^m,\qquad |m|\le N,$$
and
$$e_m(z)=\sqrt{(-1)^N\frac{\Gamma(n+|m|+\alpha+1)}{m!\,\Gamma(n+
\alpha+1)}}\,z^m,\qquad |m|>N.$$
Using the mutinomial formula (\ref{eq1}), we find that the corresponding reproducing 
kernel for $A^2_\alpha$ is given by
\begin{eqnarray*}
K_\alpha(z,w)&=&\sum_{|m|\le N}e_m(z)\overline{e_m(w)}+
\sum_{|m|>N}e_m(z)\overline{e_m(w)}\\
&=&\sum_{|m|\le N}\frac1{c_m}z^m\overline w^m+(-1)^N\sum_{|m|>N}
\frac{\Gamma(n+|m|+\alpha+1)}{m!\,\Gamma(n+\alpha+1)}\,z^m\overline w^m\\
&=&\sum_{|m|\le N}\omega_mz^m\overline w^m+\frac{(-1)^N}{(1-\langle z,
w\rangle)^{n+1+\alpha}},
\end{eqnarray*}
where
$$\omega_m=\frac1{c_m}-(-1)^N\frac{\Gamma(n+|m|+\alpha+1)}{m!\,\Gamma(n+
\alpha+1)}>(-1)^{N+1}\frac{\Gamma(n+|m|+\alpha+1)}{m!\,\Gamma(n+\alpha+1)}.$$
This completes the proof of the theorem.
\end{proof}

The appearance of the sign $(-1)^N$ in (\ref{eq23}) is a little peculiar; we do not
know if there is any simple explanation for it. We also note in passing that the
reproducing kernel given by (\ref{eq23}) is bounded.

It remains for us to consider the case in which $n+1+\alpha=-N$ is a negative integer.
The principal part of the reproducing kernel in this case will be shown to be
the function
$$(\langle z,w\rangle-1)^N\log\frac1{1-\langle z,w\rangle}.$$

Thus for every positive integer $N$ we consider the function
$$f_N(z)=(z-1)^N\log\frac1{1-z},\qquad z\in\D.$$
It is clear that each $f_N$ is analytic in the unit disk $\D$ and
$$f'_{N+1}(z)=(N+1)f_N(z)-(z-1)^N.$$
In particular,
$$f^{(k+1)}_{N+1}(z)=(N+1)f^{(k)}_N(z),\qquad k>N.$$
It follows from this and induction that $f^{(k)}_N(0)>0$ for all $k>N$.
Also observe that the $N$th derivative of $f_N$ is $-\log(1-z)$ plus a polynomial, 
so the Taylor coefficients of $f_N$ has the property that
$$\frac{f_N^{(k)}(0)}{k!}\sim\frac1{k^{N+1}}$$
as $k\to\infty$.

\begin{thm}
Suppose $n+1+\alpha=-N$ for some positive integer $N$ and
$$(z-1)^N\log\frac1{1-z}=\sum_{k=0}^\infty A_kz^k.$$
Then for any polynomial
$$Q(z,w)=\sum_{|m|\le N}\omega_mz^m\overline w^m$$
with the property that
$$\omega_m>-\frac{|m|!\,A_{|m|}}{m!},$$
we can equip $A^2_\alpha$ with an inner product such that the associated 
reproducing kernel is
\begin{equation}
K_\alpha(z,w)=Q(z,w)+(\langle z,w\rangle-1)^N\log\frac1{1-\langle z,w\rangle}.
\label{eq25}
\end{equation}
\label{44}
\end{thm}

\begin{proof}
It follows from (\ref{eq18}) and Stirling's formula that a function $f(z)=\sum_m
a_mz^m$ belongs to $A^2_\alpha$ if and only if
$$\sum_m|m|^{N+1}\frac{m!}{|m|!}\,|a_m|^2<\infty,$$
which is equivalent to
$$\sum_{|m|>N}\frac{m!}{|m|!\,A_{|m|}}|a_m|^2<\infty.$$
If $c_m>0$ for $|m|\le N$, we can define an inner product on $A^2_\alpha$ as follows:
\begin{equation}
\langle f,g\rangle_\alpha=\sum_{|m|\le N}c_ma_m\overline b_m+
\sum_{|m|>N}\frac{m!}{|m|!\,A_{|m|}}\,a_m\overline b_m.
\label{eq26}
\end{equation}
where
$$f(z)=\sum_ma_mz^m,\qquad g(z)=\sum_mb_mz^m.$$
Then $A^2_\alpha$ becomes a separable Hilbert
space and the following functions form an orthonormal basis:
$$e_m(z)=\frac1{\sqrt{c_m}}\,z^m,\qquad |m|\le N,$$
and
$$e_m(z)=\sqrt{\frac{|m|!\,A_{|m|}}{m!}}\,z^m,\qquad |m|>N.$$
The associated reproducing kernel for $A^2_\alpha$ is given by
\begin{eqnarray*}
K_\alpha(z,w)&=&\sum_{|m|\le N}e_m(z)\overline{e_m(w)}+\sum_{|m|>N}
e_m(z)\overline{e_m(w)}\\
&=&\sum_{|m|\le N}\frac{z^m\overline w^m}{c_m}+\sum_{|m|>N}
\frac{|m|!\,A_{|m|}}{m!}z^m\overline w^m\\
&=&\sum_{|m|\le N}\left(\frac1{c_m}-\frac{|m|!\,A_{|m|}}{m!}\right)z^m\overline w^m+
\sum_{k=0}^\infty A_k\sum_{|m|=k}\frac{k!}{m!}z^m\overline w^m\\
&=&\sum_{|m|\le N}\omega_mz^m\overline w^m+\sum_{k=0}^\infty A_k
\langle z,w\rangle^k\\
&=&Q(z,w)+(\langle z,w\rangle-1)^N\log\frac1{1-\langle z,w\rangle},
\end{eqnarray*}
where the coefficients of 
$$Q(z,w)=\sum_{|m|\le N}\omega_mz^m\overline w^m$$
satisfy
$$\omega_m=\frac1{c_m}-\frac{|m|!\,A_{|m|}}{m!}>-\frac{|m|!\,A_{|m|}}{m!}.$$
This completes the proof of the theorem.
\end{proof}

Once again, the reproducing kernel in (\ref{eq25}) is bounded on $\bn\times\bn$.
Also notice that we can rewrite the kernel in (\ref{eq25}) as
$$K_\alpha(z,w)=Q(z,w)+(-1)^N(1-\langle z,w\rangle)^N\log\frac1{1-\langle z,
w\rangle},$$
which is probably a partial explanation for the sign $(-1)^N$ in (\ref{eq23}).

It is clear that the reproducing kernel of a Hilbert space of holomorphic functions
depends on the inner product used for the space. We close this section by examining
the reproducing kernel of $A^2_{\alpha}$ that corresponds to the following natural 
inner product which we have used in section 7:
$$
\langle f,g\rangle_\alpha
=f(0)\overline{g(0)}+\inb R^kf(z)\overline{R^k g(z)}\,dv_{2k+\alpha}(z),
$$
where $k$ is any nonnegative integer with $2k+\alpha>-1$.
This inner product gives rise to the norm
$$
\|f\|_{2,\alpha}=\left(|f(0)|^2+\inb |R^kf(z)|^2\,dv_{2k+\alpha}(z)\right)^{1/2}
$$
for $f\in A^2_{\alpha}$. For this inner product we can show that the corresponding
reproducing kernel for $A^2_{\alpha}$ is
$$
K_w^{\alpha}(z)=
K_\alpha(z,w)=
1+R^{-2k}\left(\frac1{(1-\langle z,w\rangle)^{n+1+\alpha+2k}}\right).
$$
The result is a simple consequence of the identity,
$$
f(0)\overline{g(0)}+
\inb R^kf(z) \overline{R^{-k}g(z)}\,dv_{2k+\alpha}(z)
=\inb f(z) \overline{g(z)}\,dv_{2k+\alpha}(z),
$$
which can easily be proved by the use of Taylor expansions. 
We leave the details to the interested reader.

\section{Carleson Type Measures}

The purpose of this section is to study Carleson type measures for the
Bergman spaces $\apa$. Unlike most other sections of the paper, the results
here depend very much on the various parameters.

The notion of Carleson measures was of course introduced by Carleson
\cite{carleson1}\cite{carleson2} for the unit disk. Carleson's original definition works 
well in the theory of Hardy spaces, and this can easily be seen in such classics
as Duren \cite{duren} and Garnett \cite{garnett}. The characterization of
Carleson measures for the Hardy spaces of the unit ball can be found in
H\"ormander \cite{hormander} and Power \cite{Power}. 

Later, the notion of Carleson measures was extended to the context of Bergman 
spaces. Earlier papers in this direction include Cima-Wogen \cite{CW}, 
Hastings \cite{hastings}, Luecking \cite{luecking}, Zhu \cite{zhu1}. Also, Carleson 
type measures have been studied for holomorphic Besov spaces (of which the 
Dirichlet space is a special case); see Arcozzi-Rochberg-Sawyer \cite{ars2}, 
Kaptanoglu \cite{Kap3}, Stegenga \cite{stegenga}, and Wu \cite{wu}. In particular, 
our results of this section contain several special cases that have been known before.

For any $\zeta\in\sn$ and $r>0$ let
$$Q_r(\zeta)=\{z\in\bn:|1-\langle z,\zeta\rangle|<r\}.$$
These are the high dimensional analogues of Carleson squares in the unit disk.
They are also called nonisotropic metric balls. See Rudin \cite{rudin} or
Zhu \cite{zhu6} for more information about the geometry of these nonisotropic
balls.

\begin{thm}
Suppose $n+1+\alpha>0$ and $\mu$ is a positive Borel measure on $\bn$.
Then the following conditions are equivalent.
\begin{enumerate}
\item[(a)] There exists a constant $C>0$ such that
\begin{equation}
\mu(Q_r(\zeta))\le Cr^{n+1+\alpha}
\label{eq27}
\end{equation}
for all $\zeta\in\sn$ and all $r>0$.
\item[(b)] For each $s>0$ there exists a constant $C>0$ such that
\begin{equation}
\inb\frac{(1-|z|^2)^s\,d\mu(w)}{|1-\langle z,w\rangle|^{n+1+\alpha+s}}\le C
\label{eq28}
\end{equation}
for all $z\in\bn$.
\item[(c)] For some $s>0$ there exists a constant $C>0$ such that
the inequality in (\ref{eq28}) holds for all $z\in\bn$.
\end{enumerate}
\label{45}
\end{thm}

\begin{proof}
It is obvious that condition (b) implies (c). 

Now assume that condition (c) holds, that is, there exist positive constants
$s$ and $C$ such that the inequality in (\ref{eq28}) holds. If $\zeta\in\sn$ and
$r\in(0,1)$, then
\begin{equation}
\int_{Q_r(\zeta)}\frac{(1-|z|^2)^s\,d\mu(w)}{|1-\langle z,w\rangle|^{n+1+\alpha+s}}\le C
\label{eq29}
\end{equation}
for all $z\in\bn$. If we choose $z=(1-r)\zeta$, then
$$1-\langle z,w\rangle=(1-r)(1-\langle\zeta,w\rangle)+r$$
for all $w\in\bn$, so
$$|1-\langle z,w\rangle|\le(1-r)r+r<2r$$
for all $w\in Q_r(\zeta)$, which gives
$$\frac{(1-|z|^2)^s}{|1-\langle z,w\rangle|^{n+1+\alpha+s}}\ge\frac{r^s}
{(2r)^{n+1+\alpha+s}}=\frac{2^{-(n+1+\alpha+s)}}{r^{n+1+\alpha}}$$
for all $w\in Q_r(\zeta)$. Combining this with (\ref{eq29}), we conclude that
$$\mu(Q_r(\zeta))\le 2^{n+1+\alpha+s}Cr^{n+1+\alpha}$$
for all $\zeta\in\sn$ and all $r\in(0,1)$. The case $r\ge1$ can be disposed of very
easily. This proves that condition (c) implies (a).

Next assume that condition (a) holds. In particular, $\mu$ is a finite measure, so
$$\sup_{|z|\le3/4}\inb\frac{(1-|z|^2)^s\,d\mu(w)}{|1-\langle z,w\rangle|^{n+1+\alpha+s}}
<\infty$$
for each $s>0$. We fix an arbitrary positive number $s$ and proceed to show that 
the inequality in (\ref{eq28}) must hold for $3/4<|z|<1$.

Fix some point $z\in\bn$ with $3/4<|z|<1$ and choose $\zeta=z/|z|$. For any
nonnegative integer $k$ let $r_k=2^{k+1}(1-|z|)$. We decompose the unit ball $\bn$
into the disjoint union of the following sets:
$$E_0=Q_{r_0}(\zeta),\quad E_k=Q_{r_k}(\zeta)-Q_{r_{k-1}}(\zeta),\quad 1\le k<\infty.$$
By condition (a), we have
$$\mu(E_k)\le\mu(Q_{r_k}(\zeta))\le 2^{(k+1)(n+1+\alpha)}(1-|z|)^{n+1+\alpha}C$$
for all $k\ge0$. On the other hand, if $k\ge1$ and $w\in E_k$, then
\begin{eqnarray*}
|1-\langle z,w\rangle|&=&\bigl|(1-|z|)+|z|(1-\langle\zeta,w\rangle)\bigr|\\
&\ge&|z||1-\langle\zeta,w\rangle|-(1-|z|)\\
&\ge&(3/4)2^k(1-|z|)-(1-|z|)\\
&\ge&2^{k-1}(1-|z|).
\end{eqnarray*}
This holds for $k=0$ as well, because
$$|1-\langle z,w\rangle|\ge1-|z|\ge\frac12(1-|z|).$$
It follows that
\begin{eqnarray*}
\inb\frac{(1-|z|^2)^s\,d\mu(w)}{|1-\langle z,w\rangle|^{n+1+\alpha+s}}&=&
\sum_{k=0}^\infty\int_{E_k}\frac{(1-|z|^2)^s\,d\mu(w)}{|1-\langle z,w\rangle|^{n+1
+\alpha+s}}\\
&\le&\sum_{k=0}^\infty\frac{(1-|z|^2)^s\mu(E_k)}{(2^{k-1}(1-|z|))^{n+1+\alpha+s}}\\
&\le&\sum_{k=0}^\infty\frac{2^{s+(k+1)(n+1+\alpha)}(1-|z|)^{n+1+\alpha+s}C}
{2^{(k-1)(n+1+\alpha+s)}(1-|z|)^{n+1+\alpha+s}}\\
&=&C'\sum_{k=0}^\infty\frac1{(2^s)^k}<\infty,
\end{eqnarray*}
where $C'$ is a positive constant independent of $z$. This completes the proof 
of the theorem.
\end{proof}

Our results are most complete when $0<p\le1$. The following result settles the case
$n+1+\alpha>0$, and Proposition~\ref{49} deals with the cases $n+1+\alpha\le0$.

\begin{thm}
Suppose $\alpha>-(n+1)$, $0<p\le1$, and $\mu$ is a positive Borel measure on 
$\bn$. Then the following two conditions are equivalent.
\begin{enumerate}
\item[(a)] There exists a constant $C>0$ such that
\begin{equation}
\inb|f(w)|^p\,d\mu(w)\le C\|f\|^p_{p,\alpha}
\label{eq30}
\end{equation}
for all $f\in\apa$.
\item[(b)] There exists a constant $C>0$ such that
\begin{equation}
\mu(Q_r(\zeta))\le Cr^{n+1+\alpha}
\label{eq31}
\end{equation}
for all $\zeta\in\sn$ and all $r\in(0,1)$.
\end{enumerate}
\label{46}
\end{thm}

\begin{proof}
First assume that condition (a) holds. We consider the function
$$f(w)=\frac1{(1-\langle w,z\rangle)^{(n+1+\alpha+s)/p}},\qquad w\in\bn,$$
where $s$ is positive and $z\in\bn$. If $k$ is the smallest nonnegative integer
such that $kp+\alpha>-1$, then an elementary calculation shows that
$$R^kf(w)=\frac{Q(\langle w,z\rangle)}{(1-\langle w,z\rangle)^{k+(n+1+\alpha+s)/p}},$$
where $Q$ is a polynomial of degree $k$. It follows from Proposition~\ref{7}
that there exists a constant $C>0$ (independent of $z$) such that
$$\inb\bigl|(1-|w|^2)^kR^kf(w)\bigr|^p\,dv_\alpha(w)\le\frac C{(1-|z|^2)^s}$$
for all $z\in\bn$. Combining this with condition (a), we conclude that
$$\inb\frac{(1-|z|^2)^s\,d\mu(w)}{|1-\langle z,w\rangle|^{n+1+\alpha+s}}\le C$$
for all $z\in\bn$, which, according to Theorem~\ref{45}, is equivalent to condition (b).

Next assume that condition (b) holds. Then condition (b) of Theorem~\ref{45} holds.
We proceed to prove the inequality in (\ref{eq30}). 

Given $f\in\apa$, we use the atomic decomposition for $\apa$ (see
Theorem~\ref{32}) to write
$$f(z)=\sum_{k=1}^\infty c_k\frac{(1-|a_k|^2)^{b-(n+1+\alpha)/p}}{(1-\langle z,a_k
\rangle)^b},$$
where $b$ is a sufficiently large positive number and
$$\sum_{k=1}^\infty|c_k|^p\le C\|f\|^p_{p,\alpha}$$
for some positive constant $C$ independent of $f$. Since $0<p\le1$, we have
$$|f(z)|^p\le\sum_{k=1}^\infty|c_k|^p\frac{(1-|a_k|^2)^{pb-(n+1+\alpha)}}
{|1-\langle z,a_k\rangle|^{pb}},$$
and so
$$\inb|f(z)|^p\,d\mu(z)\le\sum_{k=1}^\infty|c_k|^p(1-|a_k|^2)^{pb-(n+1+\alpha)}
\inb\frac{d\mu(z)}{|1-\langle z,a_k\rangle|^{pb}}.$$
Apply condition (b) of Theorem~\ref{45} to the last integral above. We obtain
a constant $C'>0$ (independent of $f$) such that
$$\inb|f(z)|^p\,d\mu(z)\le C'\sum_{k=1}^\infty|c_k|^p\le CC'\|f\|^p_{p,\alpha}.$$
This completes the proof of the theorem.
\end{proof}

\begin{cor}
If $\alpha>-1$ and $p>0$, then the following two conditions are equivalent for
a positive Borel measure $\mu$ on $\bn$.
\begin{enumerate}
\item[(a)] There exists a constant $C>0$ such that
$$\inb|f(z)|^p\,d\mu(z)\le C\inb|f(z)|^p\,dv_\alpha(z)$$
for all $f\in\apa$.
\item[(b)] There exists a constant $C>0$ such that
$$\mu(Q_r(\zeta))\le Cr^{n+1+\alpha}$$
for all $r>0$ and $\zeta\in\sn$.
\end{enumerate}
\label{47}
\end{cor}

\begin{proof}
That (a) implies (b) follows from the first part of the proof of Theorem~\ref{46}.
Theorem~\ref{46} also tells us that (b) implies (a) when $0<p\le1$.

Now assume that condition (b) holds and $f\in\apa$ for some $p>1$. Then
the function $g=f^N$ belongs to $A^q_\alpha$, where $N$ is any positive
integer and $q=p/N$. We choose $N$ large enough so that $0<q<1$. Then
$$\inb|g(z)|^q\,d\mu(z)\le C\inb|g(z)|^q\,dv_\alpha(z),$$
where $C$ is a positive constant independent of $g$. This is the same as
$$\inb|f(z)|^p\,d\mu(z)\le C\inb|f(z)|^p\,dv_\alpha(z),$$
and the corollary is proved.
\end{proof}

Let $\beta(z,w)$ be the distance between $z$ and $w$ in the Bergman metric
of $\bn$. For any $R>0$ and $a\in\bn$ we write
$$D(a,R)=\{z\in\bn:\beta(z,a)<R\}.$$
When $\alpha>-1$, the condition
$$\mu(Q_r(\zeta))\le Cr^{n+1+\alpha},\qquad r>0,\zeta\in\sn,$$
is equivalent to the condition
$$\mu(D(a,R))\le C_R(1-|a|^2)^{n+1+\alpha},\qquad a\in\bn.$$
See Lemma 5.23 and Corollary 5.24 of \cite{zhu6} (Note that the definition of
$Q_r(\zeta)$ in \cite{zhu6} is different from its definition in this paper). It can be 
shown that these two conditions are no longer equivalent when $\alpha\le-1$.
In fact, if $f$ is a function in the Bloch space that is not in BMOA, then the
measure
$$d\mu(z)=|Rf(z)|^2(1-|z|^2)\,dv(z)$$
satisfies
$$\sup_{r,\zeta}\frac{\mu(Q_r(\zeta))}{r^n}=\infty$$
and
$$\sup_{a\in\bn}\frac{\mu(D(a,R))}{(1-|a|^2)^n}<\infty.$$

Recall that the Hardy space $H^p$, where $0<p<\infty$, consists of holomorphic
functions $f$ in $\bn$ such that
$$\|f\|^p_p=\sup_{0<r<1}\ins|f(r\zeta)|^p\,d\sigma(\zeta)<\infty,$$
where $d\sigma$ is the normalized surface area measure on $\sn$.
It is well known that every function $f\in H^p$ has a finite radial limit at almost 
every point on $\sn$. If we write
$$f(\zeta)=\lim_{r\to1^-}f(r\zeta),\qquad \zeta\in\sn,$$
then we actually have
$$\|f\|^p_p=\ins|f(\zeta)|^p\,d\sigma(\zeta).$$
It is known that the following two conditions are equivalent for a positive
Borel measure $\mu$ on $\bn$; see H\"ormander \cite{hormander},
Power \cite{Power}, or Zhu \cite{zhu6}.
\begin{enumerate}
\item[(a)] There exists a constant $C>0$ such that
$$\inb|f(z)|^p\,d\mu(z)\le C\ins|f(\zeta)|^p\,d\sigma(\zeta)$$
for all $f\in H^p$.
\item[(b)] There exists a constant $C>0$ such that
$$\mu(Q_r(\zeta))\le Cr^n$$
for all $r>0$ and $\zeta\in\sn$.
\end{enumerate}

\begin{cor}
Suppose $\alpha=-1$, $0<p\le2$, and $\mu$ is a positive Borel measure on $\bn$.
Then the following two conditions are equivalent.
\begin{enumerate}
\item[(a)] There exists a constant $C>0$ such that
$$\inb|f(z)|^p\,d\mu(z)\le C\|f\|^p_{p,\alpha}$$
for all $f\in\apa$.
\item[(b)] There exists a constant $C>0$ such that
$$\mu(Q_r(\zeta))\le Cr^n$$
for all $r>0$ and $\zeta\in\sn$.
\end{enumerate}
\label{48}
\end{cor}

\begin{proof}
That (a) implies (b) follows from the first part of the proof of Theorem~\ref{46}.

To show that condition (b) implies (a), we notice that $A^2_{-1}=H^2$, so the case 
$p=2$ follows from the characterization of Carleson measures for Hardy spaces. 
The case $0<p\le1$ follows from Theorem~\ref{46}. The case of $1\le p\le2$ then 
follows from complex interpolation.
\end{proof}

\begin{prop}
Let $\mu$ be a positive Borel measure on $\bn$. If $n+1+\alpha<0$ and
$0<p<\infty$, or if $n+1+\alpha=0$ and $0<p\le1$, then
the following two conditions are equivalent.
\begin{enumerate}
\item[(a)] There exists a constant $C>0$ such that
$$\inb|f(z)|^p\,d\mu(z)\le C\|f\|^p_{p,\alpha}$$
for all $f\in\apa$.
\item[(b)] The measure $\mu$ is finite.
\end{enumerate}
\label{49}
\end{prop}

\begin{proof}
Since $\apa$ contains all constant functions, it is clear that condition (a)
implies (b).

On the other hand, if $\mu$ is a finite positive Borel measure, it follows from 
Theorems~\ref{21} and \ref{22} that $\apa$ is contained in $L^p(\bn,d\mu)$. By 
the closed graph theorem, $\apa$ is continuously contained in $L^p(\bn,d\mu)$, 
so there exists a positive constant $C>0$ such that
$$\inb|f|^p\,d\mu\le C\|f\|^p_{p,\alpha}$$
for all $f\in\apa$.
\end{proof}

As far as the condition
$$\inb|f(z)|^p\,d\mu(z)\le C\|f\|^p_{p,\alpha},\qquad f\in\apa,$$
is concerned, the most difficult case is probably when $\alpha=-(n+1)$ and $1<p<\infty$. 
This case is considered in Arcozzi-Rochberg-Sawyer \cite{ars2} and complete results are 
obtained in the range $1<p<2+1/(n-1)$. Earlier results dealing with the Besov spaces 
include Arcozzi \cite{arcozzi}, Arcozzi-Rochberg-Sawyer \cite{ars1},
Stegenga \cite{stegenga}, and Wu \cite{wu}.

\begin{thm}
Suppose $0<p\le q<\infty$, $\alpha$ is real, and $\mu$ is a positive Borel measure 
on $\bn$. Then for any nonnegative integer $k$ with $\alpha+kp>-1$
the following conditions are equivalent.
\begin{enumerate}
\item[(a)] There is a contant $C>0$ such that
$$\inb|R^kf(w)|^q\,d\mu(w)\le C\|f\|_{p,\alpha}^q$$
for all $f\in A^p_{\alpha}$.
\item[(b)] For each (or some) $s>0$ there is a constant $C>0$ such that
$$\inb\frac{(1-|z|^2)^s}{|1-\langle z,w\rangle|^{s+(n+1+\alpha+kp)q/p}}\,d\mu(w)\le C$$
for all $z\in\bn$.
\item[(c)] There is a constant $C>0$ such that
$$\mu(Q_r(\zeta))\le Cr^{(n+1+\alpha+kp)q/p}$$
for all $r>0$ and $\zeta\in(0,1)$. 
\item[(d)] For each (or some) $R>0$ there exists a constant $C>0$ such that
$$\mu(D(a,R))\le C(1-|a|^2)^{(n+1+\alpha+kp)q/p}$$
for all $a\in\bn$.
\end{enumerate}
\label{50}
\end{thm}

\begin{proof}
First suppose (a) holds. Applying (a) to the functions $f_k(z)=z_k$, $1\le k\le n$,
we see that $\mu$ is a finite measure. For a fixed $z\in\bn$ let 
$$h_z(w)=\frac{(1-|z|^2)^{s/q}}{(1-\langle w,z\rangle)^{s/q+(n+1+\alpha+k)/p}}$$
and let $f_z(w)$ be an analytic function on $\bn$ such that
$$R^kf_z(w)=h_z(w)-h_z(0),$$
where 
$$h_z(0)=(1-|z|^2)^{s/q}\le 1.$$

It follows from Proposition~\ref{7} that 
$$\sup_{z\in\bn}\|f_z\|_{p,\alpha}\le C.$$
Applying (a) to $f_z$ yields
$$\inb|h_z(w)-h_z(0)|^q\,d\mu(w)\le C.$$
It follows from the elementary inequality
$$|h_z(w)|^q\le2^q\left(|h_z(w)-h_z(0)|^q+|h_z(0)|^q\right)$$
that
$$\inb|h_z(w)|^q\,d\mu(w)\le2^q(C+\mu(\bn)),$$
which gives us (b).

Next assume that (b) holds. Recall that $D(z,r)$ is the Bergman metric ball at $z$
with radius $R$. By Lemmas 2.24 and 2.20 of \cite{zhu6}, we have
\begin{eqnarray*}
|R^k f(z)|^p&\le& \frac{C}{(1-|z|^2)^{n+1+\alpha+kp}}\int_{D(z,r)}|R^kf(w)|^p\,dv_{\alpha+kp}(w)\\
&\le& C\int_{D(z,r)}\frac{|R^kf(w)|^p(1-|w|^2)^{sp/q+kp}}
{|1-\langle w,z\rangle|^{sp/q+n+1+\alpha+kp}}\,dv_{\alpha}(w)\\
&\le& C\inb\frac{(1-|w|^2)^{sp/q}\,d\lambda(w)}{|1-\langle z,w\rangle|^{sp/q+n+1+\alpha+kp}},
\end{eqnarray*}
where
$$d\lambda(w)=|R^kf(w)|^p(1-|w|^2)^{kp}\,dv_\alpha(w)$$
is a finite measure on $\bn$ whenever $f\in\apa$. In fact, $\lambda(\bn)\le C\|f\|^p_{p,\alpha}$
for some constant independent of $f$.

If $p=q$, an application of Fubini's theorem to the estimate in the previous paragraph
shows that (b) implies (a). If $p<q$, we write $p'=q/p$ and $1/p'+1/q'=1$, and apply 
H\"older's inequality to the estimate in the previous paragraph. The result is
$$|R^kf(z)|^p\le C\left[\frac{(1-|w|^2)^s\,d\lambda(w)}{|1-\langle z,w\rangle|^{s+(n+1+\alpha
+kp)q/p}}\right]^{1/p'}\left[\lambda(\bn)\right]^{1/q'}.$$
It follows that
$$|R^kf(z)|^q\le C\left(\lambda(\bn)\right)^{\frac{p'}{q'}}\inb\frac{(1-|w|^2)^s\,d\lambda(w)}{
|1-\langle z,w\rangle|^{s+(n+1+\alpha+kp)q/p}}.$$
We now integrate against the measure $d\mu$, apply Fubini's theorem, and use 
condition (b) to obtain
$$\inb|R^kf(z)|^q\,d\mu(z)\le C\left(\lambda(\bn)\right)^{1+\frac{p'}{q'}}.$$
Since $\lambda(\bn)\le C\|f\|^p_{p,\alpha}$, we get
$$\inb|R^kf(z)|^q\,d\mu(z)\le C\|f\|^q_{p,\alpha}.$$
This shows that (b) implies (a).

The equivalence of (b) and (c) has already been proved in Theorem~\ref{45}. Since
$$(n+1+\alpha+kp)q/p\ge n+1+\alpha+kp>n,$$
the equivalence of (c) and (d) follows from the remarks after Corollary~\ref{47}.
This completes the proof of the theorem.
\end{proof}

A similar result can be obtained in terms of fractional radial differential operators
$R^{s,t}$ instead of $R^k$ above. We omit the details.

Once a certain result concerning Carleson measures is established, it is then
relatively easy to formulate and prove its little oh version. For example, with the
same assumptions in Theorem~\ref{50}, we can show that the following four
conditions are equivalent.
\begin{enumerate}
\item[(a)] If $\{f_j\}$ is a bounded sequence in $\apa$ and $f_j(z)\to0$ for every $z\in\bn$,
then
$$\lim_{j\to\infty}\inb|R^kf_j(z)|^q\,d\mu(z)=0.$$
\item[(b)] For every (or some) $s>0$ we have
$$\lim_{|z|\to1^-}\inb\frac{(1-|z|^2)^s\,d\mu(w)}{|1-\langle z,w\rangle|^{s+(n+1+\alpha
+kp)q/p}}=0.$$
\item[(c)] The limit
$$\lim_{r\to0^+}\frac{\mu(Q_r(\zeta))}{r^{(n+1+\alpha+kp)q/p}}=0$$
holds uniformly for $\zeta\in\sn$.
\item[(d)] For every (or some) $R>0$ we have
$$\lim_{|a|\to1^-}\frac{\mu(D(a,R))}{(1-|a|^2)^{(n+1+\alpha+kp)q/p}}=0.$$
\end{enumerate}
The interested reader should have no trouble filling in the details.

As our next theorem shows, the assumption that $p\le q$ is essential for Theorem~\ref{50}. 
To deal with the case $p>q$, we associate two functions to any positive Borel measure
$\mu$ on $\bn$. More specifically, for any real $\gamma$ and $s$ we define
$$B_{s,\gamma}(\mu)(z)=\inb\frac{(1-|z|^2)^s\,d\mu(w)}{|1-\langle z,w\rangle|^{n+1
+s+\gamma}},\qquad z\in\bn,$$
and for any real $\gamma$ and positive $R$ we define
$$\widehat\mu_{R,\gamma}(z)=\frac{\mu(D(z,R))}{(1-|z|^2)^{n+1+\gamma}},
\qquad z\in\bn.$$
If $d\mu(z)=h(z)\,dv_\gamma(z)$, we use the convention that
$$B_{s,\gamma}(h)(z)=B_{s,\gamma}(\mu),\qquad \widehat h_{R,\gamma}(z)=
\widehat\mu_{R,\gamma}(z).$$
It is clear that $B_{s,\gamma}(\mu)(z)$ and $\widehat\mu_{R,\gamma}(z)$ are
certain averages of $\mu$ near the point $z$. The function $B_{s,\gamma}(\mu)$
is sometimes called a Berezin transform of $\mu$.

\begin{lemma}
Suppose $\mu$ is a positive Borel measure on $\bn$. If $\gamma$ is real,
$p>0$, and $R>0$, then there exists a positive constant $C$ such that
$$\inb|g(z)|^p\,d\mu(z)\le C\inb|g(z)|^p\widehat\mu_{R,\gamma}(z)\,dv_\gamma(z)$$
for all holomorphic functions $g$ in $\bn$.
\label{add1}
\end{lemma}

\begin{proof}
It follows from Lemma 2.20 and Corollary 2.21 of \cite{zhu6} that
$$(1-|z|^2)^{n+1}\sim(1-|w|^2)^{n+1}\sim v(D(z,R))\sim v(D(w,R))$$
for $w\in D(z,R)$. We use Lemma 2.24 of \cite{zhu6} and Fubini's theorem to obtain
\begin{eqnarray*}
\inb|g(z)|^p\,d\mu(z)&\le&C_1\inb\frac{d\mu(z)}{(1-|z|^2)^{n+1}}\int_{D(z,R)}|g(w)|^p\,dv(w)\\
&\le&C_2\inb\,d\mu(z)\int_{D(z,R)}\frac{|g(w)|^p\,dv_\gamma(w)}{(1-|w|^2)^{n+1+\gamma}}\\
&=&C_2\inb\frac{|g(w)|^p\,dv_\gamma(w)}{(1-|w|^2)^{n+1+\gamma}}\int_{D(w,R)}\,d\mu(z)\\
&=& C_3\inb|g(w)|^p\widehat\mu_{R,\gamma}(w)\,dv_\gamma(w),
\end{eqnarray*}
which proves the desired estimate.
\end{proof}

\begin{lemma}
Let $\mu$ be a positive Borel measure on $\bn$. If $\gamma$ is real, $s$ is real, and
$R>0$, then there exists a constant $C>0$ such that $B_{s,\gamma}(\mu)
\le CB_{s,\gamma}(\widehat\mu_{R,\gamma})$.
\label{add2}
\end{lemma}

\begin{proof}
For $w\in\bn$, apply Lemma~\ref{add1} to the function
$$g(z)=\frac{(1-|w|^2)^s}{(1-\langle z,w\rangle)^{n+1+s}}$$
with $p=1$. The desired result follows.
\end{proof}

\begin{lemma}
Let $\mu$ be a positive Borel measure on $\bn$. If $\gamma$ and $s$ are real and
$R$ is positive, then there exists a positive constant $C$ such that
$\widehat\mu_{R,\gamma}\le CB_{s,\gamma}(\mu)$.
\label{add3}
\end{lemma}

\begin{proof}
Once again, we have 
$$1-|z|^2\sim1-|w|^2\sim|1-\langle z,w\rangle|$$
for $w\in D(z,R)$. It follows that
\begin{eqnarray*}
\widehat\mu_{R,\gamma}(z)&=&\frac{\mu(D(z,R))}{(1-|z|^2)^{n+1+\gamma}}\\
&\le&
C\int_{D(z,R)}\frac{(1-|z|^2)^s\,d\mu(w)}{|1-\langle z,w\rangle|^{n+1+s+\gamma}}\\
&\le& CB_{s,\gamma}(\mu)(z),
\end{eqnarray*}
proving the desired estimate.
\end{proof}

\begin{thm}
Let $0<q<p<\infty$ and $\alpha$ be any real number, and let $\mu$ be a positive 
Borel measure on $\bn$. Then for any nonnegative integer $k$ with $\alpha+kp>-1$
the following conditions are equivalent.
\begin{enumerate}
\item[(a)] There is a constant $C>0$ such that
$$\inb|R^kf(w)|^q\,d\mu(w)\le C\|f\|_{p,\alpha}^q$$
for all $f\in A^p_{\alpha}$.
\item[(b)] For any bounded sequence $\{f_j\}$ in $A^p_{\alpha}$
with $f_j(z)\to0$ for every $z\in\bn$,
$$\lim_{j\to\infty}\inb|R^kf_j(w)|^q\,d\mu(w)=0.$$
\item[(c)] For any fixed $r>0$ the function $\widehat\mu_{r,\gamma}$
is in $L^{p/(p-q)}(\bn,dv_\gamma)$, where $\gamma=\alpha+kp$.
\item[(d)] For any fixed $s>0$ the function $B_{s,\gamma}(\mu)$ is in
$L^{p/(p-q)}(\bn,dv_\gamma)$, where $\gamma=\alpha+kp$.
\end{enumerate}
\label{51}
\end{thm}

\begin{proof} 
Let $s>0$ satisfy $s+\alpha+kq>-1$. 
It follows from Lemmas 2.24 and 2.20 of \cite{zhu6} that
\begin{eqnarray*}
|R^k f(z)|^q&\le& \frac{C}{(1-|z|^2)^{n+1+s+\alpha+kq}}\int_{D(z,r)}|R^kf(w)|^q\,
dv_{s+\alpha+kq}(w)\\
&\le& C\int_{D(z,r)}\frac{|R^kf(w)|^q(1-|w|^2)^{s+kq}}
{|1-\langle w,z\rangle|^{n+1+s+\alpha+kq}}\,dv_{\alpha}(w)\\
&=& C\inb\frac{|R^kf(w)|^q(1-|w|^2)^{s+kq}}{|1-\langle w,z\rangle|^{n+1+s+\alpha+kq}}
\chi_{D(z,r)}(w)\,dv_{\alpha}(w),
\end{eqnarray*}
where $\chi_E(z)$ denotes the characteristic function of a set $E$.
Integrate with respect to $d\mu$, apply Fubini's theorem, and use Lemma 2.20
of \cite{zhu6}. We see that the integral
$$\inb|R^k f(z)|^q\,d\mu(z)$$
is dominated by
$$\inb\frac{\mu(D(w,r))}{(1-|w|^2)^{n+1+\alpha+kq}}
|R^kf(w)|^q(1-|w|^2)^{kq}\,dv_\alpha(w).$$
If condition (c) holds, then an application of H\"older's inequality yields
\begin{eqnarray*}
\inb|R^k f|^q\,d\mu&\le& C\|f\|_{p,\alpha}^q\left[\inb\left(\frac{\mu(D(w,r))}{(1-
|w|^2)^{n+1+\alpha+kq}}\right)^{\frac p{p-q}}\,dv_{\alpha}\right]^{1-\frac qp}\\
&=& C\|f\|_{p,\alpha}^q\left[\inb\left(\frac{\mu(D(w,r))}{(1-|w|^2)^{n+1+\alpha+kp}}
\right)^{\frac p{p-q}}\,dv_{\alpha+kp}\right]^{1-\frac qp}\\
&\le& C\|f\|_{p,\alpha}^q.
\end{eqnarray*}
This proves that (c) implies (a).

Since $1-|z|\sim1-|w|$ for $z\in D(w,r)$ (see Lemma 2.20 of \cite{zhu6}), there 
exists a constant $\delta>0$ such that
$$\delta^{-1}\le\frac{1-|z|^2}{1-|w|^2}\le\delta$$
for all $z\in D(w,r)$. For $0<t<1$ let 
$$A_t=\{z\in\bn:\, 1-|z|^2<t\}.$$
Then the conditions $z\in A_t$ and $w\in D(z,r)$ imply that $w\in A_{\delta t}$.

Let $\{f_j\}$ be a bounded sequence in $A^p_{\alpha}$ with $f_j(z)\to0$ for every 
$z\in\bn$. Then a normal family argument shows that $f_j(z)\to0$ uniformly on every 
compact subset of $\bn$. Using the estimate from the first paragraph of this proof
together with Fubini's theorem, we see that the integral
$$\int_{A_t}|R^k f_j(z)|^q\,d\mu(z)$$
is dominated by
$$\inb|R^kf_j(w)|^q(1-|w|^2)^{s+kq}\,dv_\alpha(w)
\int_{A_t}\frac{\chi_{D(w,r)}\,d\mu(z)}{|1-\langle z,w\rangle|^{n+1+s+\alpha+kq}}.$$
According to the previous paragraph,
$$\chi_{A_t\cap D(w,r)}(z)=0,\qquad z\in\bn,$$
unless $w\in A_{\delta t}$. It follows that the integral
$$\int_{A_t}|R^k f_j(z)|^q\,d\mu(z)$$
is dominated by
$$\int_{A_{\delta t}}|R^kf_j(w)|^q(1-|w|^2)^{s+kq}\,dv_\alpha(w)
\int_{A_t}\frac{\chi_{D(w,r)}\,d\mu(z)}{|1-\langle z,w\rangle|^{n+1+s+\alpha+kq}}.$$
Since $|1-\langle z,w\rangle|$ is comparable to $1-|w|^2$ whenever $z\in D(w,r)$, we get
$$\int_{A_t}|R^kf_j|^q\,d\mu\le C\int_{A_{\delta t}}\frac{\mu(D(w,r))}{(1-|w|^2)^{n+1+\alpha
+kq}}|R^kf_j(w)|^q(1-|w|^2)^{kq}\,dv_\alpha.$$
By H\"older's inequlity,
$$\int_{A_t}|R^k f_j|^q\,d\mu\le C\|f_j\|_{p,\alpha}^q
\left[\int_{A_{\delta t}}\left(\frac{\mu(D(w,r))}{(1-|w|^2)^{n+1+\alpha+kp}}
\right)^{\frac p{p-q}}\,dv_{\alpha+kp}\right]^{1-\frac qp}.$$
If the function
$$\widehat\mu_{r,\gamma}(z)=\frac{\mu(D(z,r))}{(1-|z|^2)^{n+1+\alpha+kp}}$$
is in $L^{p/(p-q)}(\bn,dv_{\alpha+kp})$, then
for any given $\epsilon>0$ there is a $t\in (0,1)$ such that
$$\int_{A_{\delta t}}\left(\frac{\mu(D(w,r))}{(1-|w|^2)^{n+1+\alpha+kp}}
\right)^{p/(p-q)}\,dv_{\alpha+kp}(w)<\epsilon^{p/(p-q)}.$$
Thus for such $t$,
$$\int_{A_t}|R^k f_j(z)|^q\,d\mu(z)\le C\epsilon.$$
Since $\bn\setminus A_t$ is a compact subset of $\bn$
and  $f_j\to0$ uniformly on every compact subset of $\bn$, we have
$$\lim_{j\to\infty}\int_{\bn\setminus A_t}|R^k f_j(z)|^q\,d\mu(z)=0.$$ 
Combining this with an earlier estimate we get
$$\limsup_{j\to\infty}\inb|R^k f_j(z)|^q\,d\mu(z)\le C\epsilon.$$
Since $\epsilon$ is arbitrary, we must have
$$\lim_{j\to\infty}\inb|R^kf_j(z)|^q\,d\mu(z)=0.$$
This shows that (c) implies (b) as well.

The proof of that (b) implies (a) is standard. In fact, if (a) is not true, then there is a 
sequence $\{f_j\}$ in $A^p_{\alpha}$ such that $\|f_j\|_{p,\alpha}\le 1$ and
\begin{equation}
\lim_{j\to\infty}\inb|R^kf_j(w)|^q\,d\mu(w)=\infty. 
\label{eq32}
\end{equation}
Since  $\|f_j\|_{p,\alpha}\le 1$,  $\{f_j\}$ is uniformly bounded on compact subsets of 
$\bn$. By Montel's Theorem, there is a subsequence of $\{f_j\}$, which we still 
denote by $\{f_j\}$, that converges uniformly on compact subsets of $\bn$
to a holomorphic function $f$ in $\bn$. It follows from Fatou's lemma that 
$f\in A^p_{\alpha}$ with $\|f\|_{p,\alpha}\le 1$. In particular,
$$\|f_j-f\|_{p,\alpha}\le \max(2,2^{1/p})$$
and $f_j-f\to0$ uniformly on compact subsets of $\bn$. If condition (b) holds, then
$$\lim_{j\to\infty}\inb|R^kf_j(w)-R^kf(w)|^q\,d\mu(w)=0,$$
which contradicts (\ref{eq32}). This shows that (b) implies (a).

To prove that (a) implies (c), we follow the proof of Theorem 1 in Luecking \cite{luecking2}. 
Let $\{a_j\}$ be the sequence of points in $\bn$ from Theorem 2.30 in \cite{zhu6}.
Let $b$ be a real number such that
$$b>n\max\left(1,\frac1p\right)+\frac{1+\alpha}{p}.$$ 
Let 
$$g_j(z)=\frac{(1-|a_j|^2)^{b-(n+1+\alpha)/p}}{(1-\langle z,a_j\rangle)^{b+k}}
=\frac{(1-|a_j|^2)^{(b+k)-(n+1+\alpha+kp)/p}}{(1-\langle z,a_j\rangle)^{b+k}}.$$
Let $\{c_j\}\in l^p$. Then by Theorem 2.30 of \cite{zhu6}, we have
$$\sum_{j=1}^{\infty}c_jg_j(z)\in A^p_{\alpha+kp}.$$
Let
$$h_j(z)=R^{-k}(g_j(z)-g_j(0))$$ 
and
$$f(z)=\sum_{j=1}^{\infty} c_jh_j(z).$$
Then 
$$R^kf(z)=\sum_{j=1}^{\infty} c_j R^k h_j(z)=\sum_{j=1}^{\infty} c_j (g_j(z)-g_j(0)).$$
It is clear that $R^kf\in A^p_{\alpha+kp}$, and so $f\in A^p_{\alpha}$. Moreover,
$$\|f\|_{p,\alpha}^p\le C\sum_{j=1}^\infty |c_j|^p,$$
where $C$ is a positive constant independent of $f$. If condition (a) holds, then
$$\inb|R^kf(z)|^q\,d\mu(z)\le C\|f\|_{p,\alpha}^q\le C\left(\sum_{j=1}^\infty |c_j|^p
\right)^{q/p}.$$
Therefore,
\begin{eqnarray*}
&&\inb\left|\sum_{j=1}^{\infty} c_jg_j\right|^q\,d\mu\\
&\le& 2^q\left[\inb\left|\sum_{j=1}^{\infty} c_jg_j-\sum_{j=1}^{\infty} c_jg_j(0)
\right|^q\,d\mu+\inb\left|\sum_{j=1}^{\infty} c_jg_j(0)\right|^q\,d\mu\right]\\
&\le& 2^q\inb|R^kf(z)|^q\,d\mu(z)+2^q\mu(\bn)\left(\sum_{j=1}^\infty |c_j|^p\right)^{q/p}\\
&\le& C\left(\sum_{j=1}^\infty |c_j|^p\right)^{q/p}
\end{eqnarray*}
Let $r_j(t)$ be a sequence of Rademacher functions (see page 336 of
Luecking \cite{luecking2}).
If we replace $c_j$ by $r_j(t)c_j$, the above inequality is still true, so
$$\inb\left|\sum_{j=1}^{\infty} r_j(t)c_jg_j(z)\right|^q\,d\mu(z)
\le C\left(\sum_{j=1}^\infty |c_j|^p\right)^{q/p}.$$
Integrating with respect to $t$ from $0$ to $1$, applying Fubini's theorem, and
invoking Khinchine's inequality (see Luecking \cite{luecking2}), we obtain
$$A_q\inb\left(\sum_{j=1}^{\infty}|c_j|^2|g_j(z)|^2\right)^{q/2}\,d\mu(z)
\le C\left(\sum_{j=1}^\infty |c_j|^p\right)^{q/p},$$
where $A_p$ is the constant that appears in Khinchine's inequality.
The rest of the proof is exactly the same as the one in Luecking \cite{luecking2}.

The condition in (d) first appeared in Choe-Koo-Yi \cite{cky}, where it was used for the 
embedding of harmonic Bergman spaces into $L^q(d\mu)$. Our proof of the equivalence 
of (c) and (d) follows the method in \cite{cky}. In fact, if $\widehat\mu_{r,\gamma}$ is in
$L^{p/(p-q)}(\bn,dv_\gamma)$, then an application of Proposition~\ref{8} shows that
the function $B_{s,\gamma}(\widehat\mu_{r,\gamma})$ is also in 
$L^{p/(p-q)}(\bn,dv_\gamma)$. By Lemma~\ref{add2}, we must have
$B_{s,\gamma}(\mu)\in L^{p/(p-q)}(\bn,dv_\gamma)$. This proves that (c) implies (d).
That (d) implies (c) is a direct consequence of Lemma~\ref{add3}. The proof of the
theorem is now complete.
\end{proof}

\section{Coefficient Multipliers}

Recall from Theorem~\ref{12} that for $t=(\alpha-\beta)/p$, the operator 
$R_{s,t}$ maps $A^p_\alpha$ boundedly onto $A^p_\beta$. In terms of
Taylor coefficients, we have
$$R_{s,t}\left(\sum_ma_mz^m\right)=\sum_mc_ma_mz^m,$$
where
$$c_m=\frac{\Gamma(n+1+s+t)\Gamma(n+1+|m|+s)}{\Gamma(n+1+s)
\Gamma(n+1+|m|+s+t)}.$$
Therefore, the operator $R_{s,t}$ is just a coefficient multiplier on
holomorphic functions in $\bn$. When $\alpha$ and $\beta$ are real,
an application of Stirling's formula shows that
$$c_m\sim\frac1{|m|^t}$$
as $|m|\to\infty$. We are going to show that this result still holds if
we replace the multiplier sequence $\{c_m\}$ above by the more
explicit multiplier sequence $\{|m|^{(\beta-\alpha)/p}\}$. A similar
result will be proved for the generalized Lipschitz spaces $\Lambda_\alpha$.

We introduce two methods, one based on complex interpolation and 
the other based on atomic decomposition.

\begin{lemma}
Suppose $t$ is complex and $k$ is a postive integer large enough so
that $k+\re t>0$. There exists a constant $c$ such that
$$\int_0^1R^kf(rz)\left(\log\frac1r\right)^{t+k-1}\frac{dr}r=
c\sum_{|m|>0}|m|^{-t}a_mz^m$$
for all holomorphic
$$f(z)=\sum_ma_mz^m$$
in $\bn$.
\label{52}
\end{lemma}

\begin{proof}
Fix $z\in\bn$. We want to evaluate the integral
$$I(f,z)=\int_0^1R^kf(rz)\left(\log\frac1r\right)^{t+k-1}\frac{dr}r$$
in terms of the Taylor expansion of $f$. If $f(z)=\sum_ma_mz^m$, then
$$R^kf(z)=\sum_{|m|>0}|m|^ka_mz^m,$$
so
$$I(f,z)=\sum_{|m|>0}|m|^ka_mz^m\int_0^1r^{|m|-1}\left(\log\frac1r
\right)^{t+k-1}\,dr.$$
Making the change of variables $r=e^{-s}$, we obtain
$$I(f,z)=\sum_{|m|>0}|m|^ka_mz^m\int_0^\infty e^{-|m|s}s^{t+k-1}\,ds.$$
Let $u=|m|s$. Then
$$I(f,z)=c\sum_{|m|>0}|m|^{-t}a_mz^m,$$
where
$$c=\int_0^\infty e^{-u}u^{t+k-1}\,du.$$
\end{proof}

Given any real $\alpha$ and $\beta$, we are going to fix a sufficiently large 
positive integer $k$ and consider operators on $H(\bn)$ of the following form:
$$T_\zeta f(z)=\int_0^1R^kf(rz)\left(\log\frac1r\right)^{(\alpha-\beta)(1-\zeta)
+k-1}\frac{dr}r,$$
where $0\le\re\zeta\le1$.

\begin{lemma}
If $\re\zeta=0$, the operator $T_\zeta$ maps $A^1_\alpha$ boundedly 
into $A^1_\beta$.
\label{53}
\end{lemma}

\begin{proof}
Let $N$ be a sufficiently large positive integer. We have
$$R^NT_\zeta f(z)=\int_0^1R^{N+k}f(rz)\left(\log\frac1r\right)^{(\alpha-\beta)
(1-\zeta)+k-1}\frac{dr}r.$$
If $\re\zeta=0$, it follows from Fubini's theorem that the integral
$$I=\inb|R^NT_\zeta f(z)|(1-|z|^2)^{N+\beta}\,dv(z)$$
does not exceed 
$$\int_0^1\left(\log\frac1r\right)^{\alpha-\beta+k-1}\frac{dr}r\inb
|R^{N+k}f(rz)|(1-|z|^2)^{N+\beta}\,dv(z).$$
Let $w=rz$ in the inner integral. Then $I$ does not exceed the integral
$$\int_0^1\left(\log\frac1r\right)^{\alpha-\beta+k-1}\frac{dr}{r^{2n+1}}
\int_{|w|<r}|R^{N+k}f(w)|\left(1-\frac{|w|^2}{r^2}\right)^{N+\beta}\,dv(w).$$
Since
$$1-\frac{|w|^2}{r^2}\le1-|w|^2$$
for all $|w|<r$, we have
$$I\le\int_0^1\left(\log\frac1r\right)^{\alpha-\beta+k-1}\frac{dr}{r^{2n+1}}
\int_{|w|<r}|R^{N+k}f(w)|(1-|w|^2)^{N+\beta}\,dv(w).$$
We interchange the order of integration and obtain
$$I\le\inb|R^{N+k}f(w)|(1-|w|^2)^{N+\beta}\,dv(w)\int_{|w|}^1\left(
\log\frac1r\right)^{\alpha-\beta+k-1}\frac{dr}{r^{2n+1}}.$$
It is easy to see that there exists a constant $C>0$ such that
$$\int_{|w|}^1\left(\log\frac1r\right)^{\alpha-\beta+k-1}\frac{dr}{r^{2n+1}}
\le C\frac{(1-|w|^2)^{\alpha-\beta+k}}{|w|^{2n}}$$
for all $w\in\bn$, so
$$I\le C\inb|R^{N+k}f(w)|(1-|w|^2)^{N+k+\alpha}\frac{dv(w)}{|w|^{2n}}.$$
Since $|R^{N+k}f(w)|\le C|w|$ near the origin and
$$\inb\frac{dv(w)}{|w|^{2n-1}}<\infty$$
by polar coordinates, we can find another constant $C'>0$, independent 
of $f$, such that
$$I\le C\inb|R^{N+k}f(w)|(1-|w|^2)^{N+k}\,dv_\alpha(w).$$
This completes the proof of the lemma.
\end{proof}

\begin{lemma}
If $\re\zeta=1$, the operator $T_\zeta$ is bounded on the Bloch space $\bloch$.
\label{54}
\end{lemma}

\begin{proof}
We have
$$RT_\zeta f(z)=\int_0^1R^{k+1}f(rz)\left(\log\frac1r\right)^{(\alpha-\beta)(1
-\zeta)+k-1}\frac{dr}r.$$
If $\re\zeta=1$ and $f\in\bloch$, then
\begin{eqnarray*}
|RT_\zeta f(z)|&\le&\int_0^1|R^{k+1}f(rz)|\left(\log\frac1r\right)^{k-1}\frac{dr}r\\
 &\le& C\int_0^1(1-r^2|z|^2)^{-(k+1)}\left(\log\frac1r\right)^{k-1}\,dr\\
 &\le& C'(1-|z|^2)^{-1},
\end{eqnarray*}
where $C$ and $C'$ are positive constants independent of $z$. This shows
that $T_\zeta f$ is in the Bloch space.
\end{proof}

\begin{lemma}
Suppose $s>-1$, $t$ is a positive integer, and
$$I(z)=\int_0^1\frac{(1-x)^s\,dx}{(1-xz)^{s+t+1}},\qquad z\in\D.$$
There exists a polynomial $p(z)$ such that
$$I(z)=\frac{p(z)}{(1-z)^t},\qquad z\in\D.$$
\label{55}
\end{lemma}

\begin{proof}
We compute the integral $I(z)$ with the help of Taylor expansion.
$$I(z)=\sum_{k=0}^\infty\frac{\Gamma(k+s+t+1)}{k!\,\Gamma(s+t+1)}z^k
\int_0^1x^k(1-x)^s\,dx.$$
Since
$$\int_0^1x^k(1-x)^s\,dx=\frac{\Gamma(k+1)\Gamma(s+1)}{\Gamma(k+s+2)},$$
we have
\begin{eqnarray*}
I(z)&=&\frac1{s+1}\sum_{k=0}^\infty\frac{\Gamma(s+2)\Gamma(k+s+t+1)}
{\Gamma(s+t+1)\Gamma(k+s+2)}z^k\\
&=&\frac1{s+1}R^{s,t-1}\sum_{k=0}^\infty z^k\\
&=&\frac1{s+1}R^{s,t-1}\frac1{1-z}.
\end{eqnarray*}
Since $t$ is a positive integer, the operator $R^{s,t-1}$ is a linear
differential operator of order $t-1$ on $H(\D)$ with polynomial coefficients 
(see Proposition~\ref{4}). It follows that there exists a polynomial $p(z)$ 
such that
$$I(z)=\frac{p(z)}{(1-z)^t}.$$
This completes the proof of the lemma.
\end{proof}

We can now prove the first main result of the section.

\begin{thm}
Suppose $\alpha$ is real, $\beta$ is real, and $p>0$. Then
the operator $T$ defined on $H(\bn)$ by
$$Tf(z)=f(0)+\sum_{|m|>0}|m|^{(\beta-\alpha)/p}a_mz^m,
\quad f(z)=\sum_ma_mz^m,$$
maps $A^p_\alpha$ boundedly onto $A^p_\beta$.
\label{56}
\end{thm}

\begin{proof}
By switching the roles of $\alpha$ and $\beta$, it is enough for us to show
that the operator $T$ maps $A^p_\alpha$ into $A^p_\beta$. When $p=1$,
the desired result follows from Lemmas~\ref{52} and \ref{53}. 

First suppose that $1<p<\infty$ with $1/p+1/q=1$.  Let $\theta=1/q$. Then
$$\frac1p=\frac{1-\theta}1+\frac\theta\infty.$$
Because the dual space of $A^p_\beta$ can be identified with $A^q_\beta$ 
under the integral pairing
$$\langle f,g\rangle=\inb(1-|z|^2)^NR^{s,N}f(z)\,\overline{(1-|z|^2)^NR^{s,N}g(z)}
\,dv_\beta(z),$$
where $N$ is a sufficiently large positive number, it suffices for us to show 
that there exists a constant $C>0$, independent 
of $f$ and $g$, such that
\begin{equation}
|\langle Tf,g\rangle|\le C\|f\|_{p,\alpha}\|g\|_{q,\beta}
\label{eq33}
\end{equation}
for all $f\in A^p_\alpha$ and $g\in A^q_\beta$.

Fix a unit vector $f$ in $A^p_\alpha$ and fix a polynomial $g$ that is a unit vector
in $A^q_\beta$ (recall that the polynomials are dense in $A^q_\beta$).
It follows from the complex interpolation relation (see Theorems~\ref{38})
$\left[A^1_\alpha,\bloch\right]_\theta=A^p_\alpha$
that there exist functions $f_\zeta$, where $\zeta\in\overline S$, such that
\begin{enumerate}
\item[(a)] $f_\theta=f$.
\item[(b)] $\zeta\mapsto f_\zeta$ is a bounded continuous function from
$\overline S$ into $A^1_\alpha+\bloch$ whose restriction to $S$ is analytic.
\item[(c)] $\zeta\mapsto f_\zeta$ is a bounded continuous function from
$L(S)$ into $A^1_\alpha$ with $\|f_\zeta\|_{1,\alpha}\le C$.
\item[(d)] For $\zeta\mapsto f_\zeta$ is a bounded continuous function
from $R(S)$ into $\bloch$ with $\|f_\zeta\|_{\bloch}\le C$.
\end{enumerate}
Here $C$ is a positive constant independent of $f$.

Consider the function
$$F(\zeta)=\inb (1-|z|^2)^NR^{s,N}T_\zeta f_\zeta(z)g_\zeta(z)\,dv_\beta(z),$$
where $\zeta\in\overline S$ and
$$g_\zeta(z)=\frac{\overline{R^{s,N}g(z)}}{|R^{s,N}g(z)|}\left[(1-|z|^2)^N
|R^{s,N}g(z)|\right]^{q\zeta}.$$
Because $g$ is a polynomial, the function $F$ is bounded and continuous on
$\overline S$ and its restriction to $S$ is analytic. When $\zeta=\theta$, it
follows from Lemma~\ref{52} that $F(\theta)=\langle Tf,g\rangle$. 

When $\re\zeta=0$, it follows from Lemma~\ref{53} that $T_\zeta$ maps
$A^1_\alpha$ boundedly into $A^1_\beta$, so there exists a positive
constant $C_0$ such that
$$\|T_\zeta f_\zeta\|_{1,\beta}\le C_0\|f_\zeta\|_{1,\alpha}\le C_0C$$
for all $\re\zeta=0$. Thus there exists a constant $M_0>0$ (independent of
$f$, $g$, and $\zeta$) such that
$$|F(\zeta)|\le\inb(1-|z|^2)^N|R^{s,N}T_\zeta f_\zeta(z)|\,dv_\beta(z)\le M_0$$
for all $\re\zeta=0$.

When $\re\zeta=1$, it follows from Lemma~\ref{54} that $T_\zeta$ is bounded
on the Bloch space, so there exists a positive constant $C_1$ such that
$$\|T_\zeta f_\zeta\|_{\bloch}\le C_1\|f_\zeta\|_{\bloch}\le C_1C$$
for all $\re\zeta=1$. We can then find a positive constant $M_1$ (independent
of $f$, $g$, and $\zeta$) such that
$$|F(\zeta)|\le C_1\inb(1-|z|^2)^{Nq}|R^{s,N}g(z)|^q\,dv_\beta(z)\le M_1$$
for all $\re\zeta=1$.

It follows from Hadamard's three lines theorem that 
$$|F(\theta)|\le M_0^{1-\theta}M_1^\theta.$$
Since $M_0$ and $M_1$ are independent of $f$ and $g$, this yields the 
estimate (\ref{eq33}) and proves the theorem for $1<p<\infty$.

Next assume that $0<p\le1$. By Theorem~\ref{32}, there exists a positive number $b$
(we can choose $b$ to be as large as we want) and a sequence $\{a_k\}$ in $\bn$ 
such that every function $f\in\apa$ can be written as
$$f(z)=\sum_{k=1}^\infty c_kf_k(z),$$
with
$$\sum_{k=1}^\infty|c_k|^p\le C\|f\|^p_{p,\alpha},$$
where $C$ is a positive constant independent of $f$ and
$$f_k(z)=\frac{(1-|a_k|^2)^{b-(n+1+\alpha)/p}}{(1-\langle z,a_k\rangle)^b}.$$
By first considering finite sums and then taking a limit, we may assume that
$$Tf=\sum_{k=1}^\infty c_kTf_k.$$
Since $0<p\le1$, we must have
$$\|Tf\|^p_{p,\beta}\le\sum_{k=1}^\infty|c_k|^p\|Tf_k\|^p_{p,\beta}.$$
Since the sequence $\{f_k\}$ is bounded in $\apa$, the proof of the theorem will be
complete if we can show that there exists a constant $C>0$ such that
$$\|Tf\|_{p,\beta}\le C\|f\|_{p,\alpha}$$
for functions of the form
\begin{equation}
f(z)=\frac1{(1-\langle z,a\rangle)^b},\qquad a\in\bn.
\label{eq34}
\end{equation}

We fix a sufficiently large positive integer $k$ and apply Lemma~\ref{52} to
represent the operator $T$ as
$$Tf(z)=f(0)+c\int_0^1R^kf(rz)\left(\log\frac1r\right)^{k+(\alpha-\beta)/p-1}\frac{dr}r.$$
Write $R^k=R^{k-1}R$ and take the factor $R^{k-1}$ out of the integral sign. Then
$$Tf(z)=f(0)+cR^{k-1}\int_0^1\frac{Rf(rz)}r\left(\log\frac1r\right)^{k+(\alpha-\beta)/p-1}\,dr.$$

We assume that $b$ is chosen so that
$$b-k-\frac{\alpha-\beta}p$$
is a sufficiently large positive integer. It is easy to see that
$$\left(\frac{\log\frac1r}{1-r}\right)^{k+(\alpha-\beta)/p-1}=1+\sum_{j=1}^Lb_j(1-r)^j+
H(r),$$
where $H(r)=O((1-r)^L)$ as $r\to1$. It follows that
$$T=T_0+T_1+\cdots+T_L+T_{L+1},$$
where
$$T_0f(z)=f(0)+cR^{k-1}\int_0^1\frac{Rf(rz)}r(1-r)^{k+(\alpha-\beta)-1}\,dr,$$
and
$$T_jf(z)=cb_jR^{k-1}\int_0^1\frac{Rf(rz)}r(1-r)^{k+j+(\alpha-\beta)/p-1}\,dr,\quad 1\le j\le L,$$
and
$$T_{L+1}f(z)=c\int_0^1R^kf(rz)(1-r)^{k+(\alpha-\beta)/p-1}H(r)\frac{dr}r.$$
It then follows from Lemma~\ref{55} that there exists a constant $C>0$ such that
$$\|T_jf\|_{p,\beta}\le C\|f\|_{p,\alpha}$$
for all $0\le j\le L$ and all functions $f$ given in (\ref{eq34}). The same estimate holds 
for the operator $T_{L+1}$ as well, except this time we do not use Lemma~\ref{55}, but 
use the assumption that $L$ is large enough so that 
$$R^NT_{L+1}f(z)=c\int_0^1R^{N+k}f(rz)(1-r)^{k+(\alpha-\beta)/p-1}H(r)\frac{dr}r$$
is bounded, where $N$ is any nonnegative integer with $pN+\beta>-1$ and
$f$ is given by (\ref{eq34}). This proves the case $0<p\le1$ and completes the proof 
of the theorem.
\end{proof}

As the second main result of this section we establish an isomorphism between 
$\Lambda_\alpha$ and $\Lambda_\beta$ by a simple coefficient multiplier.

\begin{thm}
Suppose $\alpha$ and $\beta$ are real. Then the operator $T$ defined by 
$$f(z)=\sum_ma_mz^m\mapsto Tf(z)=f(0)+\sum_{|m|>0}a_m|m|^{\alpha-\beta}z^m$$
is an invertible operator from $\Lambda_\alpha$ onto $\Lambda_\beta$.
\label{57}
\end{thm}

\begin{proof}
By reversing the role of $\alpha$ and $\beta$, it suffices for us to show that
the operator $T$ maps $\Lambda_\alpha$ boundedly into $\Lambda_\beta$.

Given $f\in\Lambda_\alpha$, we fix a sufficiently large positive integer $k$
and use Lemma~\ref{52} to write
$$Tf(z)=f(0)+c\int_0^1R^kf(rz)\left(\log\frac1r\right)^{k-\alpha+\beta-1}\frac{dr}r.$$
If $N$ is another sufficiently large positive integer, then
$$R^NTf(z)=\int_0^1R^{N+k}f(rz)\left(\log\frac1r\right)^{k-\alpha+\beta-1}
\frac{dr}r.$$
Since $f\in\Lambda_\alpha$, it follows from Lemma~\ref{15} that
$$\sup_{z\in\bn}(1-|z|^2)^{N+k-\alpha}|R^{N+k}f(z)|<\infty.$$
But $R^{N+k}f(0)=0$, we must also have
$$\sup_{z\in\bn}(1-|z|^2)^{N+k-\alpha}\frac{|R^{N+k}f(z)|}{|z|}<\infty.$$
So there exists a constant $C>0$ such that
$$|R^NTf(z)|\le C\int_0^1\left(\log\frac1r\right)^{k-\alpha+\beta-1}\frac{dr}
{(1-r^2|z|^2)^{N+k-\alpha}}.$$
Now the above integral clearly converges near $r=0$. When $r$ is away from
$0$, $\log\frac1r$ is comparable to $1-r^2$. So there exists another constant $C>0$
such that
$$|R^NTf(z)|\le C\int_0^1\frac{(1-r^2)^{k-\alpha+\beta-1}\,dr}{(1-r^2|z|^2)^{N+k-\alpha}}.$$
An elementary estimate then shows that
$$|R^NTf(z)|\le\frac C{(1-|z|^2)^{N-\beta}}$$
for some constant $C>0$ and all $z\in\bn$. This shows that $Tf$ is in $\Lambda_\beta$
and completes the proof of the theorem.
\end{proof}

We mention that, at least in the case $n=1$, the theorem above also follows from
Theorem~\ref{12} and the asymptotic expansion of a ratio of two gamma functions
as given in Tricomi-Erdelyi \cite{te}. In fact, in the one-dimensional case, it is easy 
to see that if
$$f(z)=\sum_{k=0}^\infty a_kz^k$$
is a function in $\Lambda_\alpha$, then the sequence $\{k^\alpha a_k\}$ is
bounded. It is also easy to show that if the sequence $\{k^{\alpha+1}a_k\}$ is
bounded, then the function $f$ is in $\Lambda_\alpha$. This together with the
main result of Tricomi-Erdelyi \cite{te} easily gives Theorem~\ref{56} above.
Coefficients of functions in Bloch and Lipschitz spaces are also studied in
Bennet-Stegenga-Timoney \cite{BST}.

\section{Lacunary Series}

One way to construct concrete examples in certain spaces of analytic functions
is by using lacunary series. 
In this section we characterize lacunary series in weighted Bergman spaces
and Lipschitz spaces. 

We say that an analytic function $f$ on $\bn$ has a lacunary homogeneous expansion
if its homogeneous expansion is of the form 
$$
f(z)=\sum_{k=1}^{\infty} f_{m_k}(z),
$$
where $m_k$ satisfies the condition
$$
\inf_{k}\frac{m_{k+1}}{m_k}=\lambda>1.
$$
If $n=1$, the lacunary homogeneous expansions are just lacunary series in the unit disk.
When $n>1$, we say a lacunary homogeneous expansion is a lacunary series if
every homogeneous polynomial $f_{m_k}$ consists of just one term.

Our first result characterizes a lacunary homogeneous expansion in $\apa$.

\begin{prop} 
Let $0<p<\infty$, $\alpha$ be any real number, and
$$f(z)=\sum_{k=1}^{\infty}f_{m_k}(z)$$ 
be a lacunary homogeneous expansion. Then
$f\in \apa$ if and only if
$$\sum_{k=1}^{\infty}m_k^{-1-\alpha}\|f_{m_k}\|_{H^p}^p<\infty,$$
where 
$$\|f\|_{H^p}=\left(\ins |f(\zeta)|^p\,d\sigma(\zeta)\right)^{1/p}$$
denotes the $H^p$-norm of $f$. 
\label{61}
\end{prop}

\begin{proof} 
By Proposition 3 in Yang-Ouyang \cite{yo}, if  
$$g(z)=\sum_{k=1}^{\infty}g_{m_k}(z)$$
is a lacunary homogeneous expansion, then
$g\in A^p$ if and only if
$$\sum_{k=1}^{\infty}m_k^{-1}\|g_{m_k}\|_{H^p}^p<\infty.$$
Let $f\in \apa$. By Theorem~\ref{10}, if $s$ is a real number such that neither 
$n+s$ nor $n+s+(\alpha/p)$ is a negative integer, then
$f\in\apa$ if and only if $R_{s,\alpha/p}f\in A^p$, 
which, by the above result, is equivalent to 
$$\sum_{k=1}^{\infty}m_k^{-1}\|c_{m_k}f_{m_k}\|_{H^p}^p<\infty,$$
where 
$$c_{m_k}=\frac{\Gamma(n+1+s+\alpha/p)\Gamma(n+1+m_k+s)}
{\Gamma(n+1+s)\Gamma(n+1+m_k+s+\alpha/p)}.$$
It follows from Stirling's formula that
$$c_{m_k}\sim m_k^{-\alpha/p}$$
as $k\to\infty$. Thus the above condition is equivalent to
$$\sum_{k=1}^{\infty}m_k^{-1-\alpha}\|f_{m_k}\|_{H^p}^p<\infty.$$
The proof is complete.
\end{proof}

The next result characterizes a lacunary series in $\apa$.

\begin{prop} 
Let $0<p<\infty$, $\alpha$ be any real number, and
$$f(z)=\sum_{k=1}^{\infty}f_{m_k}(z)$$ 
be a lacunary series, where
$$f_{m_k}(z)=a_kz_1^{m_{k_1}}\cdots z_n^{m_{k_n}}.$$
Then $f\in \apa$ if and ony if
$$\sum_{k=1}^{\infty}\frac{|a_k|^p\prod_{i=1}^n\Gamma(\frac{m_{k_i}p}{2}+1)}
{m_k^{1+\alpha}\Gamma(\frac{m_kp}{2}+n)} <\infty.$$
\label{62}
\end{prop}

\begin{proof} 
Let $\zeta^m=\zeta_1^{m_1}\cdots\zeta_n^{m_n}$ and $|m|=m_1+\cdots m_n$.
An easy modification of the proof of Lemma 1.11 in \cite{zhu6} shows that
$$\|\zeta^m\|_{H^p}^p=\frac{(n-1)!\prod_{i=1}^n\Gamma(\frac{m_{i}p}{2}+1)}
{\Gamma(\frac{|m|p}{2}+n)}. $$
Combining this identity and Proposition~\ref{61}, we get the desired result.
\end{proof}

\begin{prop}  
Let $\alpha$ be any real number, let
$$f(z)=\sum_{k=1}^{\infty}f_{m_k}(z)$$ 
be a lacunary homogeneous expansion, and denote by 
$$\|f_{m_k}\|_{H^{\infty}}=\sup_{\zeta\in\sn}|f_{m_k}(\zeta)|.$$
Then
\begin{enumerate}
\item[(a)]
$f\in \Lambda_{\alpha}$ if and only if
$$\sup_{k\ge 1}m_k^{\alpha}\|f_{m_k}\|_{H^{\infty}}<\infty.$$
\item[(b)] $f\in \Lambda_{\alpha,0}$ if and only if
$$\lim_{k\to\infty}m_k^{\alpha}\|f_{m_k}\|_{H^{\infty}}=0.$$
\end{enumerate}
\label{63}
\end{prop}

\begin{proof} 
The results follow easily from Theorem~\ref{16}, the corresponding result for 
$\Lambda_{\alpha,0}$, and Propositions 2 and 3 in Wulan-Zhu \cite{wz}.
We leave the details to the interested reader.
\end{proof}

\begin{prop} 
Let $\alpha$ be any real number and
$$f(z)=\sum_{k=1}^{\infty}f_{m_k}(z)$$
be a lacunary series, where
$$f_{m_k}(z)=a_kz_1^{m_1}\cdots z_n^{m_{k_n}}.$$
Then 
\begin{enumerate}
\item[(a)]
$f\in \Lambda_{\alpha}$ if and only if
$$\sup_{k\ge 1}m_k^{\alpha}|a_k|\sqrt{\frac{m_{k_1}^{m_{k_1}}
\cdots m_{k_n}^{m_{k_n}}}{m_k^{m_k}}}<\infty.$$
\item[(b)] $f\in \Lambda_{\alpha,0}$ if and only if
$$\lim_{k\to\infty}m_k^{\alpha}|a_k|\sqrt{\frac{m_{k_1}^{m_{k_1}}
\cdots m_{k_n}^{m_{k_n}}}{m_k^{m_k}}}=0.$$
\end{enumerate}
\label{64}
\end{prop}

\begin{proof} 
The results follow directly from Proposition~\ref{63} and Lemma 4 in 
Wulan-Zhu \cite{wz}.
\end{proof}

Several special cases of the main results of this section are known. For example, 
lacunary series in the Bloch space of the unit disk are described in 
Anderson-Clunie-Pommerenke \cite{ACP}, lacunary series in weighted Bergman
spaces $\apa$ of the unit ball, where $\alpha>-1$, are described in Stevi\'c \cite{Stevic},
and lacunary series in Bloch and certain Lipschitz spaces of the unit ball are characterized
in Wulan-Zhu \cite{wz}.

\section{Inclusion Relations}

In this section we study inclusion relations among weighted Bergman spaces and 
Lipschitz spaces. From the definition and Proposition~\ref{64} it is very easy to see 
that if $\alpha>\beta$ then $\Lambda_{\alpha}\subset\Lambda_{\beta}$,
and the inclusion is strict. 

The inclusion relations between weighted Bergman spaces are more complicated
in general. Several embedding theorems have been known before, and our results
here overlap with some of them; see Aleksandrov \cite{Alex}, 
Beatrous-Burbea \cite{BB2}, Graham \cite{Graham}, Luecking \cite{luecking2}, and 
Rochberg \cite{Roch}. We begin with the following simple case.

\begin{prop} 
Let $0<p<\infty$, and let $\alpha$ and $\beta$
be any two real numbers satisfying $\alpha<\beta$. Then
$$\apa\subset A^p_{\beta},$$
and the inclusion is strict.
\label{65}
\end{prop}

\begin{proof} 
The inclusion is obvious. To prove that the inclusion is strict, we only need to
test functions of the form $f_t(z)=(1-z_1)^t$. See Yang-Ouyang \cite{yo} for
a similar argument.
\end{proof}

To better describe the inclusion relations of Bergman spaces, we introduce the
notion of Lipschitz stretch first. More specifically, if $X$ is a space of analytic 
functions, we define the {\it Lipschitz stretch} of $X$ as follows:
$$\Lambda(X)=\inf\{\beta-\alpha: \Lambda_{-\alpha}\subset X\subset
\Lambda_{-\beta}\}.$$
We also call the constants 
$$\alpha_0=\sup\{\alpha: \Lambda_{-\alpha}\subset X\},\quad
\beta_0=\inf\{\beta: X\subset \Lambda_{-\beta}\}$$
the lower and upper bounds of the Lipschitz stretch, respectively. A similar concept 
using Bloch type spaces was introduced in Zhao \cite{zhao} for spaces of analytic
functions in the unit disk.

\begin{thm} 
Let $0<p<\infty$ and let $\alpha$ be any real number. 
Then for any $\gamma<(1+\alpha)/p$ we have
$$\Lambda_{-\gamma}\subset \apa\subset\Lambda_{-(n+1+\alpha)/p}.$$
Both inclusions are strict and  best possible, where``best possible'' means that, for 
each $p$ and $\alpha$, the index $\gamma$ of $\Lambda_{-\gamma}$ on the left-hand 
side cannot be replaced by a larger number, and the index $(n+1+\alpha)/p$ on the 
right-hand side cannot be replaced by a smaller one.
\label{66}
\end{thm}

\begin{proof} 
Suppose $f\in \apa$. Then $R^kf\in A^p_{pk+\alpha}$, where $k$ is a nonnegative 
integer satisfying $pk+\alpha>-1$. By Theorem~\ref{20}, there exists a positive
constant $C$ such that
$$(1-|z|^2)^{k+(n+1+\alpha)/p}|R^kf(z)|\le C$$
for all $z\in\bn$. This means $f\in \Lambda_{-(n+1+\alpha)/p}$, so
$\apa\subset\Lambda_{-(n+1+\alpha)/p}$.

Next suppose $\gamma<(1+\alpha)/p$ and $f\in \Lambda_{-\gamma}$.
Let $k$ be a nonnegative integer such that $k+\gamma>0$. Then $kp+\alpha>-1$ and 
$\alpha-p\gamma>-1$, so
\begin{eqnarray*}
&&\int_{\bn}(1-|z|^2)^{pk}|R^kf(z)|^p\,dv_{\alpha}(z)\\
&\le&\sup_{z\in\bn}\left((1-|z|^2)^{k+\gamma}|R^kf(z)|\right)^p
\int_{\bn}(1-|z|^2)^{-p\gamma}\,dv_{\alpha}(z)\\
&\le& C\sup_{z\in\bn}\left((1-|z|^2)^{k+\gamma}|R^kf(z)|\right)^p.
\end{eqnarray*}
Thus $\Lambda_{-\gamma}\subset\apa$.

We only give a sketch of the rest of the proof since it is similar to the argument
used in Yang-Ouyang \cite{yo}. For $t>0$ let $k$ be a nonnegative integer such that 
$k+\gamma>0$. Since the radial derivative is an invertible operator on the space
of holomorphic functions in $\bn$ that vanish at the origin, we can define a holomorphic
function $f_t$ in $\bn$ by
$$f_t(z)=R^{-k}\left[(1-z_1)^{-t-k}-1\right].$$
Then
$$R^kf_t(z)=(1-z_1)^{-t-k}-1,$$
so for any $z\in\bn$ and $t\le \gamma$,
\begin{eqnarray*}
(1-|z|^2)^{k+\gamma}|R^kf_t(z)|
&\le& (1-|z|^2)^{k+\gamma}\left(|1-z_1|^{-t-k}+1\right)\\
&\le& C(1-|z|^2)^{\gamma-t}\le C.
\end{eqnarray*}
On the other hand, if $t>\gamma$, then we take $z=(x,0,...,0)$, where $x$ is a 
real number between $0$ and $1$, to obtain
$$(1-|z|^2)^{k+\gamma}|R^kf_t(z)|=(1-x^2)^{k+\gamma}((1-x)^{-t-k}-1)
\ge (1-x)^{\gamma-t}\to\infty$$
as $x\to1$. 
Thus 
\begin{equation}
f_t\in\Lambda_{-\gamma} \quad {\rm if \ and\ only\ if } \quad t\le \gamma.
\label{eq35}
\end{equation}
By a similar computation as used in Yang-Ouyang \cite{yo}, we see that
\begin{equation}
f_t\in\apa \quad {\rm when } \quad t<\frac{n+1+\alpha}{p},
\label{eq36}
\end{equation}
and 
\begin{equation}
f_t\not\in\apa \quad {\rm when } \quad t=\frac{n+1+\alpha}{p}.
\label{eq37}
\end{equation}

For any $\varepsilon>0$ let $t=(n+1+\alpha)/p-\varepsilon/2$. Then 
$$(n+1+\alpha)/p-\varepsilon<t<(n+1+\alpha)/p.$$ 
By (\ref{eq36}) and (\ref{eq35}) we hav
$$f_t\in\apa \quad {\rm but}\quad f_t\not\in \Lambda_{-((n+1+\alpha)/p-\varepsilon)}.$$
This shows that the inclusion $\apa\subset\Lambda_{(n+1+\alpha)/p}$
is the best possible. At the same time it also shows that
the inclusion $\Lambda_{-\gamma}\subset\apa$ is strict, since
$$\Lambda_{-\gamma}\subset\Lambda_{-((n+1+\alpha)/p-\varepsilon)}$$
as $\varepsilon\le n/p$.

Let $t=(n+1+\alpha)/p$. Then by  (\ref{eq37}) and (\ref{eq35}),
$f_t\not\in\apa$ but $f_t\in \Lambda_{-(n+1+\alpha)/p}$, so
the inclusion $\apa\subset\Lambda_{-(n+1+\alpha)/p}$ is strict.

To show that the left inclusion is the best possible, we let
$$f_{p,\alpha}(z)=\sum_{k=1}^{\infty}f_{m_k}(z)
=\sum_{k=1}^{\infty}2^{k(1+\alpha)/p}W_{2^k}(z),$$
where $\{W_{2^k}\}$ is a sequence of polynomials with Hadamard gaps as in 
Theorem 1.2 of Ryll-Wojtaszczyk \cite{rw} and Corollary 1 of Ullrich \cite{ullrich} with 
the following properties: 
$$\|W_{2^k}\|_{H^{\infty}}=1,\qquad \|W_{2^k}\|_{H^p}>C(n,p),$$
where $C(n,p)$ is a constant depending only on $n$ and $p$.

From Proposition~\ref{63} and Proposition~\ref{61} we easily deduce that
$f_{p,\alpha}\in\Lambda_{-(1+\alpha)/p}$ but $f_{p,\alpha}\not\in\apa$.
Thus the inclusion $\Lambda_{-\gamma}\subset \apa$
is best possible. The proof is complete.
\end{proof}
 
As a direct consequence we obtain the Lipschitz stretch of $\apa$.

\begin{cor} 
Let $0<p<\infty$ and let $\alpha$ be any real number. Then
$$\Lambda(\apa)=\frac 1p$$
with lower bound $(1+\alpha)/p$ and upper bound $(n+1+\alpha)/p$. 
\label{67}
\end{cor}

\begin{cor} 
All weighted Bergman spaces are different, that is, $A^p_\alpha\not=A^q_\beta$
whenever $(p,\alpha)\not=(q,\beta)$.
\label{68}
\end{cor}

\begin{proof} 
If $p=q$ but $\alpha\neq \beta$, then by Proposition~\ref{65}, $\apa$ and $A^q_{\beta}$ 
are different. If $p\neq q$, then Corollary~\ref{67} tells us that $\Lambda(\apa)=1/p$, 
while $\Lambda(A^q_{\beta})=1/q$. Thus $\apa$ and $A^q_{\beta}$ have different
Lipschitz stretchs, so they must be different.
\end{proof}

The following two theorems completely describe the inclusion relations between two 
weighted Bergman spaces. Beatrous-Burbea \cite{BB2} contains several related results.

\begin{thm} 
Let $0<p\le q<\infty$. Then $\apa\subset A^q_{\beta}$ if and only if
$$\frac{n+1+\alpha}{p}\le \frac{n+1+\beta}{q},$$
and in this case the inclusion is strict.
\label{69}
\end{thm}

\begin{proof} 
Let $0<p\le q<\infty$ and  $f\in \apa$. Let $k$ be a nonnegative integer such that 
$pk+\alpha>-1$ and $kq+\beta>-1$. It follows from the closed graph theorem that
the inclusion $\apa\subset A^q_{\beta}$ is equivalent to
\begin{equation}
\int_{\bn}|R^kf(z)|^q\,dv_{kq+\beta}(z)\le C\|f\|_{p,\alpha}^q.
\label{eq38}
\end{equation}
Let $s>0$ be a real number which is sufficiently large.
By Theorem~\ref{50}, the inequality in (\ref{eq38}) is equivalent to
$$\sup_{z\in\bn}\int_{\bn}\frac{(1-|z|^2)^s}{|1-\langle z,w\rangle|^{s+(n+1+\alpha+
kp)q/p}}\,dv_{kq+\beta}(w)<\infty,$$ 
or
$$\sup_{z\in\bn}(1-|z|^2)^s\int_{\bn}\frac{(1-|w|^2)^{kq+\beta}}
{|1-\langle z,w\rangle|^{s+kq+(n+1+\alpha)q/p}}\,dv(w)<\infty.$$
By Proposition~\ref{7}, the inequality above holds if and only if
$$c=s+kq+(n+1+\alpha)q/p-(n+1)-(kq+\beta)\le s,$$
which is easily seen to be equivalent to
$$\frac{n+1+\alpha}{p}\le \frac{n+1+\beta}{q}.$$
In view of Corollary~\ref{68} the proof is now complete. 
\end{proof}

\begin{thm} 
Let $0<q<p<\infty$. Then $\apa\subset A^q_{\beta}$ if and only if
$$\frac{1+\alpha}{p}< \frac{1+\beta}{q},$$
and in this case the inclusion is strict.
\label{70}
\end{thm}

\begin{proof} 
Let $0<q<p<\infty$ and $f\in \apa$. Let $k$ be a nonnegative integer such that 
$pk+\alpha>-1$ and $kq+\beta>-1$. Once again, the closed graph theorem tells
us that the inclusion $\apa\subset A^q_{\beta}$ is equivalent to
\begin{equation}
\int_{\bn}|R^kf(z)|^q\,dv_{kq+\beta}(z)\le C\|f\|_{p,\alpha}^q.
\label{eq39}
\end{equation}
Let $s>0$ be a real number which is sufficiently large.
By Theorem~\ref{51}, the inequality in (\ref{eq39}) is equivalent to
\begin{equation}
B_{s,kp+\alpha}(v_{kq+\beta})\in L^{p/(p-q)}(\bn, dv_{kp+\alpha}).
\label{eq40}
\end{equation}
If $s$ is large enough, then by Proposition~\ref{7},
\begin{eqnarray*}
B_{s,kp+\alpha}(v_{kq+\beta})(z)
&=&\int_{\bn}\frac{(1-|z|^2)^s(1-|w|^2)^{kq+\beta}}
{|1-\langle w,z\rangle|^{n+1+s+kp+\alpha}}\,dv(w)\\
&\sim&
(1-|z|^2)^{-k(p-q)-(\alpha-\beta)},
\end{eqnarray*}
as $|z|$ approaches $1$. Thus (\ref{eq40}) is equivalent to
$$\int_{\bn}(1-|z|^2)^{-(k(p-q)+(\alpha-\beta))\frac{p}{p-q}+kp+\alpha}\,dv(z)<\infty,$$
which is equivalent to
$$\frac{1+\alpha}{p}<\frac{1+\beta}{q}.$$
This along with Corollary~\ref{68} finishes the proof.
\end{proof}

\section{Special Cases}

In this section we point out to the reader the various special cases of
the spaces $A^p_\alpha$ and $\Lambda_\alpha$.

As was mentioned in the introduction, when $\alpha>-1$, the spaces $A^p_\alpha$
are traditionally called weighted Bergman spaces. In this case, a holomorphic
function $f$ in $\bn$ belongs to $\apa$ if and only if 
$$\inb|f(z)|^p(1-|z|^2)^\alpha\,dv(z)<\infty.$$
See Zhu \cite{zhu6}.

When $\alpha=-(n+1)$, or $n+\alpha+1=0$, we have mentioned several times
earlier that the space $A^p_\alpha$ is traditionally denoted by $B_p$ and is called a 
diagonal Besov space. Alternatively, a holomorphic function $f$ in $\bn$ belongs to 
the Besov space $B_p$ if and only if
$$\inb\bigl|(1-|z|^2)^kR^kf(z)\bigr|^p\,d\tau(z)<\infty,$$
where $k$ is any positive integer with $pk>n$ and
$$d\tau(z)=\frac{dv(z)}{(1-|z|^2)^{n+1}}$$
is the M\"obius invariant measure on $\bn$. See Zhu \cite{zhu6}.

When $\alpha=-1$ and $p=2$, the space $A^p_\alpha$ coincides with the classical
Hardy space $H^2$. See (\ref{eq20}) and (1.22) of Zhu \cite{zhu6}. Recall that $H^p$ 
consists of holomorphic functions $f$ in $\bn$ such that
$$\sup_{0<r<1}\ins|f(r\zeta)|^p\,d\sigma(\zeta)<\infty,$$
where $d\sigma$ is the normalized surface measure on the unit sphere $\sn$.

When $\alpha=-n$ and $p=2$, the space $A^p_\alpha$ is the so-called Arveson
space, which is usually defined as the Hilbert space of holomorphic functions
in $\bn$ whose reproducing kernel is given by
$$K(z,w)=\frac1{1-\langle z,w\rangle};$$
see Theorem~\ref{41}. This space has attracted much attention lately in the study 
of multi-variable operator theory. We mention Arveson's influential paper \cite{arveson} 
and the recent monograph \cite{cg} by Chen and Guo.

When $0<\alpha<1$, the space $\Lambda_\alpha$ is the classical Lipschitz
space of holomorphic functions $f$ in $\bn$ satisfying the condition
$$\sup\left\{\frac{|f(z)-f(w)|}{|z-w|^\alpha}:z,w\in\bn,z\not=w\right\}<\infty.$$
See Section 6.4 of Rudin \cite{rudin}. The space $\Lambda_1$ is also called the 
Zygmund class, especially in the case when $n=1$.

When $\alpha=0$, the space $\Lambda_\alpha$ is just the classical Bloch space,
consisting of functions $f\in H(\bn)$ such that
$$\sup_{z\in\bn}(1-|z|^2)|Rf(z)|<\infty.$$

When $\alpha<0$, the spaces $\Lambda_\alpha$ have appeared in the literature
under the name of growth spaces. In this case, a holomorphic function $f$ in $\bn$
belongs to $\Lambda_\alpha$ if and only if
$$\sup_{z\in\bn}(1-|z|^2)^{|\alpha|}|f(z)|<\infty.$$

The term ``Bloch type spaces" or $\alpha$-Bloch spaces can also be found in
recent literature. More specifically, for any $\alpha>0$ the $\alpha$-Bloch space is 
denoted by $\bloch_\alpha$ and consists of holomorphic functions $f$ in $\bn$
such that
$$\sup_{z\in\bn}(1-|z|^2)^\alpha|Rf(z)|<\infty.$$
It is then clear that the $\alpha$-Bloch space $\bloch_\alpha$ is the same as
our generalized Lipschitz space $\Lambda_{1-\alpha}$. See Zhu \cite{zhu6}.

\section{Further Remarks}

Unless $p=2$, the space $A^p_{-1}$ is not the same as the Hardy
space $H^p$, although in many situations it is useful to think of $H^p$
as the limit of $A^p_\alpha$ as $\alpha$ approaches $-1$. One particular
problem here is to identify the complex interpolation space between
$H^p$ and $A^p_\alpha$ when $\alpha>-1$ and $p\ge1$. It is also interesting
to ask for the complex interpolation space between $H^p$ and $\Lambda_\alpha$.

The spaces $A^p_\alpha$ when $\alpha$ is a negative integer appear 
to be very special. It would be interesting to see some ``singular properties'' of
these spaces.

Finally, we conjecture that Theorem~\ref{46} remains valid when $1<p<\infty$. This 
is probably very difficult, since an affirmative answer would characterize Carleson 
measures for the Hardy space $H^2$ as a special case, and it is well known that 
the characterization of Carleson measures for Hardy spaces is very technical.

\end{document}